\documentclass[preprint]{siamart171218}

\usepackage[utf8]{inputenc}
\usepackage{amssymb}
\usepackage{amsmath}
\usepackage{algpseudocode,algorithm,algorithmicx}

\usepackage{bm}

\usepackage{hyperref}
\usepackage{afterpage}

\usepackage[section]{placeins}

\makeatletter
\AtBeginDocument{%
  \expandafter\renewcommand\expandafter\subsection\expandafter{%
    \expandafter\@fb@secFB\subsection
  }%
}
\makeatother

\usepackage{multirow}
\usepackage{color}

\newtheorem{remark}{Remark}

\usepackage{geometry}
\usepackage{float}
\usepackage{subcaption}
\usepackage{mathtools}

\renewcommand{\chi}{\mathcal{X}}

\newcommand{\He}{\mathcal{H}}

\newcommand{\gG}{{\bm G}}

\newcommand{\gM}{{\bm M}}

\newcommand{\gV}{{\bm V}}

\newcommand{\fC}{\mathcal{C}}
\newcommand{\fK}{\mathcal{K}}
\newcommand{\fV}{\mathcal{V}}
\newcommand{\fF}{\mathcal{F}}

\newcommand{\bigO}{\mathcal{O}}
\newcommand{\incr}{ \Delta}

\newcommand{\Lie}{\mathcal{L}}
\newcommand{\indi}{\mathbf{1}}

\newcommand{\vu}{{\bm u}}
\newcommand{\vv}{{\bm v}}
\newcommand{\vw}{{\bm w}}
\newcommand{\vtx}{{\omega}}
\newcommand{\vx}{{\bm x}}
\newcommand{\vy}{{\bm y}}
\newcommand{\vxi}{{\bm \xi}}

\newcommand{\vg}{{\bm \gamma}}
\newcommand{\vs}{{\bm s}}
\newcommand{\vi}{{\bm i}}
\newcommand{\vr}{{\bm r}}
\newcommand{\va}{{\bm a}}
\newcommand{\vb}{{\bm b}}

\newcommand{\vn}{{\bm n}}
\newcommand{\vT}{{\bm T}}
\newcommand{\vR}{{\bm R}}

\newcommand{\vX}{{\bm X}}
\newcommand{\vI}{{\bm I}}
\newcommand{\vtX}{\tilde{\bm X}}
\newcommand{\vtu}{\tilde{\bm u}}
\newcommand{\vtw}{\tilde{\bm w}}
\newcommand{\vhX}{{\bm \chi}}

\newcommand{\vPsi}{{\bm \Psi}}
\newcommand{\vPhi}{{\bm \Phi}}

\newcommand{\vphi}{{\bm \varphi}}
\newcommand{\vpsi}{{\bm \psi}}

\newcommand{\vhE}{{\bm E}}
\newcommand{\vhe}{{\bm \epsilon}}
\newcommand{\vhxi}{{\bm \xi}}
\newcommand{\vheta}{{\bm \eta}}

\renewcommand{\d}{d}

\newcommand{\R}{{\mathbb{R}}}

\newcommand{\curl}{\nabla \times}
\renewcommand{\div}{\nabla \cdot}
\newcommand{\grad}{\nabla}
\newcommand{\Laplace}{\Delta}
\newcommand{\PSolve}{\Laplace^{-1}}

\definecolor{myPurp}{rgb}{0.5, 0.0, 1.0}

\title{A Characteristic Mapping Method for the three-dimensional incompressible Euler equations}

\author{Xi-Yuan Yin\thanks{Department of Mathematics and Statistics, McGill University, Montr\'{e}al, Qu\'{e}bec H3A 0B9, Canada}
\and Kai Schneider\thanks{Institut de Math\'{e}matiques de Marseille, Aix-Marseille Universit\'{e}, CNRS, 13453 Marseille Cedex 13, France}
\and Jean-Christophe Nave\footnotemark[1]  \thanks{Corresponding author. E-mail address: \email{jcnave@math.mcgill.ca}.}}

\begin{document}

\maketitle

\begin{abstract}
We propose an efficient semi-Lagrangian Characteristic Mapping (CM) method for solving the three-dimensional (3D) incompressible Euler equations. This method evolves advected quantities by discretizing the flow map associated with the velocity field. Using the properties of the Lie group of volume preserving diffeomorphisms \textit{SDiff}, long-time deformations are computed from a composition of short-time submaps which can be accurately evolved on coarse grids. This method is a fundamental extension to the CM method for two-dimensional incompressible Euler equations \cite{CME}. We take a geometric approach in the 3D case where the vorticity is not a scalar advected quantity, but can be computed as a differential 2-form through the pullback of the initial condition by the characteristic map. This formulation is based on the Kelvin circulation theorem and gives point-wise a Lagrangian description of the vorticity field. We demonstrate through numerical experiments the validity of the method and show that energy is not dissipated through artificial viscosity and small scales of the solution are preserved. We provide error estimates and numerical convergence tests showing that the method is globally third-order accurate.
\end{abstract}

\section{Introduction}
Turbulence remains one of the oldest and most challenging research problems in both pure and applied science; the high Reynolds number limit of the Navier-Stokes equations is of particular interest and developments in scientific computing are useful in advancing the frontiers of our understanding of fluid dynamics in  a highly nonlinear regime. Efficient and precise numerical schemes for simulating incompressible inviscid fluids in three-dimensional space are an essential component, however
the multiscale nature of turbulent flows along with the computationally demanding high dimensionality requires the development of specialized algorithms. Different directions have been pursued so far. 

Among the Eulerian ansatz using a fixed computational grid, Fourier pseudo-spectral methods are certainly very attractive discretizations \cite{ishihara2009}, requiring nevertheless some viscous or hyperviscous regularization, see e.g. the discussion in \cite{farge2017euler}. However one drawback of Eulerian schemes is the Courant-Friedrichs-Lewy (CFL) condition which constraints the size of time steps in relation to the spatial discretization scales. This implies that the complexity of the simulations with $N$ grid points in each spatial direction is proportional to $N^4$ or even worse \cite{schneider2013cfl}. The progress of high-resolution numerical simulation using Fast Fourier Transforms is hence limited and directly linked to the development of supercomputers following Moore’s law. Furthermore, Eulerian methods are also prone to artificial dissipation and special care needs to be taken in the spatial resolution of the solution and when using spectral dealiasing, a comparative study of spatial discretization methods for the Euler equations can be found in \cite{grafke2008numerical}. Generally, high-resolution computational grids are needed to reduce the effects of dissipation, currently pseudo-spectral simulations with up to $12288^3$ grid points can be performed \cite{ishihara2020second}.

The Lagrangian ansatz, and in particular the semi-Lagrangian approach which combines Lagrangian time integration with Eulerian
grids, does not suffer from a time-step restriction due to the CFL condition, see e.g. Staniforth and C\^ot\'e \cite{staniforth1991semi} (1991) and references therein.
Hence they are well suited for advection-dominated problems. Purely Lagrangian approaches include the vortex blobs methods \cite{chorin1973numerical, beale1982vortex,  hou1990convergence, pelz1997locally, cottet2000vortex, oliver2001vortex}, and the vortex particle and filament methods \cite{rosenhead1932point, hou1990convergence, cottet1991convergence, winckelmans1993contributions, rossinelli2010gpu}. These methods are characterized by a particle-based discretization of the vorticity field; the motion of the fluid is idealized as the transport of a collection of point vortices or compactly supported vortex blobs and the velocity of each particle can be computed from vortices using the Biot-Savart law. These methods are inherently spatially adaptive since the representation of the vorticity field is reduced to a collection of point vortices concentrated where the vorticity is important. Furthermore, they are more effective in avoiding artificial viscous dissipation compared to their Eulerian counterparts. Some drawbacks include the difficulty in the representation and controlled resolution of Eulerian quantities. Methods for transferring Lagrangian quantities to fixed Eulerian grids include the vortex-in-cell methods \cite{christiansen1973numerical, couet1981simulation, cottet2004advances, sbalzarini2006ppm} and the Cauchy-Lagrangian frameworks \cite{frisch}. 

Various Eulerian or semi-Lagrangian methods have been used to provide evidence of singularity for the 3D Euler equations \cite{ashurst1987numerical, grauer1998adaptive, bustamante20083d, hou2018potential} or nonsingular super-exponential growth in the maximum vorticity \cite{pumir1990collapsing, brachet1992numerical, hou2006dynamic}. 
For instance a pseudospectral computation of an axisymmetric solution suggesting finite-time blow-up was performed by Kerr in 1993 \cite{kerr1993evidence}. Some recent computations including new test cases for potential singularities in 3D Euler have been proposed in \cite{moffatt2020towards, yao2020physical}.

%Originality/goal: proof of concept of the CM method, potential to investigate singularities in 3d Euler.
The work in this paper can be seen as a proof of concept for a novel semi-Lagrangian numerical method, which allows a very detailed investigation of singularities in 3D Euler.
To this end we propose a geometric method for the 3D incompressible Euler equations in its vorticity form. The method uses the numerical framework of the Gradient-Augmented Level-Set methods \cite{nave2010gradient} and Reference Map methods \cite{kohno2013new} and consists in a semi-Lagrangian discretization of the backward flow map, called the characteristic map, generated by the velocity field. This is based on a previous work on the 2D incompressible Euler equations using the Characteristic Mapping method \cite{CME}, which we extend and generalize here to the three-dimensional case. Compared to the 2D equations, the 3D Euler equations present several significant challenges. Firstly, the presence of an additional vortex stretching term requires a more geometric formulation of the CM method in order to be seamlessly incorporated in the framework: a direct treatment of the vortex stretching as a source term would not conform to the characteristic structure of the method, nullifying its numerical qualities. Secondly, due to the vortex stretching term, we no longer have conservation of any $L^p$ norms, including the $L^\infty$ norm. Indeed, in the 2D case, the scalar vorticity is an advected quantity with infinitely many Casimir invariants, which can be preserved numerically \cite{bowman2015fully}, and small scale features appear only from rapid growth in vorticity gradient. In the 3D case, the question of finite-time blow-up in the solutions of incompressible Euler equations with smooth initial data is a notoriously difficult open problem in the theory of PDEs and is related to the Clay institute Millenium problem on the Navier-Stokes equations \cite{fefferman2006}. 
From a numerical point of view, the rapid growth in both the magnitude of the vorticity and its gradient further increases the difficulty in providing sufficient spatial resolution of the solution. For the CM method, the dependency of vorticity on the spatial deformations or rate of strain tensor would then involve the Jacobian of the characteristic map in the computation of the vorticity, thus increasing the regularity requirements on the method. Lastly, the higher dimensionality further emphasizes on the computational efficiency of the method and on high order accuracy.

The rest of the paper is organized as follows: in section \ref{sec:mathForm} we recall the CM framework presented in \cite{CME} and generalize it through a geometric reformulation in the context of the 3D incompressible Euler equations. In section \ref{sec:numImpl}, we provide some details on the numerical implementation of the method together with formal error estimates supported by convergence tests. Section \ref{sec:numTests} contains numerical tests and discussions; In particular we simulate the anti-parallel axisymmetric perturbed vortex tubes tests similar to those appearing in \cite{kerr1993evidence} and \cite{hou2006dynamic}. Finally, in section \ref{sec:conclu} we make some concluding remarks and propose future directions of work.

\pagebreak
\section{Mathematical Formulation} \label{sec:mathForm}

We present here the Characteristic Mapping method for the incompressible Euler equations in three-dimensional space. This work is the natural continuation of the framework presented in \cite{kohno2013new, CM, CME} and extends the Characteristic Mapping method for the 2D Euler equations in \cite{CME} to the 3D case. The main challenge in the 3D case is the inclusion of the vortex stretching term in a way that is compatible with the CM method and preserves its arbitrary resolution and low-dissipation properties. For this, we expand on the CM framework by including a more geometric formulation of the problem in terms of differential forms. 

\subsection{Characteristic Mapping Method} \label{sec:CMM}

The Characteristic Mapping (CM) method consists in computing the diffeomorphic transformations of the domain generated by a given velocity field. For a given smooth and divergence-free velocity field $\vu$ on a three-dimensional domain $U \in \R^3$, we consider the family of characteristic curves $\vg(t)$ parametrized by their initial condition $\vg_0$:
\begin{subequations} \label{eqs:charCurves}
\begin{gather} 
\frac{d}{dt} \vg(t) = \vu (\vg(t), t)  , \\
\vg(0) = \vg_0 . 
\end{gather}
\end{subequations}

We define the characteristic map $\vX_{[t_1, t_2]}$ associated with the velocity $\vu$ to be the solution operator for the characteristic curves, that is,
\begin{gather}
\vX_{[t_1, t_2]} (\vg(t_1)) = \vg(t_2)
\end{gather}
for all times $t_1, t_2$ and for all characteristic curves $\vg$. One can check that, keeping $t_0$ fixed, the characteristic map satisfies the following equations
\begin{subequations} \label{eqs:charEqns}
\begin{gather}
\partial_t \vX_{[t_0, t]} = \vu( \vX_{[t_0, t]}, t) ,  \label{eq:CMfwd} \\
(\partial_t + \vu \cdot \grad)  \vX_{[t, t_0]} = 0 .  \label{eq:CMbkw}
\end{gather}
\end{subequations}

The map $\vX_{[t_1, t_2]}$ can be thought of as a transformation the space from time $t_1$ to $t_2$ following the flow, i.e. vortex lines are transported. There is no requirement that $t_1 < t_2$, if $t_1 <t_2$, the map is forward in time and we will call it the forward map, if $t_1 > t_2$, we call it the backward map. It is straightforward to check the following properties of the characteristic maps:
\begin{subequations} \label{eqs:CMGroup}
\begin{gather}
\vX_{[t_1, t_2]} \circ \vX_{[t_0, t_1]} = \vX_{[t_0, t_2]} , \label{eq:groupCompose} \\
\vX_{[t_0, t_1]}^{-1} = \vX_{[t_1, t_0]} , \label{eq:groupInverse} \\
\vX_{[t_0, t_0]} = \vx ,
\end{gather}
\end{subequations}
for arbitrary $t_0, t_1, t_2$. Indeed, for a given divergence-free velocity field $\vu$, the characteristic maps $\vX$ are elements of $\textit{SDiff}(U)$, the Lie group of volume-preserving diffeomorphisms of the domain $U$, with the space of divergence-free vector fields as its Lie algebra. The underlying theory for the characterization of the Euler equations as geodesic flow in the space of volume preserving diffeomorphisms can be found in the works of Arnold \cite{arnold1966geometrie}. For simplicity of notation, we will denote the forward map $\vX_{[0, t]} (\vx)$ as $\vX_F(\vx, t)$ and the backward map $\vX_{[t, 0]} (\vx)$ as $\vX_B(\vx, t)$ when the time-interval of mapping is not emphasized. As a function of $t$, we can formally see $\vX_F$ as the integral curve of the time-dependent velocity field $\vu$ on $\textit{SDiff}(U)$ starting from identity; $\vX_B$ is the corresponding inverse element of the group for each time $t$.

The characteristic maps act as solution operators to the transport equations. Consider the advection equation for a scalar $\phi$ under the velocity field $\vu$:
\begin{gather}
(\partial_t + \vu \cdot \grad) \phi = 0 , \\
\phi(\vx, 0) = \phi_0(\vx) .
\end{gather}
From the method of characteristics, we know that $\frac{d}{dt} \phi(\vg(t), t) = 0$ for any characteristic $\vg$ given by \eqref{eqs:charCurves}. It follows that
\begin{gather}
\phi(\vx, t) = \phi_0 (\vX_B(\vx, t)) .
\end{gather}

In geometric terms, since $\phi$ is Lie advected by $\vu$, i.e. $(\partial_t + \Lie_\vu ) \phi = 0$, we have that $\phi$ is given by the pullback ${\vX_B}^* \phi_0$. This is also called the relabelling symmetry or back-to-label map. Indeed, $\vX_B$ allows us to switch between Lagrangian and Eulerian frames. One can think of $\vX_B$ as identifying the characteristic curve passing through $\vx$ at time $t$ and returning the location of the corresponding particle in $U$ at time $0$ which, by convention, we use as the Lagrangian reference space. Here, $\phi_0$, the initial condition of the scalar, is a 0-form, however this is true for higher degree forms and will be our main tool for solving the Euler equations in vorticity form.

\subsection{The Euler Equations}

We consider the incompressible Euler equations on a three-dimensional domain $U$, for simplicity, we take $U$ to be the periodic cube $\mathbb{T}^3$.
\begin{subequations}
\begin{gather}
\partial_t \vu + (\vu \cdot \grad) \vu = \frac{1}{\rho} \grad p , \\
\partial_t \rho + \div (\rho \vu) = 0  \label{eq:rhoCont} ,
\end{gather}
\end{subequations}
where $\rho$ is the scalar density and $p$ is the pressure. For the incompressible equations, the density is assumed to be constant, in which case the continuity equation \eqref{eq:rhoCont} reduces to the divergence-free condition
\begin{gather}
\div \vu = 0 .
\end{gather}

We let $\vw = \curl \vu$ be the vorticity vector field, which is divergence free by construction. The vorticity equation can be expressed as follows
\begin{subequations}
\begin{gather}
\partial_t \vw + (\vu \cdot \grad) \vw = (\vw \cdot \grad)\vu - \vw (\div \vu), \label{eq:vtxEqn} \\
\div \vu = 0  .
\end{gather}
\end{subequations}

Since we take a geometric approach, it is more convenient to express the vorticity equation in terms of differential forms. Define the vorticity 2-form $\vtx = \star(\vw^\flat)$, where $\flat$ is the lowering of the tensor index and $\star$ is the Hodge-star operator (see Lang \cite{lang2012fundamentals} page 418). 
%\Kai{Could be useful citing a textbook for the notation or add a one line definition.} 
Formally, the 2-form $\vtx$ performs linear measurements on infinitesimal 2D surfaces by dot product of $\vw$ with the surface normal; by definition of curl, this measurement yields the total circulation of $\vu$ along the surface boundary. The vorticity equation \eqref{eq:vtxEqn} is equivalent to the Lie-transport equation of the vorticity 2-form $\vtx$:
\begin{subequations} \label{eqs:LieAdvVtx}
\begin{gather} 
\partial_t \vtx + \Lie_\vu \vtx = 0 , \\
\partial_t \rho + \div(\rho \vu) = 0 .
\end{gather}
\end{subequations}

The Cartan formula for the Lie derivative: $\Lie_\vv \phi = d (\iota_\vv \phi) + \iota_\vv d\phi$, where $\iota_\vv$ is the interior product, is used to recuperate the original vorticity vector equation\footnotemark:
\footnotetext{Using the following identities:
\begin{align*}
&( \iota_\vv \phi )^\sharp  = \vv \cdot \phi^\sharp \quad \text{when $\phi$ is a 1-form} , \quad \quad  &&( \iota_\vv \phi )^\sharp  = (\star\phi)^\sharp \times \vv \quad \text{when $\phi$ is a 2-form} , \\
&(d \phi )^\sharp  = \grad \phi \quad \text{when $\phi$ is a 0-form, and}  \quad \quad &&(d \phi )^\sharp  = \curl \phi^\sharp \quad \text{when $\phi$ is a 1-form} .
\end{align*}
}
\begin{gather} \label{eq:deriveLieAdv}
\left(\star \left( \partial_t  \vtx +  \Lie_\vu  \vtx \right) \right)^\sharp =  \partial_t \vw + \left( \star d (\iota_\vu \vtx) \right)^\sharp + \left( \star \iota_\vu d  \vtx \right)^\sharp = \partial_t \vw + \curl \left(  \vw \times \vu \right) \nonumber  \\
= \partial_t \vw + \vw (\div \vu) - (\vw \cdot \grad ) \vu + (\vu \cdot \grad) \vw = 0 ,
\end{gather}
where $\vw \times \vu$ is the Lamb vector. This means that the vorticity 2-form is conserved, i.e. it is Lie-advected or ``frozen into'' the flow. We note that this uses the general barotropic compressible version of the vorticity equation. In principle, in the incompressible case, the compression term $\vw (\div \vu)$ vanishes, however, this would mean that the Lie-advected vorticity would be allowed to intensify due to numerical errors on the divergence-free condition.

\begin{remark}
One can alternatively check that with $\vw$ evolving strictly under the incompressible equations, the 2-form $\rho \vtx$ is Lie-advected, that is\footnotemark[\value{footnote}], 
\[\left(\star \left( \partial_t \rho \vtx + \Lie_\vu  \rho \vtx \right) \right)^\sharp = \rho ( \partial_t \vw + (\vu \cdot \grad) \vw - (\vw \cdot \grad)\vu) = 0,\]
where $\rho$ is still assumed to satisfy the continuity equation \eqref{eq:rhoCont} (in case the discretized $\vu$ is not exactly divergence-free). Then the vorticity field can be obtained from the Lie-advected $\rho \vtx$ by scaling by $\rho^{-1}$, this cancels any stretching of vortices due to artificial volume compression from numerical errors in $(\div \vu)$.
\end{remark}

In the context of incompressible fluids, the above simply reduces to the statement that the vorticity 2-form is Lie-advected by the velocity field. This gives us an expression of the vorticity as the pullback of the initial condition by the characteristic map:
\begin{gather} \label{eq:vtxPB}
\vtx (\cdot, t) = {\vX_{[t, 0]}}^*\vtx_0 ,
\end{gather}
where the superscript asterisk denotes pullback. For a mapping $F : U \to U$ the pullback $F^*$ it is the dual operator to the pushforward operator denoted by the subscript asterisk $F_*$. The pullback of a $k$-form $\eta$ is defined by $(F^*\eta) (\vv) = \eta \left(F_*\vv  \right)$ where $\vv$ is an arbitrary $k$-vector representing an infinitessimal $k$-dimensional oriented parallelogram and the pushforward $F_*\vv $ is its image under the mapping $F$. Hence, for the 2-form $\vtx$, by the generalized Stokes' theorem, equation \eqref{eq:vtxPB} is equivalent to the conservation of circulation along all closed curves transported by the forward flow map.

Equation \eqref{eq:vtxPB} provides many simplifications both numerically and in the analysis mainly due to the fact that pullback commutes with exterior derivatives. For instance, in the study of the Euler equations through Clebsch variables, one makes the simplifying assumption that the initial velocity 1-form is given by $f dg + d \psi$ for some scalar functions $f$, $g$ and $\psi$. The initial vorticity is then given by $df \wedge dg$. Applying \eqref{eq:vtxPB} to this initial condition and commuting pullback and $d$, we get that the vorticity 2-form at time $t$ is given by
\begin{gather}
\vtx(\cdot, t) = d \left( f \circ \vX_{[t, 0]} \right) \wedge d  \left( g \circ \vX_{[t, 0]} \right) ,
\end{gather}
that is, it is sufficient to solve the advection equations for $f$ and $g$ and reconstruct the vorticity by a cross product of their gradients. We note that the helicity scalar field is defined as $h = \vu \cdot \vw$ which corresponds to the volume form $\vu^\flat \wedge \vtx$. This implies that the Clebsch variable representation is limited to cases where $h$ is exact, i.e. $h = d \phi$ for some 2-form $\phi$. It follows that total helicity is 0 for flows admitting Clebsch variables, i.e. non-helical flows. A generalized version of the Clebsch approach has been studied in \cite{deng2005level}, these Generalized Clebsch variables can be used to represent any initial condition, including helical flows. In fact, the initial velocity expansion \eqref{eq:initVeloExpansion} used in this paper can be seen as a special case of these variables. 

For the numerical method described here, we proceed in the following general setting. We assume that there exists closed 1-forms denoted (by abuse of notation) $d\theta_1, d\theta_2, \ldots, d\theta_n$ and scalar functions $u_1, u_2, \ldots, u_n$ such that the initial velocity 1-form $\vu^\flat$ can be expressed as
\begin{gather}
\vu^\flat = \sum_{k=1}^n u^k d \theta_k .
\end{gather}
Then, the initial vorticity form is given by
\begin{gather}
\vtx_0 = \sum_{k=1}^n d u^k \wedge d \theta_k .
\end{gather}
This gives us a closed expression for the vorticity depending only on $\vX_{[t, 0]}$:
\begin{gather} \label{eq:vtxFormPB}
\vtx(\cdot, t) = \sum_{k=1}^n d \left( u^k \circ \vX_{[t, 0]} \right) \wedge {\vX_{[t, 0]}}^* d \theta_k .
\end{gather}

We do not require that the $d \theta_k$ 1-forms be exact, as long as the pullback is easy to compute. In fact, for the algorithm implemented in this work, on the 3D torus, the $d \theta_k$ forms, with $n=3$, are simply the coordinate covectors $(1,0,0)$, $(0,1,0)$ and $(0,0,1)$, $u_k$ are the corresponding coordinate values of $\vu_0$ and the pullback ${\vX_{[t, 0]}}^* d \theta_k$ is given by $\partial_k \vX_{[t, 0]}$.

Using the following expansion for the initial velocity
\begin{gather} \label{eq:initVeloExpansion}
\vu_0^\flat = u^1 (1,0,0) + u^2 (0,1,0) + u^3(0,0,1) ,
\end{gather}
we get that the vorticity vector at time $t$ is given by
\begin{gather} \label{eq:vtxVeloPB}
\vw (\cdot, t) = \sum_{k=1}^3 \left( \grad u^k \cdot \grad \vX_{[t, 0]} \right) \times \grad \vX_{[t, 0]},
\end{gather}
which further simplifies to
\begin{gather}
w^i (\cdot, t) = \epsilon_{ijk} \epsilon_{abc} w^a_0 \partial_j \vX_{[t,0]}^b \partial_k \vX_{[t,0]}^c .
\end{gather}
in summation notation, where $\epsilon$ are the Levi-Civita symbols.

Upon further inspection, the above expression is %the 
Cramer's rule expansion of 
\begin{gather} \label{eq:vtxVecPB}
\vw (\cdot, t) = \det \left( \grad \vX_{[t,0]} \right)  \left( \grad \vX_{[t,0]} \right)^{-1} \vw_0 (\vX_{[t,0]}) =  \left( \grad \vX_{[t,0]} \right)^{-1} \vw_0 (\vX_{[t,0]}) ,
\end{gather}
where the determinant factor can be omitted since the maps are volume preserving.

\begin{remark}
The basis 1-forms $d \theta_k$ are chosen here to express general initial conditions on the torus. In specific cases, for instance in the presence of Clebsch variables, the number of basis 1-forms can be reduced to improve computational performance. That is, the computation of the characteristic map allows for a flexible framework where the vorticity at time $t$ can be constructed by a pullback formula \eqref{eq:vtxFormPB}, not limited to the formula in \eqref{eq:vtxVecPB}.
\end{remark}

\begin{remark}
An equivalent formulation can be obtained from a Lagrangian perspective by considering the forward map. Following characteristic curves $\vg$, we see that the vorticity field satisfies
\begin{gather}
\frac{d}{dt} \vw(\vg(t), t) = \grad \vu \cdot \vw(\vg(t), t).
\end{gather}
Noticing that the gradient of the forward map evolves according to $\partial_t \grad \vX_F = \grad \vu \cdot \grad \vX_F$, one can show that the vorticity field at $\vX_F$ is given by
\begin{gather}
\vw(\vX_F, t) = \grad \vX_F \cdot \vw_0(\vx) .
\end{gather}
This is Cauchy's Lagrangian formula used in many Lagrangian particle approaches \cite{beale1982vortex, constantin2001eulerian}. Composing the above equation with $\vX_B$ to return to Eulerian frame and applying the inverse function theorem we get
\begin{gather}
\vw (\cdot, t) = \left( \grad \vX_B \right)^{-1} \vw_0 (\vX_B) .
\end{gather}
The factor $\det \left( \grad \vX_B \right)$ discrepancy with \eqref{eq:vtxVecPB} does not show up in the incompressible case since the transformations are volume preserving. In fact, using that $\rho \vtx$ is Lie-advected and $\rho(\vx, t) = \rho_0 (\vX_B) \det \left( \grad \vX_B \right)$, isolating $\vtx$ from $\rho \vtx$ would remove the determinant factor. The CM method for a compressible flow has been studied in \cite{CMdiff} in the context of diffusion-driven density transport.
\end{remark}

\begin{remark}
The evolution of the vorticity 2-form through pullback by $\vX_B$ is the infinitessimal expression of the Kelvin circulation theorem which states that the total circulation along a closed curve passively carried by the fluid flow is constant. As a matter of fact, equation \eqref{eq:vtxVeloPB} can be obtained directly by applying the Kelvin circulation theorem to the definition of the curl operator. Formally, considering that $\vw  \cdot \vn = \lim_{|S| \to 0} \frac{1}{|S|} \int_{\partial S} \vu \cdot \d \vs $ for some infinitessimal surface $S$ with unit normal $\vn$, we apply the Kelvin circulation theorem to $\partial S$, moving it to its position and shape at time $t=0$, to obtain $\vw \cdot \vn = \lim_{|S| \to 0} \frac{1}{|S|} \int_{\vX_B(\partial S)} \vu_0 \cdot \d \vs =  \lim_{|S| \to 0} \frac{1}{|S|} \int_{\partial S} \vu_0(\vX_B) \cdot \d \vX_B(\vs)$ which yields equation \eqref{eq:vtxVeloPB} after taking the limit.

This also relates the CM method to the Kelvin-filtered turbulence models which can roughly be summarized by the vorticity equation
\begin{gather}
\partial \vw  + (\vv \cdot \grad)\vw = (\vw \cdot \grad) \vv
\end{gather}
where $\vv$ is a filtered version of the velocity field $\vu$, for instance $\vv = (I - \alpha^2 \Laplace)^{-1} \vu$. 
We
refer to \cite{foias2001navier} for a review on the theory of these nonlinearly dispersive equations. In the inviscid case, the vorticity is given through pullback by a modified flow map; for the Kelvin filtered equations, the modified flow map is obtained from the filtered velocity field $\vv$, in the CM framework, the modification on the flow map is a result of a combination of filtering and numerical errors.
\end{remark}

\begin{remark}
The pullback formulation allows us to quickly check the conservation of total helicity. Indeed, given that the vorticity is the curl of the velocity, i.e. $d \vu^\flat = \vtx$, we have that by Helmholtz, there exists a 1-form $\eta$ and a scalar $\psi$ such that $\vu^\flat = \eta + d \psi$ and consequently, the vorticity can be written as $\vtx = d \eta$. 
Furthermore, since $\vtx = {\vX_B}^* \vtx_0$, we also have that there exists some $\eta_0$ such that $\eta = {\vX_B}^* \eta_0$ and $d\eta_0 = \vtx_0$. The local helicity is given by $\vu \cdot \vw$ which corresponds to the 3-form $\vu^\flat \wedge \vtx = \eta \wedge d \eta + d \psi \wedge d \eta$. Noting that $d \psi \wedge d \eta = d(\psi d \eta)$ is exact and so has vanishing total integral, we have
\begin{gather}
\int_U \vu^\flat \wedge \vtx = \int_U \eta \wedge d \eta = \int_U {\vX_B}^* (\eta_0 \wedge d\eta_0) = \int_{\vX_B(U)}\eta_0 \wedge d\eta_0 = \int_U \vu_0^\flat \wedge \vtx_0 .
\end{gather}
This property relies only on the fact that both $\eta$ and $\vtx$ evolve through pullback by the same $\vX_B$.
\end{remark}

The derivations in this section allow us to express the evolution of the vorticity field using the characteristic map. The evolution of the characteristic maps is in turn given by the velocity field, which we compute from the vorticity field using the Biot-Savart law $\vu = - \PSolve \curl \vw$. The fully coupled vorticity-characteristic map equations are as follows
\begin{subequations} \label{eqs:CMsummary}
\begin{gather}
(\partial_t + \vu \cdot \grad)  \vX_B = 0 \\
\vw (\cdot, t) = \left( \grad \vX_B \right)^{-1} \vw_0 (\vX_B) \\
\vu = - \PSolve \curl \vw .
\end{gather}
\end{subequations}

\section{Numerical Implementation} \label{sec:numImpl}

The numerical approach of the CM method for the 3D Euler equations largely follows the framework of the Gradient-Augmented Level-Set \cite{nave2010gradient} and Jet-Scheme methods \cite{seibold2011jet}. In this section, we will first present the general numerical framework for the CM method. We will then discuss the specific implementation details used for the numerical experiments in this work.

The computational method is mainly based on the three equations in \eqref{eqs:CMsummary}. Numerically, this corresponds to evolving a characteristic map $\vhX_{[t_n, 0]}$ in a finite dimensional discretization space $\fV $ approximating \textit{Diff}($U$), the space of diffeomorphisms on $U$. Then at every time step, the discretized vorticity $\vtw^n$ and velocity $\vtu^n$ will be reconstructed by pullback using $\vhX_{[t_n, 0]}$. For the rest of this paper, we will denote by the script letter $\vhX$ the approximation of the characteristic map $\vX$ in $\fV$, the superscript $n$ on a variable will denote the evaluation of said variable at time $t_n$ and the tilda indicates an approximation or modified equation. The numerical method comprises the three following parts:
\begin{enumerate}
\item[1.] A discretized velocity field $\vtu^n$ at time $t_n$, (assuming the characteristic map $\vhX_{[t_n, 0]}$ is known). This is given by the Biot-Savart law 
\begin{gather}
\vtu^n = -  \PSolve \curl \left( (\grad \vhX_{[t_n,0]})^{-1} \vw_0(\vhX_{[t_n, 0]} ) \right)
\end{gather} 
computed using spectral Fast-Fourier transform methods.
\item[2.] A numerical approximation $\vtX_{[t_{n+1}, t_n]}$ of the one-step map $\vX_{[t_{n+1}, t_n]}$. For instance, a first order approximation would give 
\begin{gather} \label{eq:egOSmap}
\vtX_{[t_{n+1}, t_n]}(\vx) \approx \vx - \incr{t} \vtu^n(\vx).
\end{gather}
In the current implementation, we use a third-order Runge-Kutta backward in time integrator with a Hermite cubic time-interpolation of the velocity field. A Hermite in time interpolation is used instead of the Lagrange interpolation in \cite{CME} to improve accuracy. 
\item[3.] A time update for the characteristic map based on the group property \eqref{eq:groupCompose}. The map $\vhX_{[t_{n+1}, 0]}$ at time $t_{n+1}$ is given by
\begin{gather} \label{eq:egMapUpdate}
\vhX_{[t_{n+1}, 0]} = \He \left[\vhX_{[t_n, 0]} \circ  \vtX_{[t_{n+1}, t_n]} \right] ,
\end{gather}
for some interpolation operator $\He : \textit{Diff}(U) \to \fV$.
\end{enumerate}

In the following sections, we will examine each component of the algorithm in detail.

\subsection{Spatial Discretization}
In the method presented here, the domain $U = \mathbb{T}^3$ is discretized using a Cartesian meshgrid denoted $\gG$. In 3D, we label the grid points of $\gG$ by $\vx_{i,j,k}$, with each $\vx_{i,j,k}$ located at the lower corner of the cell $C_{i,j,k}$. Given a grid $\gG$, we define the Hermite cubic interpolation operator (for details see e.g. \cite{kolomenskiy2016adaptive, nave2010gradient}) using piecewise smooth basis functions which are tricubic in each cell $C_{i,j,k}$. These basis functions are constructed by tensor product of the 1D Hermite cubic basis functions $Q_i$, with $Q_0$ interpolating function value, and $Q_1$, the derivative:
\begin{subequations} \label{eqs:HermiteBasis1D}
\begin{align}
& q_0 (s) = (1 + 2|s|)(1 - |s|)^2 , \\
& q_1 (s) = s(1 - |s|)^2  , \\
& Q_i(s) = q_i(s) \indi_{[-1, 1]}(s) ,
\end{align}
\end{subequations}
where $\indi$ denotes the indicator function. The $Q_0$ and $Q_1$ functions form the shape functions corresponding to the function value and derivative interpolation on a 1D grid. They have the property that $\partial^a Q_b (s) = \delta^a_b \delta^0_s$ for $a, b \in \{0,1\}$ and $s \in \{-1, 0, 1\}$. The 3D shape functions on the grid $\gG$ are then defined by
\begin{gather}
H_{\va}^{\vi} (\vx) = \prod_{m=1}^3 Q_{a_m} \left( \frac{x^m - x^m_{\vi}}{\incr{x}^m} \right) (\incr{x}^m)^{a_m}
\end{gather}
where $\vx = (x^1, x^2, x^3)$, $\vi = (i_1, i_2, i_3)$ and $\va = (a_1, a_2, a_3)$.

On grid points $\vx_\vr$ of $\gG$, these basis functions satisfy
\begin{gather}
\partial^{\vb}  H_{\va}^{\vi} (\vx_\vr) = \delta_\va^\vb \delta_\vr^\vi .
\end{gather}

The $H_{\va}^{\vi}$ functions are supported on the 8 cells surrounding the grid point $\vx_\vi$, are tricubic in each cell and are globally $\fC^1$ with continuous mixed derivatives of the form $\partial^\vb$ with $b_m \leq 1$ for $m=1,2,3$. We define the Hermite cubic interpolation operator for a smooth function $f$ to be
\begin{gather}
\He_\gG [ f ] (\vx) =  \sum_{\va \in \{0,1\}^3} \sum_{\vs \in \gG} \partial^\va f (\vx_\vi) H_{\va}^{\vi} (\vx) .
\end{gather}
The minimum regularity requirement for this to be well-defined is that $\partial^\vb f$ is continuous for all mixed partial derivatives involving at most one derivative in each Cartesian coordinate. We denote $\fK^1(U) = \left\{ f \in \fC^1(U) \: | \: \partial^\vb f \in \fC^0(U) \text{ for } \vb \in \{0,1\}^3 \right\}$. Then $\He_\gG : \fC^1(U) \to \fK^1(U)$ is a projection operator. Furthermore, for $f \in \fC^4(U)$, it is known that
\begin{gather}
\| \partial^\va f - \partial^\va \He_\gG [f] \|_\infty \lesssim \incr{x}^{4 - | \va |} | f |_{\fC^4} ,
\end{gather}
where $| \va | = a_1 + a_2 + a_3$ and $| f |_{\fC^4} = \sum_{| \vb | = 4 } \| \partial^\vb f \|_{\infty}$ is the $\fC^4$ seminorm; more detailed error bounds can be found in \cite{birkhoff1968piecewise}.

For the characteristic maps, $\vhX$ is given by three coordinate functions. We defined the discretization subspace $\fV = \left( \fK^1 (U) \right)^3$ and use coordinate-wise Hermite cubic interpolation as interpolation operator for the evolution of the characteristic maps. Here we note that $\fV$ is not necessarily contained in $\textit{SDiff}(U)$ or even $\textit{Diff}(U)$ as there is no \emph{a priori} constraint on the determinant or the invertibility of the interpolated map. However, when interpolating a volume-preserving diffeomorphism, the error on the Jacobian determinant is $\bigO(\incr{x}^3)$ and therefore, for sufficiently well-behaved maps, the interpolant will be a diffeomorphism. For longer time simulations, the characteristic map will develop strong small scale features which cannot be resolved using a fixed grid, in those cases, a remapping method will be employed to decompose the transformation; we will examine this in later sections. 

\subsection{Velocity Interpolation}
Using the Hermite cubic interpolation in the previous section, the numerical characteristic map $\vhX_{[t_n, 0]}$ is defined as a diffeomorphism of the domain $U$. We define the numerical vorticity $\vtw^n$ through pullback by $\vhX_{[t_n, 0]}$:
\begin{gather} \label{eq:defNumVort}
\vtw^n (\vx) = \left(\grad \vhX_{[t_n, 0]} \right)^{-1} \vw_0 (\vhX_{[t_n, 0]}(\vx) ) ,
\end{gather}
where the gradient of $\vhX_{[t_n,0]}$ is directly evaluated from the interpolant. This defines $\vtw^n$ as a $\fC^0$ vector field on $U$. The numerical velocity $\vtu^n$ is in turn computed from the convolution of the Biot-Savart kernel with $\vtw^n$; we will do this using Fourier spectral methods. We will discretize the velocity field on a grid $\gV$ and denote by $\fF_\gV$ the discrete Fourier transform computed by FFT on the grid $\gV$. Sampling $\vtw^n$ on $\gV$ using \eqref{eq:defNumVort} and applying a forward Fourier transform yields a truncated Fourier series for the vorticity field. Computing the Biot-Savart kernel in frequency space then gives us a Fourier series representation for the velocity. Finally, we define the numerical velocity $\vtu^n$ as the Hermite interpolant of this truncated Fourier series; this allows us to evaluate the velocity at arbitrary locations in the domain without having to compute the inverse Fourier transform at non-uniform grid points. The definition of the numerical velocity at $t_n$ can be summarized as follows:
\begin{gather} \label{eq:defNumVelo}
\vtu^n(\vx) = \He_\gV \left[ \fF^{-1}_\gV \left[ - \PSolve \curl \fF_\gV [ \vtw^n ]  \right] \right] (\vx) .
\end{gather}
In the above equation, it is understood that in order to define the Hermite cubic interpolant for the velocity, the required mixed partial spatial derivatives are computed directly from the Fourier series.

Similarly, from $\partial_t \vw = (\vw \cdot \grad) \vu - (\vu \cdot \grad)\vw$, we can also discretize the time derivative of $\vu$ at $t_n$:
\begin{gather} \label{eq:defNumdVelodt}
\partial_t \vtu^n(\vx) = \He_\gV \left[ \fF^{-1}_\gV \left[ - \PSolve \curl \fF_\gV [ (\vtw^n \cdot \grad) \vtu^n - (\vtu^n \cdot \grad)\vtw^n ]  \right] \right] (\vx) .
\end{gather}

The data $\vtu^n$ and $\partial_t \vtu^n$ at time steps $t_n$ allow us to locally approximate the velocity using a 4-dimensional time-space Hermite cubic interpolant. For the one-step map \eqref{eq:egOSmap}, we will define an approximate $\vtu$ in the interval $[t_n, t_{n+1}]$ by extending the interpolant obtained from the velocity data $\vu^{n-1}$ and $\vu^n$ in the interval $[t_{n-1}, t_n]$. This gives the following definition for the numerical velocity field:
\begin{align} \label{eq:defNumFullVelo}
\vtu (\vx, t) = & \left( q_0(t-t_{n-1}) \vtu^{n-1}(\vx) + q_0(t-t_{n}) \vtu^n(\vx) \right) \nonumber \\ 
& + \incr{t} \left( q_1(t-t_{n-1}) \partial_t \vtu^{n-1}(\vx) + q_1(t-t_{n}) \partial_t \vtu^n(\vx) \right) \quad \text{for} \quad t \in [t_n, t_{n+1}),
\end{align}
using the Hermite basis functions given in \eqref{eqs:HermiteBasis1D}. We note that since for each $n$, $\vtu^n$ and $\partial_t \vtu^n$ are Hermite interpolants of divergence-free vector fields, the modified velocity field $\vtu (\vx, t)$ is a linear combination of divergence-free velocity fields and is also divergence-free at all time up to interpolation error. This error can be reduced by refining the velocity interpolation grid which can be achieved by a zero-padding in frequency space before taking the inverse Fourier transform.

The one-step map in the interval $[t_n, t_{n+1}]$ is then obtained from the backward in time flow of the approximate velocity field $\vtu$. We define the numerical one-step map $\vtX_{[t_{n+1}, t_n]}$ pointwise using a third order Runge-Kutta integration of $\vtu$:
\begin{gather} \label{eq:defNumOSmap}
\vtX_{[t_{n+1}, t_n]} = \int_{t_{n+1}}^{t_n} \vtu \left(\vtX_{[\tau, t_n]} , \tau \right) \d \tau .
\end{gather}

The one-step map is used in the time update of the characteristic map \eqref{eq:egMapUpdate}; it is therefore only evaluated at grid points. However, in order to compute the chain rules for the derivatives required to define the Hermite interpolant, we also need the mixed partial derivatives of $\vtX_{[t_{n+1}, t_n]}$ at grid points. This is done using a $4^{th}$ order version of the $\epsilon$-difference scheme described in \cite{chidyagwai2011comparative}. The $\epsilon$-difference schemes introduce an $L^\infty$ error of order $\epsilon^4 + \delta \sum_{k=0}^3 \incr{x}^k \epsilon^{-k}$ where $\delta$ is the machine precision. For all computations presented in this paper, we used $\epsilon = 2.5 \times 10^{-3}$ which corresponds to an error term of at most $10^{-8}$ and effectively less than $10^{-11}$, if $\incr{x}<0.1$.

\begin{remark}
Evaluating the extrapolation formula \eqref{eq:defNumFullVelo} at $t = t_{n+1}$ would generally not give the same velocity as $\vtu^{n+1}$ which is obtained by vorticity pullback followed by the Biot-Savart law. Indeed, the extrapolation of the velocity in the interval $[t_n, t_{n+1})$ is used only to evolve the characteristic map to time $t_{n+1}$. The velocity at $t_{n+1}$ is then reconstructed using \eqref{eq:defNumVort} and \eqref{eq:defNumVelo}, similar to a predictor-corrector approach. This also implies that numerical errors in the extrapolation are not directly carried and amplified in the next time step.
\end{remark}

\subsection{Error Estimates} \label{sec:ErrorEst}
We will examine in this section the numerical error on the characteristic map and its relation to the error on the vorticity field. We will try to characterize the nature of the numerical error and provide some estimates. We use as starting assumption that the numerical map $\vhX_B$ is consistent with the exact map $\vX_B$ in the $\fC^{1, \alpha}$ norm for some $\alpha \in (0,1)$, that is the error is $o(1)$. This will allow us to estimate the global $\fC^{1, \alpha}$ error to third-order in time and space by omitting higher order terms in the error. The consistency assumption is then implied for short-time since the initial numerical map $\vhX_{[0,0]}= \vx$ is exact. In order to preserve the advective structure of the error, we define the following error map:
\begin{gather}
\vhE_{[0, t_n]} := \vhX_{[t_n, 0]} \circ \vX_{[t_n, 0]}^{-1} = \vhX_{[t_n, 0]} \circ \vX_{[0, t_n]} = \vhX_B \circ \vX_F,
\end{gather}
which measures by how much the diffeomorphism $\vhX_{[t_n, 0]}$ differs from the inverse of the exact forward flow map $\vX_{[0, t_n]}$. Indeed, since the composition of the forward and backward maps $\vX_{[t_n,0]} \circ \vX_{[0, t_n]} = \vx$ gives the identity map, we have that
\begin{gather}
\vhE_{[0, t_n]} = \vx + \left(\vhX_{[t_n, 0]} - \vX_{[t_n, 0]}  \right) \circ \vX_{[0, t_n]},
\end{gather}
that is, the deviation of $\vhE_{[0, t_n]}$ from the identity map is the numerical error of the map evaluated at the pushforward location. This is a Lagrangian representation of the error since $\vx - \vhE_{[0, t_n]}(\vx)$ essentially gives the time $t_n$ map error for a particle starting at $\vx$ at time $0$. We also note that since $\vhX_B$ is a $\fC^1$ diffeomorphism by construction, it follows that all maps considered here are also $\fC^1$ diffeomorphisms. Thus left and right inverses exist, are equal and are also diffeomorphisms.

We also define the auxiliary ``modified map'' 
\begin{gather}
\vtX_{[t_n, 0]} := \vtX_{[t_1, 0]} \circ  \vtX_{[t_2 , t_1]} \circ \dots \circ \vtX_{[t_n, t_{n-1}]},
\end{gather}
where $\vtX_{[t_k, t_{k-1}]}$ is the one-step map given in \eqref{eq:defNumOSmap} obtained from RK3 integration on the numerical interpolated velocity field. 

The full error map is decomposed as follows:
\begin{gather}
\vhE_{[0, t_n]} = \vhX_{[t_n, 0]} \circ \vtX_{[t_n, 0]}^{-1} \circ \vtX_{[t_n, 0]} \circ \vX_{[t_n, 0]}^{-1} .
\end{gather}

We let $\vPhi_{[0, t_n]} = \vhX_{[t_n, 0]} \circ \vtX_{[t_n, 0]}^{-1}$ and $\vPsi_{[0, t_n]} = \vtX_{[t_n, 0]} \circ \vX_{[t_n, 0]}^{-1}$ and compute their time-evolution as follows:
\begin{gather}
\vPhi_{[0, t_n]} = \He_\gM \left[ \vhX_{[t_{n-1}, 0]} \circ \vtX_{[t_n, t_{n-1}]} \right] \circ \left( \vtX_{[t_{n-1}, 0]} \circ \vtX_{[t_n, t_{n-1}]} \right)^{-1} \nonumber \\  = \He_\gM \left[ \vhX_{[t_{n-1}, 0]} \circ \vtX_{[t_n, t_{n-1}]} \right] \circ \left( \vhX_{[t_{n-1}, 0]} \circ \vtX_{[t_n, t_{n-1}]} \right)^{-1} \circ \vhX_{[t_{n-1}, 0]} \circ \vtX_{[t_{n-1}, 0]}^{-1} \nonumber \\
= \vhxi_n \circ \vPhi_{[0, t_{n-1}]} ,
\end{gather}
where we defined a one-step error $\vhxi_n := \He_\gM \left[ \vhX_{[t_{n-1}, 0]} \circ \vtX_{[t_n, t_{n-1}]} \right] \circ \left( \vhX_{[t_{n-1}, 0]} \circ \vtX_{[t_n, t_{n-1}]} \right)^{-1} $. We note that $\vhxi_n$ is the error due to Hermite interpolation since
\begin{gather}
\vhxi_n - \vx = \left( \He_\gM \left[ \vhX_{[t_{n-1}, 0]} \circ  \vtX_{[t_{n}, t_{n-1}]} \right] -  \vhX_{[t_{n-1}, 0]} \circ  \vtX_{[t_{n}, t_{n-1}]}  \right) \circ \vtX_{[t_{n}, t_{n-1}]}^{-1} \circ \vhX_{[t_{n-1}, 0]}^{-1}  \nonumber \\ = \left( \He_\gM \left[ \vhX_{[t_{n-1}, 0]} \circ  \vtX_{[t_{n}, t_{n-1}]} \right] -  \vhX_{[t_{n-1}, 0]} \circ  \vtX_{[t_{n}, t_{n-1}]}  \right) \circ  \vtX_{[t_{n}, 0]}^{-1} \circ \vPhi_{[0, t_{n-1}]}^{-1}.
\end{gather}

We define the following interpolation error
\begin{gather}
\vphi_n :=\He_\gM \left[ \vhX_{[t_{n-1}, 0]} \circ  \vtX_{[t_{n}, t_{n-1}]} \right] -  \vhX_{[t_{n-1}, 0]} \circ  \vtX_{[t_{n}, t_{n-1}]} ,
\end{gather}
we then obtain that
\begin{gather}
\vPhi_{[0, t_n]} = (\vx + (\vhxi_n - \vx) ) \circ \vPhi_{[0, t_{n-1}]} = \vPhi_{[0, t_{n-1}]} + \vphi_n \circ  \vtX_{[t_{n}, 0]}^{-1} = \vx + \sum_{k=1}^n \vphi_k \circ  \vtX_{[t_{k}, 0]}^{-1} .
\end{gather}

A similar derivation gives us
\begin{gather}
\vPsi_{[0, t_n]} = \vheta_n \circ \vPsi_{[0, t_{n-1}]}  ,
\end{gather}
where $\vheta_n := \left( \vtX_{[t_{n-1}, 0]} \circ \vtX_{[t_n, t_{n-1}]} \right) \circ \left( \vtX_{[t_{n-1}, 0]} \circ \vX_{[t_n, t_{n-1}]} \right)^{-1}$. We note that $\vheta$ is the velocity approximation error which gives the discrepancy between the true flow and the flow obtained from the modified velocity $\vtu$.
\begin{gather}
\vheta_n - \vx = \left( \vtX_{[t_{n-1}, 0]} \circ \vtX_{[t_n, t_{n-1}]} -  \vtX_{[t_{n-1}, 0]} \circ \vX_{[t_n, t_{n-1}]}  \right) \circ \vX_{[t_n, 0]}^{-1} \circ \vPsi_{[0, t_{n-1}]}^{-1} .
\end{gather}

We define the following modified flow error
\begin{gather}
\vpsi_n = \vtX_{[t_{n-1}, 0]} \circ \vtX_{[t_n, t_{n-1}]} -  \vtX_{[t_{n-1}, 0]} \circ \vX_{[t_n, t_{n-1}]} ,
\end{gather}
i.e. the map evolution error due to errors in the approximated velocity field; we note that this term is approximately $\incr{t} \grad \vtX_{[t_{n-1}, 0]} \cdot (\vtu - \vu)$. This gives us
\begin{gather}
\vPsi_{[0, t_n]} = (\vx + (\vheta_n - \vx) ) \circ \vPsi_{[0, t_{n-1}]} = \vPsi_{[0, t_{n-1}]} + \vpsi_n \circ \vX_{[t_n, 0]}^{-1} = \vx + \sum_{k=1}^n \vpsi_k \circ  \vX_{[t_{k}, 0]}^{-1}.
\end{gather}

This allows us to write the error map as
\begin{gather}
\vhE_{[0, t_n]} = \left(  \vx + \sum_{k=1}^n \vphi_k \circ  \vtX_{[t_{k}, 0]}^{-1} \right) \circ \left( \vx + \sum_{k=1}^n \vpsi_k \circ  \vX_{[t_{k}, 0]}^{-1} \right) .
\end{gather}

We use the fact that for two $\fC^{1, \alpha}$ diffeomorphisms, $\bm{f}$ and $\bm{g}$, the composition $\bm{f} \circ \bm{g}$ is also $\fC^{1, \alpha}$ with $\| \bm{f} \circ \bm{g} \|_{\fC^{1, \alpha} } \leq C \| \bm{f}  \|_{\fC^{1, \alpha} } \| \bm{g} \|_{\fC^{1, \alpha} }^{1+\alpha}$, to estimate the norm of the Lagrangian displacement error $\vhe_{[0, t_n]} := \vx - \vhE_{[0, t_n]}$:
\begin{gather} \label{eq:LErrorEst}
\| \vhe_{[0, t_n]} \|_{\fC^{1, \alpha}} \lesssim \sum_{k=1}^n  \left\|  \vpsi_k \right\|_{\fC^{1, \alpha}} +  \left( 1 +   \sum_{k=1}^n \left\| \vpsi_k\right\|_{\fC^{1, \alpha}} \right)^{1+\alpha} \sum_{k=1}^n \left\| \vphi_k\right\|_{\fC^{1, \alpha}} ,
\end{gather}
where the terms dependent on $\vX_{[t_n, 0]}$ have been absorbed in the constants of the inequality; $\vtX_{[t_n, 0]}$ is also approximated with $\vX_{[t_n, 0]}$ using the consistency assumption. Here we use the notation $A \lesssim B$ to denote that there exists some constant $c$ such that $A < cB$. In this section, the constant will depend on the dimension, the domain, the constants involved in the norms, and the solution $\vu$ as well as $\vX_B$.

The $\vphi_n$ error is an error pertaining the numerical resolution of $\vhX_B$, it can be controlled as long as we can control the higher derivatives of $\vhX_{[t_{n-1}, 0]} \circ  \vtX_{[t_{n}, t_{n-1}]}$. The stability analysis of a similar methods has been studied in \cite{goodrich2006hermite}.

The $\vpsi_n$ error is a feedback between the map error and the velocity error (and also numerical integration), with $\vpsi_n \approx \grad \vtX_{[t_{n-1}, 0]} \incr{t} (\vtu - \vu)$ and $\| \vpsi_n \|_{\fC^{1,\alpha}} \lesssim \incr{t} \| \vtu^{n-1} - \vu^{n-1} \|_{\fC^{1,\alpha}}^{1 + \alpha}$. We first bound the error on the vorticity which will require the Eulerian version of the error map. Define
\begin{gather} \label{eq:EEerror}
\vhE_{[t_n, 0]} :=  \vX_{[t_n,0]}^{-1} \circ \vhX_{[t_n,0]} = \vX_{[0, t_n]} \circ \vhE_{[0,t_n]} \circ \vX_{[t_n,0]},
\end{gather}
here the conjugation exactly serves the purpose of changing the error from a Lagrangian frame to an Eulerian frame. 

We note that 
\begin{gather} \label{eq:compareVort}
\vw(\vx,t_n) = (\grad \vX_B)^{-1} \vw_0 (\vX_B) \quad \text{and} \quad \vtw^n(\vx) = (\grad \vhX_B)^{-1} \vw_0 (\vhX_B),
\end{gather}

So we have that, using $\vhX_B = \vX_B \circ \vhE_B$ and $(\grad \vhX_B)^{-1} = \left( \grad \vhE_B \right)^{-1} (\grad \vX_B)^{-1} |_{\vhE_B^{-1}} $
\begin{gather}
\vtw^n = \left( \grad \vhE_B \right)^{-1} \vw^n ( \vhE_B ) \quad \text{i.e.} \quad \tilde{\vtx}^n = {\vhE_B}^* \vtx^n.
\end{gather}
The numerical vorticity is the pullback of the exact vorticity by the error map. This was expected since the exact map was decomposed into the composition of the numerical map and the error map. 

A first look at the velocity error, letting $\vhe_B = \vx - \vhE_B$ assuming $\vhe_B$ small, we have to leading order terms
\begin{gather}
\vw^n - \vtw^n \approx ( \vI - (\grad \vhE_B)^{-1} )\vw^n (\vhE_B) + \grad\vw^n \cdot (\vx - \vhE_B) \nonumber \\
\approx \grad \vhe_B \vw^n + \grad \vw^n \vhe_B .
\end{gather}
up to second order terms of $\bigO(\| \vhe_B \|^2)$.

This gives us $\| \vw^n - \vtw^n \|_{\fC^{0,\alpha}} \lesssim \| \vhe_B \|_{\fC^{1, \alpha}}$ for some $\alpha \in (0,1)$, which after Biot-Savart, yields
\begin{gather}
\| \vu^n - \vtu^n \|_{\fC^{1,\alpha}} \lesssim \| \vhe_B \|_{\fC^{1, \alpha}}
\end{gather}

We use this to control the $\vpsi_n$ error which is due to the difference between the numerical and the exact velocities as well as the interpolation and integration schemes both of which have $4^{th}$ order local truncation error.
\begin{gather}
\| \vpsi_n \|_{\fC^{1,\alpha}} \lesssim  \incr{t} \| \vhe_{[t_{n-1}, 0]} \|_{\fC^{1,\alpha}}^{1+\alpha} + \bigO(\incr{t}^4).
\end{gather}

Since the Lagrangian and Eulerian error maps are related by conjugation by the map $\vX_{[t_n, 0]}$, we have that $\| \vhe_{[t_n, 0]} \|_{\fC^{1,\alpha}} \lesssim  \| \vhe_{[0, t_n]} \|_{\fC^{1,\alpha}}^{1+\alpha}$ and vice-versa. Therefore, we can write the following estimate for the Eulerian error, up to leading order terms, using \eqref{eq:LErrorEst}:
\begin{gather}
\| \vhe_{[t_n, 0]} \|_{\fC^{1,\alpha}}^{1/(1+\alpha)} \lesssim  \sum_{k=1}^n \incr{t}  \| \vhe_{[t_{k-1}, 0]} \|_{\fC^{1,\alpha}}^{1+\alpha} +  \sum_{k=1}^n \left\| \vphi_k\right\|_{\fC^{1, \alpha}} .
\end{gather}
This can be majorized by the ODE by approximating the discrete sum of order $\incr{t}$ terms with the integral from 0 to $t = n\incr{t}$ and taking a time derivative.
\begin{gather} \label{eq:errorODE}
\dot{\vhe}_B = (1 + \alpha)\vhe_B^{\alpha/(1+\alpha)} ( \vhe_B^{1+\alpha} + A),
\end{gather}
where $A = \bigO(\incr{t}^3) + \frac{1}{\incr{t}} \left\| \vphi_k\right\|_{\fC^{1, \alpha}}$, which is the CM advection error and is $\bigO( \incr{t}^3+ \incr{x}^3)$.
 
We note that the H\"older-$\alpha$ norm was introduced artificially to gain the full 2 degrees of regularity from the Poisson equation in the Biot-Savart law. We can therefore pick $\alpha > 0$ arbitrarily small, in which case, the map error estimate $\vhe_B$ in \eqref{eq:errorODE} solves a regularly perturbed $1^{st}$ order linear ODE with a source term of $\bigO(\incr{t}^3 + \incr{x}^3)$. We recall that this is built on the assumption that the modified flow map $\vtX_{[t_{n}, 0]}$ is consistent with $\vX_{[t_{n}, 0]}$ so that higher order error terms can be omitted; this is true since the initial $\vtX_{[t_{0}, 0]}$ is exact and the time-evolution of the error is third-order in $\incr{x}$ and $\incr{t}$ according to the above derivations. It is also assumed that $\vhX_{[t_n,0]} \circ \vtX_{[t_{n+1}, t_n]}$ is well represented by Hermite interpolation, i.e. that the spatial resolution is high enough and that the grid data of $\vhX_{[t_i,0]}$ do not oscillate unboundedly. This can be in part controlled by having high enough resolution and also by a remapping method discussed in the following section. The stability of the CM method was also discussed in \cite{nave2010gradient, CM, CME} and convergence of similar methods using Hermite interpolation was proven in \cite{goodrich2006hermite}. Overall, the CM method should have $\bigO(\incr{x}^3 + \incr{t}^3)$ error for the map in $\fC^1$ norm which would translate to a $3^{rd}$ order global error.

We provide here two numerical tests for the error estimates derived above. As a sanity check, we test the method on the stationary Arnold-Beltrami-Childress (ABC) flow so that the numerical solution can be compared against a known exact solution. We will also perform a second test using a standard Taylor-Green vortex initial condition, numerical results will be compared against a high-resolution reference test.

Both tests are performed on a periodic domain $[-2\pi, 2\pi]^3$. The ABC flow initial condition is given by
\begin{gather}
\vw_0(x,y,z) = \frac12 \left( \begin{matrix}
\cos(y) + \sin(z) \\ \cos(z) + \sin(x) \\ \cos(x) + \sin(y) 
\end{matrix} \right), 
\end{gather}
and the Taylor-Green initial condition is given by
\begin{gather}
\vw_0 (x,y,z) = \left( \begin{matrix} \cos \left( \frac{x}{2} \right) \sin \left( \frac{y}{2} \right) \sin \left( z \right)  \\  \sin \left( \frac{x}{2} \right) \cos \left( \frac{y}{2} \right) \sin \left( z \right)  \\  -\sin \left( \frac{x}{2} \right) \sin \left( \frac{y}{2} \right) \cos \left( z \right) \end{matrix}  \right) .
\end{gather}

\begin{figure}
\centering
\includegraphics[width = 0.425\linewidth]{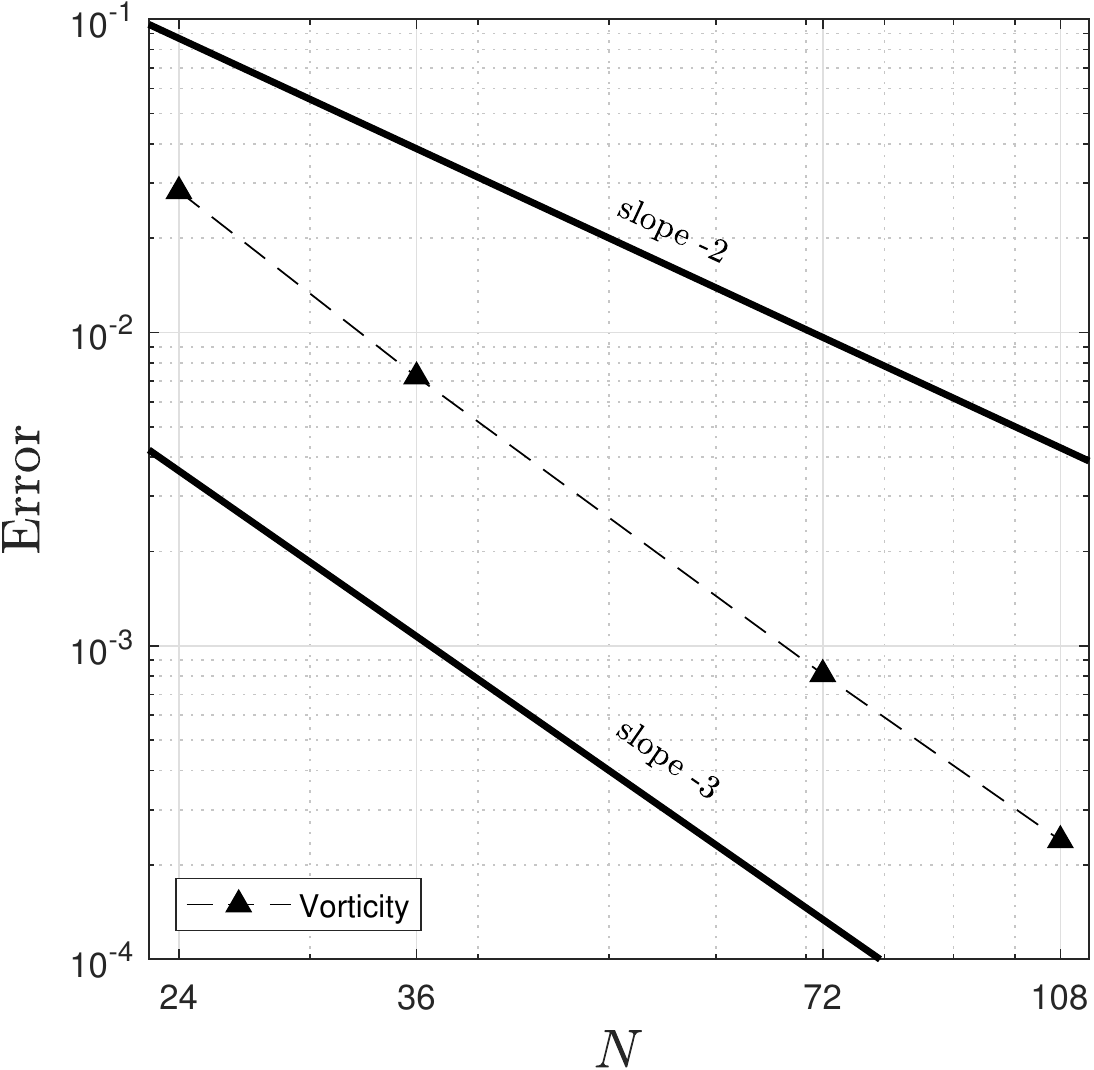}
\caption{$L^\infty$ vorticity error for the ABC test at $T_f=2$. Numerical solution is directly compared to the exact solution $\vw(\vx, t) = \vw_0(\vx)$.}
\label{fig:convPlot1}
\end{figure}

\begin{figure}
\centering
\includegraphics[width = 0.425\linewidth]{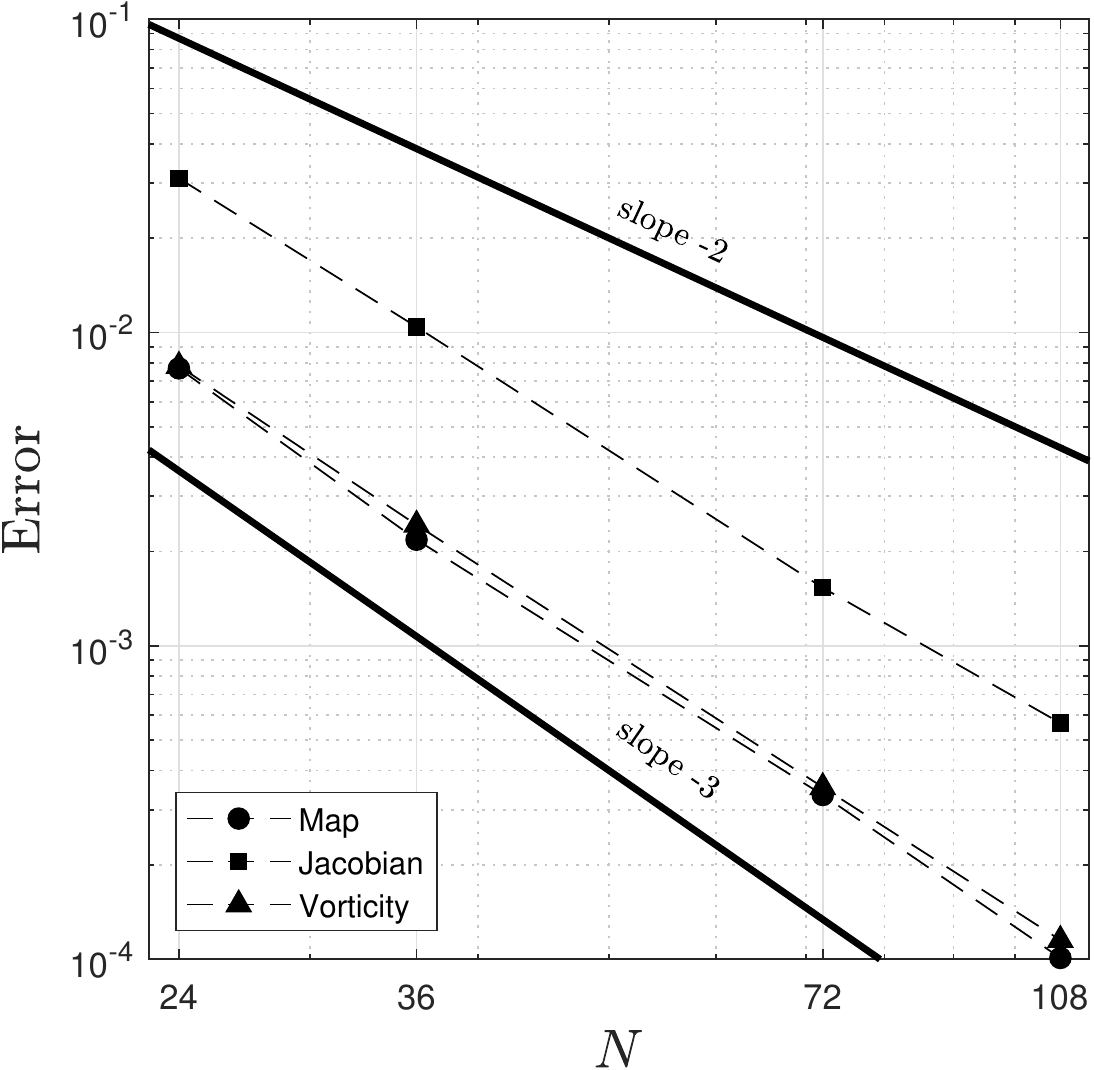}
\caption{$L^\infty$ map, Jacobian and vorticity errors for the Taylor-Green Vortex test at $T_f=2$. Error is calculated by comparing results against a $N=216$ test.}
\label{fig:convPlot2}
\end{figure}

We run each test on increasingly finer grids of $N$ cells per dimension for both $\gM$ and $\gV$ and using $N/12$ time steps to reach the final time. In both cases, the final time is $T_f = 2$, so that $\incr{x} = \frac{4 \pi}{N}$ and $\incr{t} = \frac{24}{N}$. For the ABC test, we measure the vorticity error using the exact solution $\vw^n = \vw_0$ and for the Taylor-Green vortex test, we measure the $\vw$, $\vhX_B$ and $\grad \vhX_B$ errors at grid points by comparing against a reference higher resolution test with $N = 216$. We note that for both tests the maximum velocity throughout the simulation is greater than 0.9 at all times so that $\incr{t}$ exceeds the CFL condition. The maximum vorticity for the ABC test is constant in time, for the Taylor-Green test, a $20\%$ growth is observed over the $[0, 2]$ time interval. Figures \ref{fig:convPlot1} and \ref{fig:convPlot2} show the $L^\infty$ errors for both tests at $T_f= 2$; %for each initial condition; 
the errors are computed directly from grid values and confirm the expected $3^{rd}$ order error.


\subsection{Submap Decomposition} \label{sec:submapDecomp}
The error estimates in section \ref{sec:ErrorEst} show that the error is advective in nature since the numerical solution can be written as the pullback of the exact solution by an error map; therefore, viscous and hyper-viscous type dissipation in the solution are avoided. Indeed, by computing the characteristic map, the vorticity field is provided functionally and can be evaluated anywhere on the domain by interpolating $\vhX_B$. Instead of directly evolving the vorticity on some grid, in which case $\vtw^n$ depends on the grid values of $\vtw^{n-1}$, the vorticity in the CM method is in principle ``reconstructed'' at every step using $\tilde{\vtx}^n = {\vhX_{[t_n, 0]}}^* \vtx_0$.  For traditional grid-based methods, $\vtw^n$ is obtained from grid values of $\vtw^{n-1}$ and typically carries a viscous or hyper-viscous error of the form $\epsilon \Laplace^p \vtw^{n-1}$, in the CM method, the error is instead advective with $\vtw^n = {\vhE_B}^* \vw^n$, thus preventing the loss of subgrid scales due to artificial viscosity. One implication is that the CM method can better avoid artificial merging of vortices at the subgrid scale in particular in the study of vortex tube reconnection problems. The error in this case is of elastic type.  Indeed, since $\vhX_B$ is evolved using the GALS method, the leading order spatial errors on the map are of the form $\epsilon \Laplace^p \vhX_B$. This can be seen as an elasticity term which dampens extensive deformation of the domain under the $\vhX_B$ mapping. The evolution of the elasticity error is governed by the $\vphi$ term in the error map.

The $\vphi$ map representation error depends on the $\gM$ grid and the regularity of the characteristic map $\vX_B$. For Hermite cubic interpolation, this error roughly scales with the $4^{th}$ spatial derivative of $\vX_B$. At time $t=0$, this error is $0$ since $\vX_{[0,0]}$ is the identity map, then the error increases with time as the characteristic map develops more complicated spatial features. We limit the growth of this error by periodically reinitializing the characteristic map using the group property of the flow maps. Indeed, a time $t$ characteristic map can be decomposed as follows using \eqref{eq:groupCompose}:
\begin{gather}
\vX_{[t,0]} = \vX_{[T_1, 0]} \circ \vX_{[T_2, T_1]} \circ \cdots \circ \vX_{[T_{m-1}, T_{m-2}]} \circ \vX_{[t, T_{m-1}]} ,
\end{gather}
for $0 = T_0 < T_1 < \cdots < T_{m-1} < T_m = t$. The $T_i$ are remapping times, in each interval $[T_{i}, T_{i+1}]$, the evolution of the characteristic map is given by the $T_{i+1}-T_i$ time flow of equations \eqref{eqs:CMsummary} with $\vw(\vx, T_{i})$ as initial condition. 

\begin{remark}
Note that the $T_i$ refer to the remapping times and $t_n$ are the time steps of the scheme with $t_n-t_{n-1} = \incr{t}$ small, approximating the limit to 0. On the other hand  $T_i - T_{i-1}$ is $\bigO(1)$, its purpose is to subdivide the time interval $[0, T_{final}]$ into shorter subintervals where the characteristic maps are better behaved.
\end{remark}

At any $t$, the vorticity is given by the following pullback using submap decomposition:
\begin{gather}
\vw(\vx, t) = \left(\grad \vX_{[t, T_{m-1}]}  \right)^{-1}\cdots  \left(\grad \vX_{[T_{2}, T_{1}]}  \right)^{-1} \left(\grad \vX_{[T_1, 0]}  \right)^{-1} \left( \vw_0 \circ \vX_{[T_1, 0]}  \circ \cdots  \circ \vX_{[t, T_{m-1}]}  \right) .
\end{gather}

Numerically, this means that we can compute each submap individually and use
\begin{gather} \label{eq:submapCompNum}
\vtw(\cdot, t) = \left( \prod_{i=0}^{m-1}  \left(\grad \vhX_{[T_{m-i}, T_{m-i-1}]}  \right)^{-1}  \right) \vw_0 \left( \vhX_{[T_1, 0]}  \circ \cdots  \circ \vhX_{[t, T_{m-1}]}  \right)
\end{gather}
to compute the pullback. Each submap will only perform the mapping in the time subinterval $[T_i, T_{i+1}]$ and remapping times can be either fixed or chosen dynamically so that each submap can be well represented using the grid $\gM$. Indeed, after each remapping, we essentially solve a separate Euler equation in the time interval $[T_i, T_{i+1}]$, the characteristic map is reinitialized with the identity map. This means that the spatial representation of the current submap is again exact and will accumulate error over $[T_i, T_{i+1}]$, until an error threshold is exceeded and remapping is triggered. In terms of the error estimates in section \ref{sec:ErrorEst}, the remapping resets the $\vphi$ error of the $i^{th}$ submap to 0, which prevents further accumulation of spatial interpolation errors due to the fixed resolution of the $\gM$ grid. One can therefore control the accumulation of the elasticity type error by changing the frequency of the remapping. More frequent remapping reduces the effect that the $\gM$ grid has on the spatial features of the map; in the extreme case where remapping is done at every time step, the $\vphi$ error can no longer accumulate and numerical error is reduced to $\vpsi$ only. For the numerical experiments presented in this paper, we use the volume-preservation error of the $\vhX_B$ map as remapping criterion, that is, the remapping time $T_i$ is chosen to be the first time $t$ such that the error $|\det \grad \vhX_{[t, T_{i-1}]} - 1|$ is greater than some chosen tolerance. On one hand, this serves as an \emph{a posteriori} estimate of the $\fC^1$ error of the map and on the other hand, this allows us to guarantee that the composition of all submaps yields a diffeomorphism and provides some control on the overall volume-preserving property of the characteristic map.

\subsection{Implementation Summary}
The previous subsections contain the numerical tools for implementing the CM method for the 3D incompressible Euler equations, we give here a short summary of the method in pseudocode format. We note that the method uses two discretization grids, a grid $\gM$ for representing the numerical map $\vhX_{[t_n, 0]}$ and a  grid $\gV$ for sampling the numerical vorticity $\vw^n$ and computing the Biot-Savart law using Fourier spectral methods. These two grids do not need to have the same resolution, in fact, using the submap decomposition method described in the previous subsection, a short time characteristic map can be represented on a coarse grid $\gM$. The $\gV$ grid used to represent the vorticity needs to be fine enough to avoid sampling errors. Indeed, as the flow evolves, the vorticity can develop small scale features and high gradients. If $\gV$ is not fine enough to resolve $\vw^n$, this can cause aliasing errors in the Fourier transform. One way to reduce the effect of undersampling is to mollify $\vw^n$ in Fourier space which was studied in \cite{CME}, however this also reduces the accuracy of the scheme. Other possibilities for future investigation include the use adaptive mesh for the vorticity sampling and wavelet methods for the computation of the Biot-Savart law \cite{schneider1997comparison, schneider2010wavelet}. The CM method for 3D incompressible Euler equations can be summarized with the following pseudocode algorithm.

\begin{algorithm}
  \caption{CM method 3D incompressible Euler equations.
    \label{alg:CME3D}}
  \begin{algorithmic}[1]
    \Require{Initial vorticity $\vw_0$, grids $\gM$ and $\gV$, time step $\incr{t}$, final time $T_f$}
      \State Initialize $n \gets 0$, $t_n \gets 0$, $m \gets 0$, $T_m \gets 0$.
      \While{$T_m < T_f$}
      \State $\vhX_{[t_n, T_m]} \gets \vx$ (identity map)
      \While{$\| \det \grad \vhX_{[t_n, T_m]} -1  \|_\infty < TOL $}
      \State Sample $\vw^n$ on grid points of $\gV$. \Comment{using \eqref{eq:submapCompNum}}
      \State $\fF_\gV [\vtu^n] = - \PSolve \curl \fF_\gV [ \vtw^n ]$.
      \State Compute $\partial^\vb \vtu^n$ for $\vb \in \{0,1\}^3$ in Fourier space. Define $\vtu^n(\vx)$.\Comment{using \eqref{eq:defNumVelo}}
      \State Compute $(\vtw^n \cdot \grad) \vtu^n - (\vtu^n \cdot \grad)\vtw^n$ on $\gV$ and define $\partial_t \vtu^n(\vx)$. \Comment{using \eqref{eq:defNumdVelodt}}
      \State Define $\vtu(\vx, t)$, by linear combination of spatial interpolants. \Comment{using \eqref{eq:defNumFullVelo}}
      \State Compute $\partial^\vb \vtX_{[t_n + \incr{t}, t_n]}$ on $\gM$ using RK3 integration of $\vtu$. \Comment{using \eqref{eq:defNumOSmap}}
      \State Update characteristic map $\vhX_{[t_{n+1}, T_m]} \gets \He_{\gM} \left[ \vhX_{[t_{n}, T_m]} \circ \vtX_{[t_n + \incr{t}, t_n]} \right]$.
      \State $t_{n+1} \gets t_n + \incr{t}$, $n \gets n+1$.
      \EndWhile
      \State $T_{m+1} \gets t_n$, $m \gets m+1$.
      \EndWhile
  \end{algorithmic}
\end{algorithm}
\begin{remark}
Since the discrete vorticity evaluation at line 5 in algorithm \ref{alg:CME3D} will eventually produce aliasing errors from undersampling, we will usually introduce a low-pass filter or a Fourier truncation at lines 6 and 8. This is the effective scale cut-off of the velocity field governing the discrete flow. As a rule of thumb, we pick this low-pass filter to be the coarsest scale in the discretization and the grid $\gV$ to be the finest scale with $\gM$ at an intermediate scale. This is to ensure that the map grid has enough resolution to represent a short-time deformation generated by the filtered velocity field, and also that the $\gV$ grid is fine enough so that the Hermite interpolation of the filtered velocity is accurate and preserves the divergence-free property.
\end{remark}

\section{Numerical Tests} \label{sec:numTests}
In this section, we present several numerical tests computed using the CM method for the 3D incompressible Euler equations. The algorithm is implemented in C using OpenMP parallelization and Discrete Fourier transforms are performed using the FFTW library \cite{FFTW05}. The tests in this section are performed on a laptop computer with an AMD Ryzen 7 4800H CPU with 8 cores (16 threads) and 16GB of RAM; for these tests, a wallclock computation time is recorded. The larger FFT computations for the spectrum plots are performed on a cluster computer. The application of spatial adaptivity was not studied in this work, however, the general formulation allows for the use of adaptive grids; this has been studied for the Gradient-Augmented Level-Set methods in \cite{kolomenskiy2016adaptive} which should be straightforwardly extendable to the CM methods.

\subsection{Perturbed Antiparallel Vortex Tubes} \label{sec:numTestKerr}
The question of finite-time blow-up in the solution of the 3D incompressible Euler equations is an important open problem in mathematics. One extensively studied initial condition for potentially generating a finite-time blow-up are the perturbed antiparallel vortex tubes studied by Kerr in 1993 \cite{kerr1993evidence}. In the viscous case, for the Navier-Stokes equations, this initial condition evolves into a vortex reconnection in the process of which a topological change of the vortex cores occurs.

The initial condition can be constructed as the pullback of two antiparallel vortex tubes by a shear-deformation of the $[-2\pi, 2\pi]^3$ periodic domain. The initial vorticity field is antisymmetric across the $z = 0$ plane with each half-space containing a vortex tube of opposite orientation. We construct the initial condition $\vw_0$ as follows. Consider the unperturbed vortex tube in the $z>0$ half given by
\begin{gather} \label{eq:vortexTubeIC}
\vphi_+ (x,y,z) = \exp \left( \frac{-r^2}{1-r^2} + r^4(1+ r^2 + r^4) \right) \left( \begin{matrix}
0 \\ 1 \\ 0 
\end{matrix} \right) \quad \text{ if } \quad r < 1 ,
\end{gather}
and is the zero vector if $r \geq 1$. Here, $r$ is the scaled distance from the vortex core given by
\begin{gather}
r(x,y,z) = R^{-1}\sqrt{(x-x_0)^2 + (z-z_0)^2}.
\end{gather}
This forms a vortex tube oriented in the $y$-direction centered at $x=x_0$, $z=z_0$ and supported in a tube of radius $R$. One can check that this initial condition is divergence-free and therefore a valid vorticity field on the flat 3-torus. The pair of antiparallel vortex tubes is given by
\begin{gather}
\vphi (x,y,z) = \vphi_+ (x,y,z) - \vphi_+(x,y,-z).
\end{gather}

The vortex tubes are perturbed by the following domain deformation:
\begin{gather}
\vT : \left( \begin{matrix}
x \\ y \\ z 
\end{matrix} \right) \mapsto
\left( \begin{matrix}
x + \delta_x \cos \left( \frac{\pi}{L_x} s(y) \right) \\ y \\ z + \delta_z \cos \left( \frac{\pi}{L_z} s(y) \right)
\end{matrix} \right),
\end{gather}
where
\begin{gather}
s(y) = y + L_y \delta_{y2} \sin\left( \pi y /L_y \right) + L_y \delta_{y1} \sin \left( y + L_y \delta_{y2} \sin\left( \pi y /L_y \right) \right) .
\end{gather}

The perturbed vortex tubes are then defined as the pullback $(\vT^{-1})^* \vphi $. A closed form expression can be obtained using $(\vT^{-1})^* \vphi  = \grad \vT |_{\vT^{-1}} \vphi (\vT^{-1})$ and the fact that $\vT$ is a shear-deformation, and for a fixed $y$, $\vT$ is simply a translation on the $x$-$z$ plane.

\begin{figure}
\centering
\includegraphics[width = 0.35\linewidth]{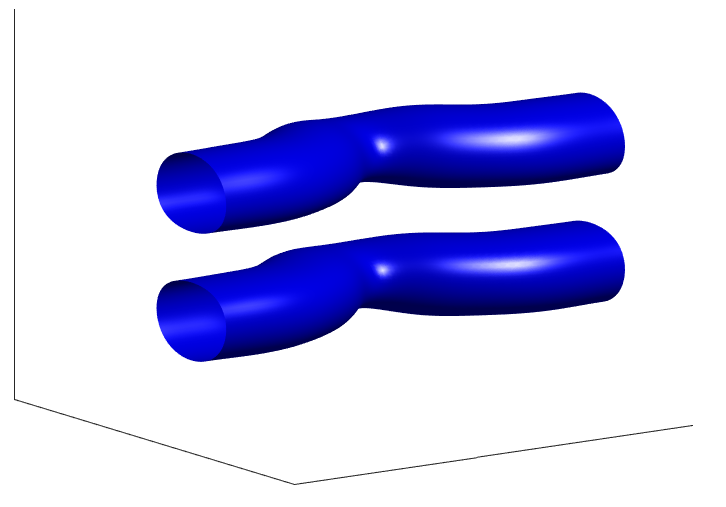}
\caption{Initial vortex tubes for the Kerr initial condition \cite{kerr1993evidence}. The figure shows the level-set surface of $|\vw_0|$ at 0.4015 which is $60\%$ of the maximum value.}
\label{fig:kerrInit}
\end{figure}

The initial vorticity $\vw_0$ used in \cite{kerr1993evidence, hou2006dynamic, yao2020physical} is defined as a filtered and rescaled version of the above perturbed vortex tubes given by
\begin{gather}
\vw_0 = 8 K * (\vT^{-1})^* \vphi .
\end{gather}
The exact expression for the filter $K$ might have been slightly different in references \cite{kerr1993evidence} and \cite{hou2006dynamic}, here we use the filter $K$ defined in Fourier space by $\hat{K} (\vxi) = \exp (-0.05(\xi_1^4 + \xi_2^4 + \xi_3^4))$ where $\xi_i$ are the integer wave numbers. The specific parameters for the initial condition, taken from \cite{hou2006dynamic}, are $R = 0.75, \, \delta_{y1} = 0.5, \, \delta_{y2} = 0.4, \, \delta_x = -1.6, \, \delta_z = 0, \, x_0 = 0, \, z_0 = 1.57, \, L_x = L_y = 4 \pi, \, L_z = 2 \pi$. Figure \ref{fig:kerrInit} shows a level-set surface of the initial condition.

\begin{figure}
\centering
\begin{subfigure}{0.32\linewidth}
\centering
\includegraphics[width = \linewidth]{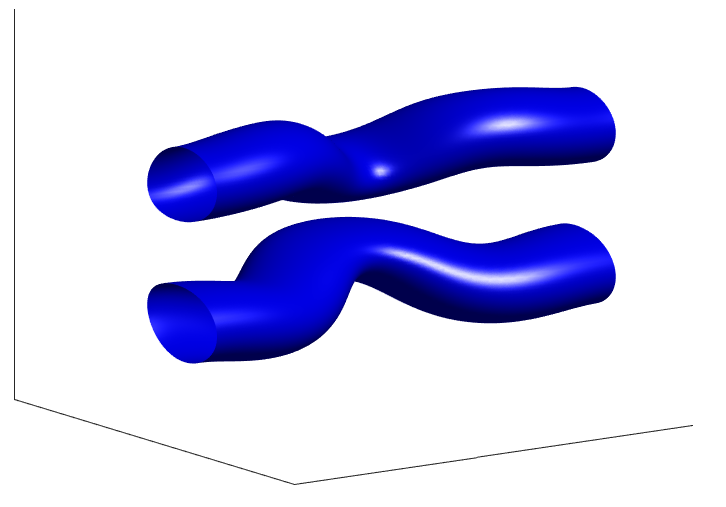}
\caption{$t=3$}
\end{subfigure}
\begin{subfigure}{0.32\linewidth}
\centering
\includegraphics[width = \linewidth]{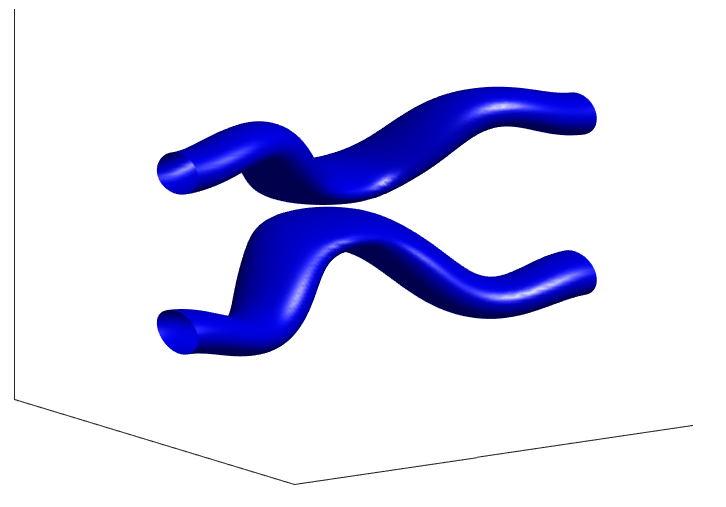}
\caption{$t=6$}
\end{subfigure}
\begin{subfigure}{0.32\linewidth}
\centering
\includegraphics[width = \linewidth]{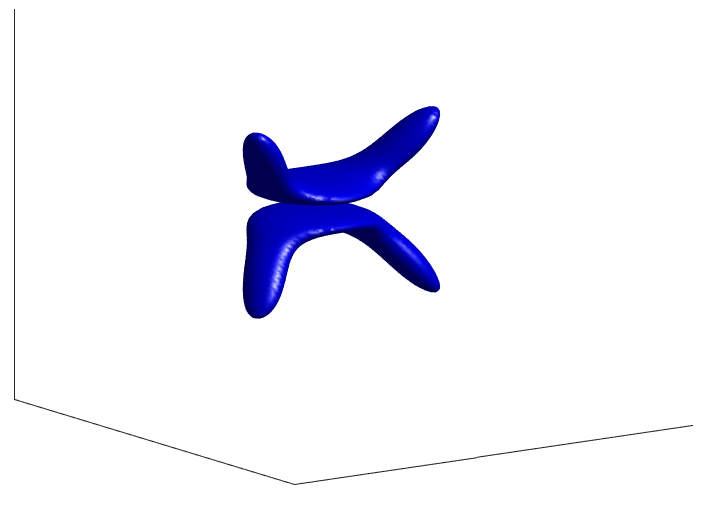}
\caption{$t=9$}
\end{subfigure}
\caption{Evolution of the vortex cores for the Kerr initial condition. We show level-set surfaces of $|\vw|$ at 0.4105, 0.5435 at 0.7519 for times 3, 6 and 9 respectively, which is $60\%$ of the maximum value. Figures are generated using a $128^3$ grid.}
\label{fig:vtxKerrEvolve}
\end{figure}

We first perform a low-resolution simulation using a laptop computer. The simulation is carried out on a $\gM$ grid of $64 \times 48 \times 32$ points and vorticity sampling is computed on a $\gV$ grid of $96 \times 72 \times 48$ points; here the $z$-direction is most finely sampled as important higher frequencies are expected to be produced in that direction, similarly, the $y$-direction has the coarsest representation. The resulting Fourier series is truncated at a radius of $32$, the time steps are fixed at $1/50$ and the Jacobian determinant error tolerance $| \det \grad \vhX_B -1|$ for the remapping is set at $10^{-3}$. The initial condition is defined on a $128^3$ grid, we then compute its mixed-partial derivatives in Fourier space in order to define a Hermite cubic interpolant. The simulation is run until time $t=17$ requiring a total of $79$ submaps, taking under an hour of wallclock computation time. Figure \ref{fig:vtxKerrEvolve} shows the evolution of the vortex cores throughout the simulation and table \ref{tab:kerrSim} contains the energy and helicity conservation errors as well as total enstrophy and maximum vorticity and velocity evaluated at regular time intervals. The energy is defined as the squared $L^2$ norm of the velocity in the $[-2 \pi, 2\pi]^3$ domain $\| \vu \|_{L^2}^2$, the total enstrophy is $\| \vw \|_{L^2}^2$ and the helicity is defined as the $L^2$ inner-product of the velocity and the vorticity, $H := (\vu, \vw)_{L^2}$.

\begin{remark}
We note here that the Fourier truncation or filtering used in the computation of the velocity field from the sampled vorticity is not related to the 2/3 rule typically used in Fourier pseudo-spectral methods, whose purpose is to dealias spurious modes generated by the frequency convolutions when computing the nonlinear term in physical space in each time step. In the CM method, there is no direct time-stepping of the velocity and vorticity fields, at each step, the vorticity is reconstructed by direct sampling of the functional expression in equation \eqref{eq:submapCompNum}. Since all map computations are carried out on a coarse grid $\gM$, the purpose of the filtering is to ensure that the velocity field defined from the sampled vorticity is sufficiently band-limited so that the backward flow map it generates is regular enough to be accurately represented on $\gM$. In the extreme cases where the vorticity field exhibits important subgrid scales, essentially discontinuous from a numerical point of view, this filtering can help prevent the Gibbs phenomenon from generating spurious oscillations in the entire domain. Ultimately, the size of the truncation would scale with the resolution of $\gM$ to maintain consistency. This filtering is not always necessary, without filtering, the effective truncation of the Fourier series for the velocity will be the grid size of $\gM$ as higher frequencies in the flow cannot be represented on $\gM$. However, our numerical experiments suggest that a small amount of smoothing generates better results with more accurate energy and helicity conservation.
\end{remark}

\begin{table}[h]
\centering
\begin{tabular}{| r || c | c | c | c | c | r | r |}
\hline
& & & & & &  & \\[-1em]
$t$ & $\| \vw \|_{L^2}^2$ & $\| \vu \|_{L^2}^2 / \| \vu_0 \|_{L^2}^{2} -1 $ & $H - H_0$ & $ \| \vw \|_{L^{\infty}} $ & $ \| \vu \|_{L^{\infty}} $ & $n_{maps}$ & time (s) \\[0.2em]
\hline
\hline
% 0 & 67.2181 & 0.000e+00 & 0.000e+00 & 0.6691 & 0.7393 & 1 & 0 \\ 
% 1 & 67.0990 & 1.341e-06 & -5.619e-14 & 0.6666 & 0.7272 & 1 & 102 \\ 
% 2 & 67.3909 & 1.744e-06 & -1.717e-13 & 0.6715 & 0.7131 & 1 & 206 \\ 
% 3 & 68.0925 & 1.989e-07 & -2.849e-13 & 0.6841 & 0.7019 & 1 & 310 \\ 
% 4 & 69.1954 & -1.406e-06 & 3.333e-13 & 0.7387 & 0.7171 & 2 & 416 \\ 
% 5 & 70.6844 & -1.663e-06 & 5.680e-13 & 0.8181 & 0.7350 & 2 & 524 \\ 
% 6 & 72.5405 & -2.369e-06 & 8.913e-13 & 0.9059 & 0.7501 & 3 & 633 \\ 
% 7 & 74.7460 & -3.141e-06 & 1.129e-12 & 1.0046 & 0.7622 & 4 & 746 \\ 
% 8 & 77.2912 & -4.061e-06 & 1.218e-12 & 1.1179 & 0.7715 & 5 & 861 \\ 
% 9 & 80.1845 & -5.020e-06 & -2.656e-12 & 1.2532 & 0.7785 & 6 & 981 \\ 
% 10 & 83.4644 & -6.245e-06 & -7.347e-12 & 1.4181 & 0.7833 & 8 & 1108 \\ 
% 11 & 87.2154 & -7.793e-06 & -4.307e-12 & 1.6296 & 0.7864 & 10 & 1240 \\ 
% 12 & 91.5883 & -1.012e-05 & -3.304e-12 & 1.9121 & 0.7885 & 13 & 1382 \\ 
% 13 & 96.8361 & -1.246e-05 & 8.511e-12 & 2.2935 & 0.7920 & 18 & 1538 \\ 
% 14 & 103.3978 & -8.680e-06 & 3.923e-11 & 2.8495 & 0.8013 & 25 & 1719 \\ 
% 15 & 112.1115 & 1.893e-05 & 4.295e-11 & 3.8549 & 0.8223 & 37 & 1935 \\ 
% 16 & 124.5780 & 9.829e-05 & -2.643e-11 & 5.6541 & 0.8545 & 56 & 2208 \\ 
% 17 & 143.8254 & 2.633e-04 & 1.910e-10 & 9.1819 & 0.9052 & 79 & 2558 \\
0 & 67.2181 & $0.000\times 10^{0}$ & $0.000\times 10^{0}$ & 0.6691 & 0.7393 & 1 & 0 \\ 
1 & 67.0990 & $1.341\times 10^{-6}$ & $-5.619\times 10^{-14}$ & 0.6666 & 0.7272 & 1 & 102 \\ 
2 & 67.3909 & $1.744\times 10^{-6}$ & $-1.717\times 10^{-13}$ & 0.6715 & 0.7131 & 1 & 206 \\ 
3 & 68.0925 & $1.989\times 10^{-7}$ & $-2.849\times 10^{-13}$ & 0.6841 & 0.7019 & 1 & 310 \\ 
4 & 69.1954 & $-1.406\times 10^{-6}$ & $3.333\times 10^{-13}$ & 0.7387 & 0.7171 & 2 & 416 \\ 
5 & 70.6844 & $-1.663\times 10^{-6}$ & $5.680\times 10^{-13}$ & 0.8181 & 0.7350 & 2 & 524 \\ 
6 & 72.5405 & $-2.369\times 10^{-6}$ & $8.913\times 10^{-13}$ & 0.9059 & 0.7501 & 3 & 633 \\ 
7 & 74.7460 & $-3.141\times 10^{-6}$ & $1.129\times 10^{-12}$ & 1.0046 & 0.7622 & 4 & 746 \\ 
8 & 77.2912 & $-4.061\times 10^{-6}$ & $1.218\times 10^{-12}$ & 1.1179 & 0.7715 & 5 & 861 \\ 
9 & 80.1845 & $-5.020\times 10^{-6}$ & $-2.656\times 10^{-12}$ & 1.2532 & 0.7785 & 6 & 981 \\ 
10 & 83.4644 & $-6.245\times 10^{-6}$ & $-7.347\times 10^{-12}$ & 1.4181 & 0.7833 & 8 & 1108 \\ 
11 & 87.2154 & $-7.793\times 10^{-6}$ & $-4.307\times 10^{-12}$ & 1.6296 & 0.7864 & 10 & 1240 \\
12 & 91.5883 & $-1.012\times 10^{-5}$ & $-3.304\times 10^{-12}$ & 1.9121 & 0.7885 & 13 & 1382 \\
13 & 96.8361 & $-1.246\times 10^{-5}$ & $8.511\times 10^{-12}$ & 2.2935 & 0.7920 & 18 & 1538 \\ 
14 & 103.3978 & $-8.680\times 10^{-6}$ & $3.923\times 10^{-11}$ & 2.8495 & 0.8013 & 25 & 1719 \\
15 & 112.1115 & $1.893\times 10^{-5}$ & $4.295\times 10^{-11}$ & 3.8549 & 0.8223 & 37 & 1935 \\ 
16 & 124.5780 & $9.829\times 10^{-5}$ & $-2.643\times 10^{-11}$ & 5.6541 & 0.8545 & 56 & 2208 \\
17 & 143.8254 & $2.633\times 10^{-4}$ & $1.910\times 10^{-10}$ & 9.1819 & 0.9052 & 79 & 2558 \\ 
\hline
\end{tabular}
\caption{Evolution of total enstrophy, energy conservation relative error (divided by initial energy), helicity conservation error (initial helicity is 0), maximum vorticity and velocity, number of remaps and wallclock computation time of the Kerr initial condition using CM method for 3D Euler. Grid resolutions: $64 \times 48 \times 32$ for $\gM$, $96 \times 72 \times 48$ for $\gV$, $\incr{t} = 1/50$, Fourier truncation at radius $32$, remapping Jacobian determinant tolerance at $10^{-3}$. All data in this table are evaluated using a grid of resolution $256^3$.}
\label{tab:kerrSim}
\end{table}

\begin{remark} \label{rmk:locatemax}
We note that the vorticity maxima $\| \vw \|_{L^{\infty}}$ shown in table \ref{tab:kerrSim} are lower bounds since they are evaluated using a $256^3$ grid. The actual maximum vorticity of the numerical solution is higher. Using the arbitrary resolution property of the method, we can refine this computation by recursively refining the vorticity sampling around the maxima. For instance, on a grid of $N^3$ points, we can locate the vorticity maximum on the grid and resample the vorticity in a region of size $3 \incr{x}$, again using a grid of $N^3$ points. With $N = 256$, 3 iterations of the above procedure allows us to estimate the vorticity maximum to 3.8618, 5.7674 and 9.5185 for times 15, 16 and 17 respectively. This shows reasonable agreement with the high resolution reference computations performed in \cite{hou2006dynamic}.
\end{remark}

\begin{figure}
\centering
\begin{subfigure}{0.485\linewidth}
\centering
\includegraphics[width = 0.45 \linewidth]{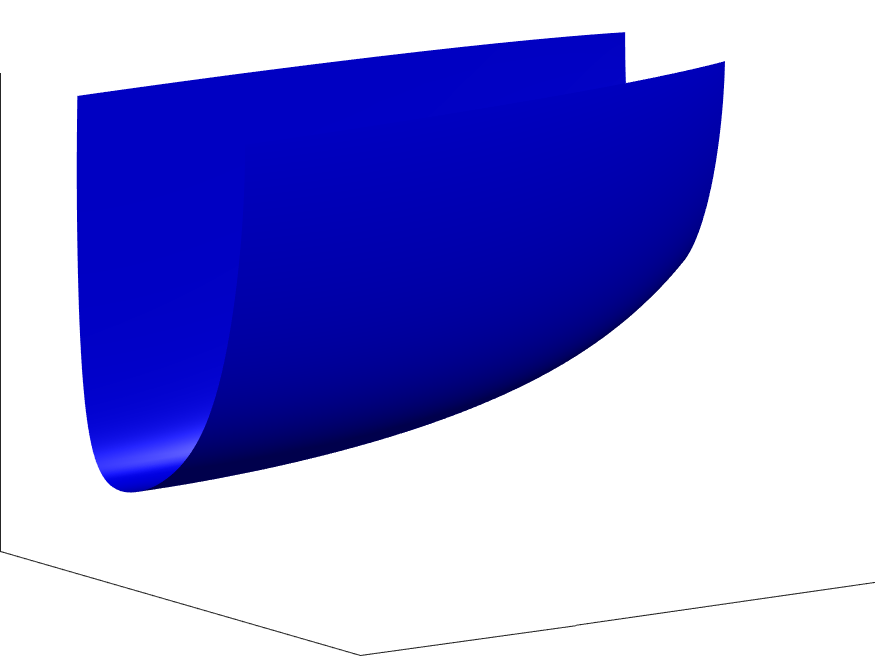}
\hspace*{5pt}
\includegraphics[width = 0.45 \linewidth]{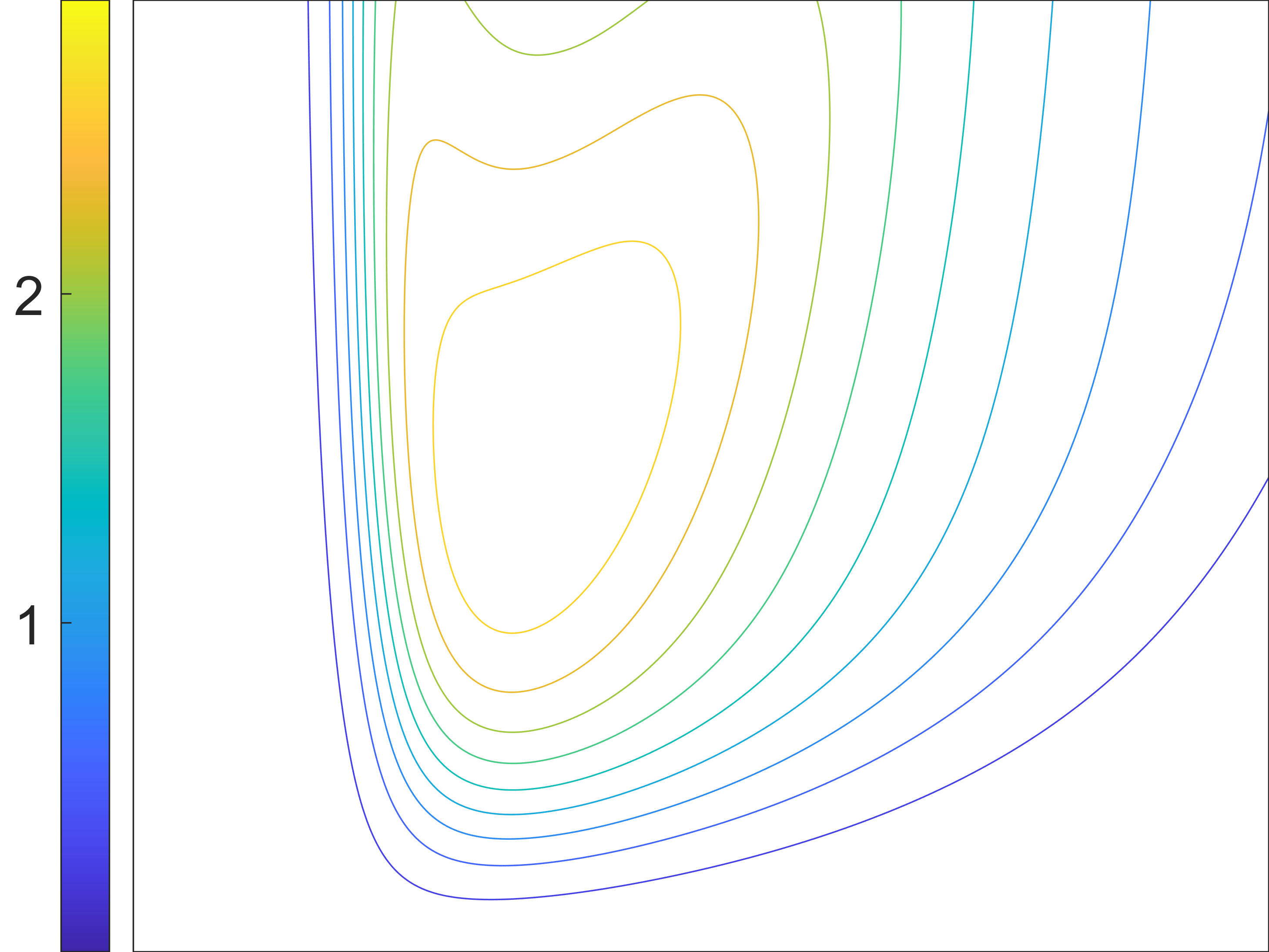}
\caption{$t=14$}
\end{subfigure}
\begin{subfigure}{0.485\linewidth}
\centering
\includegraphics[width = 0.45 \linewidth]{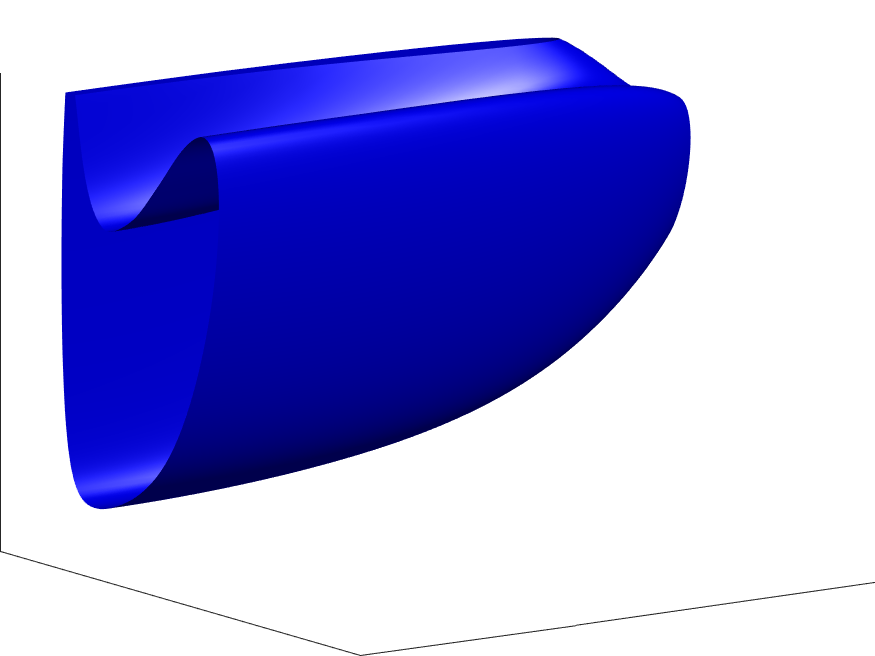}
\hspace*{5pt}
\includegraphics[width = 0.45 \linewidth]{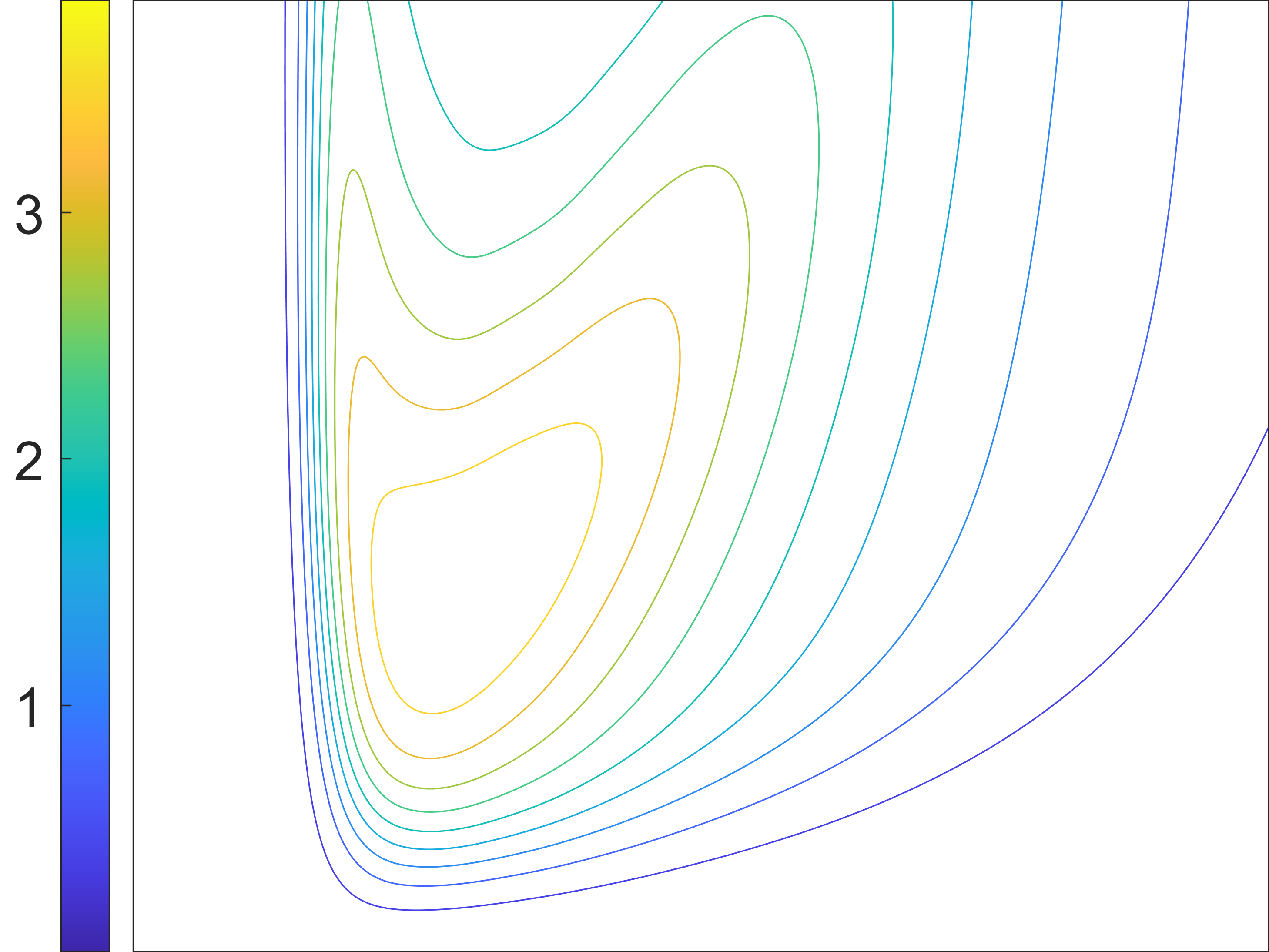}
\caption{$t=15$}
\end{subfigure}
\begin{subfigure}{0.485\linewidth}
\centering
\includegraphics[width = 0.45 \linewidth]{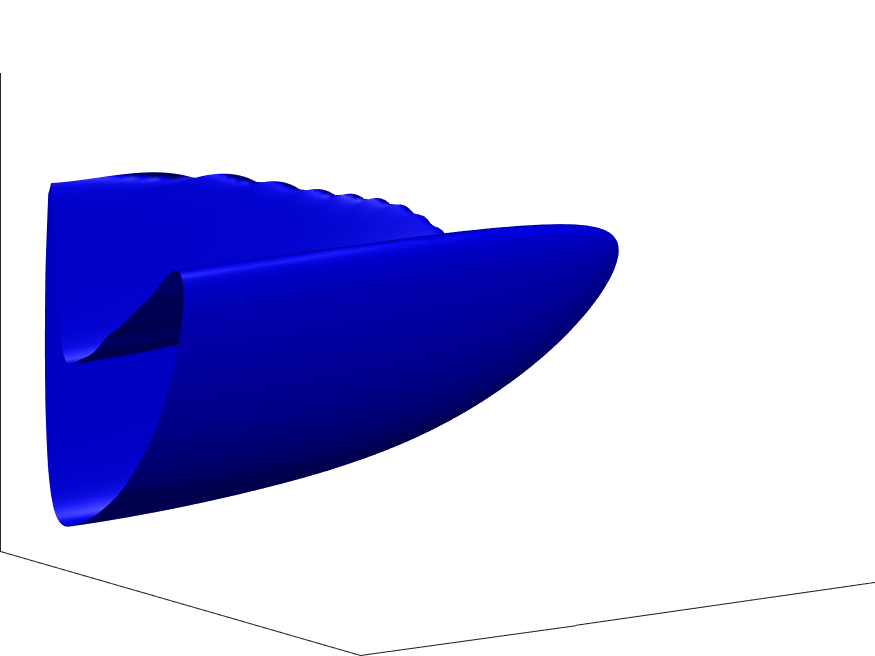}
\hspace*{5pt}
\includegraphics[width = 0.45 \linewidth]{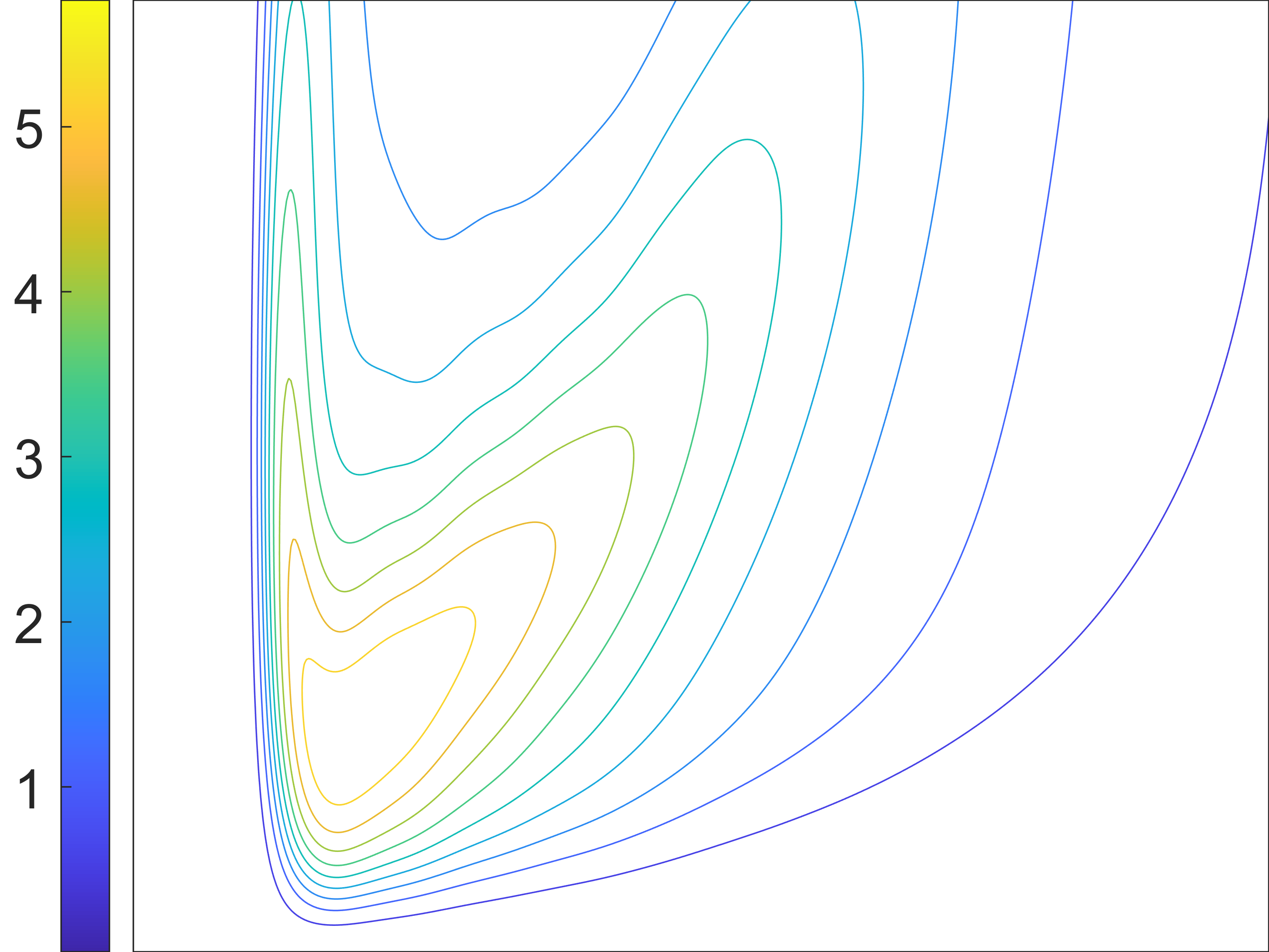}
\caption{$t=16$}
\end{subfigure}
\begin{subfigure}{0.485\linewidth}
\centering
\includegraphics[width = 0.45 \linewidth]{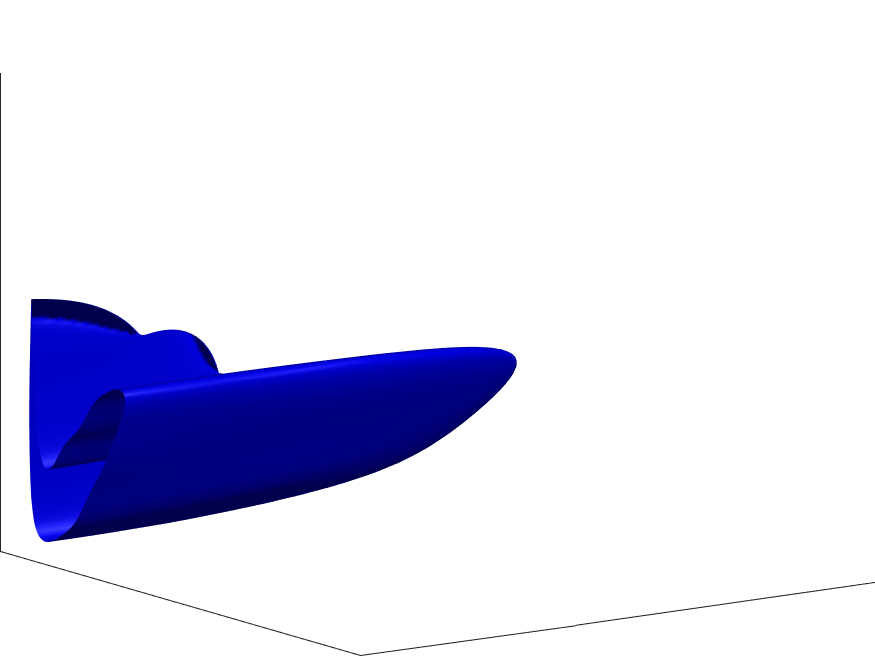}
\hspace*{5pt}
\includegraphics[width = 0.45 \linewidth]{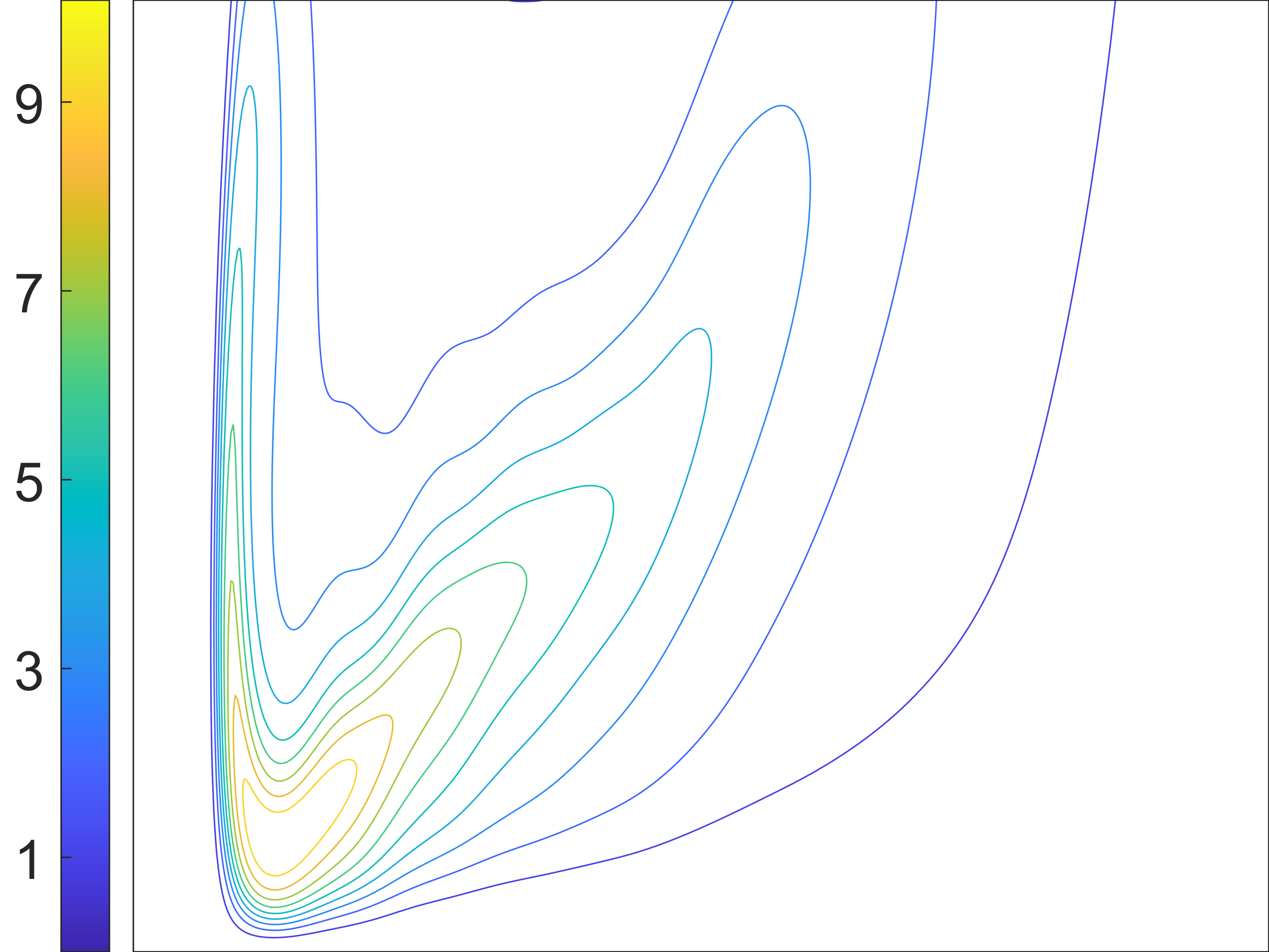}
\caption{$t=17$}
\end{subfigure}
\caption{Zoomed view of the vortex core from $t=14$ to $17$ and contour plot of the vorticity strength $|\vw|$ across the symmetry plane $y=0$; 10 isoline values evenly distributed from 0 to 2.8929, 3.8618, 5.7677 and 10.0798 are shown for times 14, 15, 16, and 17, respectively. Viewed domain: $[-5.5, 0] \times [0, 1] \times [0, 0.25]$. Contour plots are generated using $512^2$ 2D grids. We note that this is a much finer local sampling of the vorticity field compared to the data in table \ref{tab:kerrSim} and thus able to resolve a higher vorticity maximum (see remark \ref{rmk:locatemax}).}
\label{fig:vtxKerrZoom}
\end{figure}

The functional definition of the vorticity field through the pullback by $\vhX_B$ provides arbitrary resolution of the solution independently of the discretization grids. This allows us to zoom in on the solution, in particular for larger times $t$ where the vorticity starts developing significant small scale features. Figure \ref{fig:vtxKerrZoom} shows zoomed views of the vortex tubes and contour plots of the vorticity intensity across the symmetry plane $y=0$. We note that the vertical length of the viewed domain is $0.25$, a bit more than the width of a single cell of the map grid $\gM$. The domain deformation at time $t=17$ cannot be represented properly using a single map on the $\gM$ grid, however, through the dynamic remapping method, the full %time $17$ 
deformation at time $t=17$ can be represented using the composition of $79$ short-time submaps defined on a coarse $64 \times 48 \times 32$ grid. The vorticity field $\vtw$ defined by pullback  through the $79$ submaps is therefore able to exhibit small scale features and high gradients as shown in figure \ref{fig:vtxKerrZoom}.

\begin{figure}
\centering
\begin{subfigure}{0.485\linewidth}
\centering
\includegraphics[width = \linewidth]{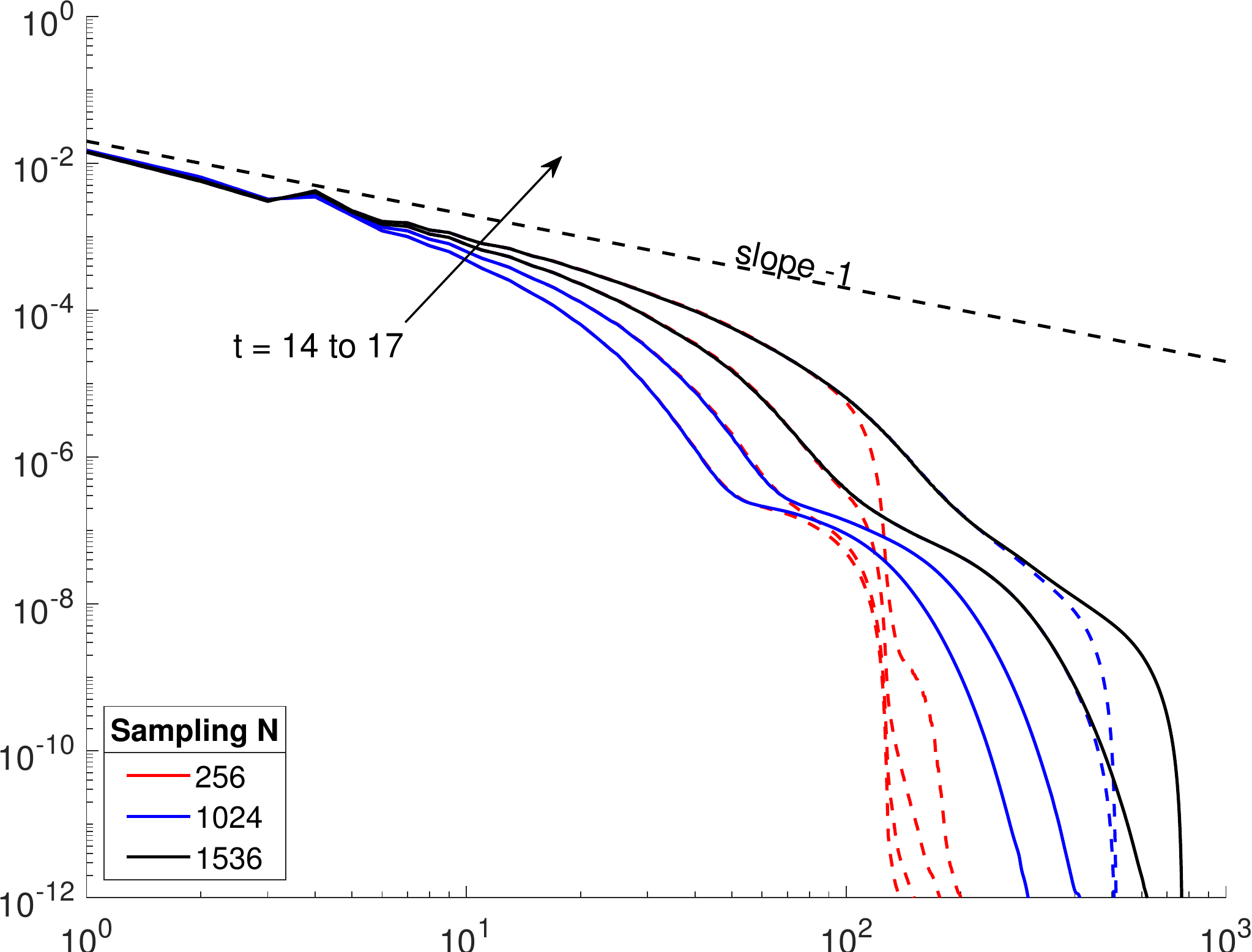}
\caption{Enstrophy}
\end{subfigure}
\begin{subfigure}{0.485\linewidth}
\centering
\includegraphics[width = \linewidth]{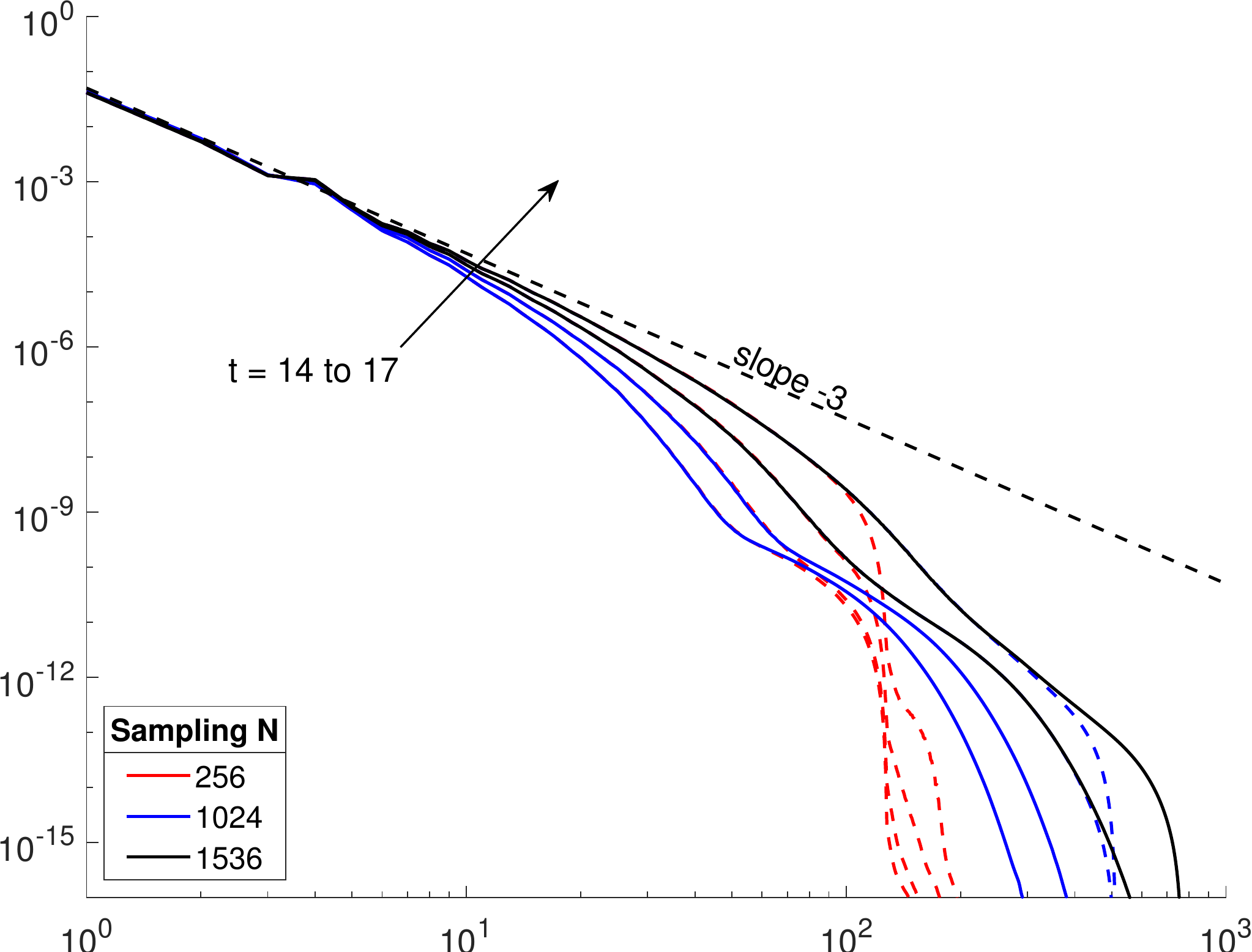}
\caption{Energy}
\end{subfigure}
\caption{Enstrophy and energy spectra at times $t= 14, 15, 16$ and $17$. Dotted black lines are the $k^{-1}$ and $k^{-3}$ curves, respectively. Successive finer sampling are shown, red curves are obtained from $256^3$ sampling, blue curves from $1024^3$ and black curves from $1536^3$ sampling. The highest resolution sampling for each time $t$ is drawn in full-line, coarser samplings for the same $t$ are shown as dotted-lines.}
\label{fig:spectra_c}
\end{figure}

As seen from table \ref{tab:kerrSim}, the remapping routine is triggered more frequently as the simulation approaches $t=17$, meaning that the Jacobian determinant error accumulates at an increasing rate in part due to a lack of spatial resolution. One often used measurement of the smoothness of the solution is the isotropic spectrum of the Fourier series given by 
\begin{subequations}
\begin{gather}
E(k) = \frac12 \sum_{|\vxi| \in [k-\frac12, k+\frac12)}  \left| \fF[\vtw](\vxi) \right|^2 \\
Z(k) = \frac12 \sum_{|\vxi| \in [k-\frac12, k+\frac12)}  \left| \fF[\vtu](\vxi) \right|^2
\end{gather}
\end{subequations}
We plot the enstrophy and energy spectra of the solution at times $t=14$ to $17$ in figure \ref{fig:spectra_c}.

The spectrum plots in figure \ref{fig:spectra_c} demonstrate a key property of the CM method: that the resolution scale of the map is not the dissipation scale of the vorticity solution, the vorticity field itself is not being dissipated in a viscous manner, and therefore can be resolved to arbitrary resolution. Indeed, although all computations were carried out on coarse grids and the Fourier support of the velocity field evolving the map is only a ball of radius 32, we can reconstruct the vorticity field by computing a fine grid pullback by evaluating the maps and applying equation \eqref{eq:submapCompNum}. Figure \ref{fig:spectra_c} shows that the vorticity and velocity fields obtained this way have the expected decay in their Fourier transforms. These can be compared with Figs. 17 and 18 in \cite{hou2006dynamic} and general agreement is found. The vorticity field $\vw^n = \left(\grad \vhX_{[t_n, 0]} \right)^{-1} \vw_0 (\vhX_{[t_n, 0]}(\vx) )$, defined functionally (i.e. $\vhX_B$ evaluated by interpolation and $\vw_0$ by direct function evaluation), contains arbitrary fine scales until round-off errors dominate. The reason for this arbitrary resolution is twofold: Firstly, since the discretized quantity is the characteristic map, the pointwise definition of the vorticity field can in fact be viewed as Lagrangian. Indeed, consider a particle starting at position $\vy$ at time $t=0$, its position at time $t$ is given by the forward map $\vX_F(\vy, t)$ and the vorticity field at the particle location is given by $\vw(\vX_F(\vy, t), t) = \grad \vX_F(\vy, t) \cdot \vw_0(\vy)$. In order to obtain the vorticity at an Eulerian point $\vx$ at time $t$, we plug in $\vy = \vX_B(\vx, t)$ and use the inverse property \eqref{eq:groupInverse} to get $\vw(\vx, t) = (\grad \vX_B)^{-1} \vw_0 (\vX_B)$. This means that by interpolating the discrete map $\vhX_B$ at an Eulerian point $\vx$, we are approximating a particle path and the associated local material deformation of a Lagrangian particle whose time $t$ position is $\vx$; the time $t$ vorticity is then directly constructed from the initial condition by applying the material deformation. This can also be expressed using the error map studied in \ref{sec:ErrorEst}, where the numerical vorticity field can be written as the pullback of the exact vorticity by an error map $\vhE_B$: $\vtw^n = {\vhE_B}^* \vw^n$. Therefore, since $\vhE_B$ is a $C^1$ diffeomorphism, as long as the error is controlled, the pullback ${\vhE_B}^* \vw^n$ will not destroy small scales. The functional definition of $\vtw^n$ by pullback therefore allows us to oversample the vorticity on a $1536^3$ grid even though all computations were carried out on much coarser grids. For traditional Eulerian methods, in order to preserve these small scales and prevent large artificial dissipation, the vorticity field will have to be discretized and evolved on a $1536^3$ grid throughout the entire computation. With the CM method, these scales are not lost to dissipation and can be obtained by a fine grid sampling. Secondly, the map error can be controlled using the submap decomposition method made possible by the group structure of the characteristic maps. Indeed, the map error $\vhE_B$ arises in part from the error in approximating \textit{SDiff}$(U)$ by a finite-dimensional interpolation space $\fV$. As the flow develops, the small scale features in $\vX_B$ not resolved in $\fV$ increase, adding to the $\vhE_B$ error. Through the remapping method, this representation error is reset to 0 for each submap since the initial condition for each new submap is the identity map, which is represented exactly in $\fV$. Appropriate remapping therefore guarantees that each submap $\vX_{[T_{i+1}, T_i]}$ can be well represented in $\fV$ and that its numerical error remains in the asymptotic regime, i.e. the omitted small scales are not significant enough to pollute the large scale, lower frequency features which carry most of the energy. The resulting global-time map $\vhX_B$ is obtained as the composition of $n_{maps}$ submaps; $\vhX_B$ can be seen as an element of $\fV^{n_{maps}}$ and therefore is able to represent the small scales features generated by the long-time flow through the composition of coarse grid maps.

\begin{figure}
\centering
\begin{subfigure}{0.24\linewidth}
\centering
\includegraphics[width = \linewidth]{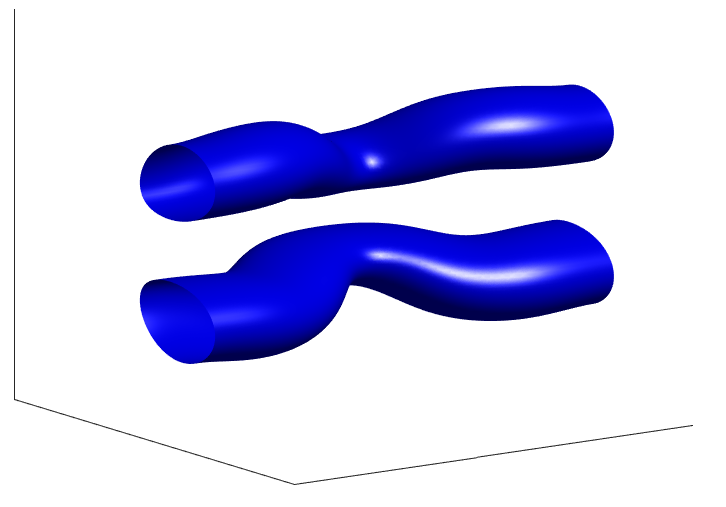}
\caption{$t=4$}
\end{subfigure}
\begin{subfigure}{0.24\linewidth}
\centering
\includegraphics[width = \linewidth]{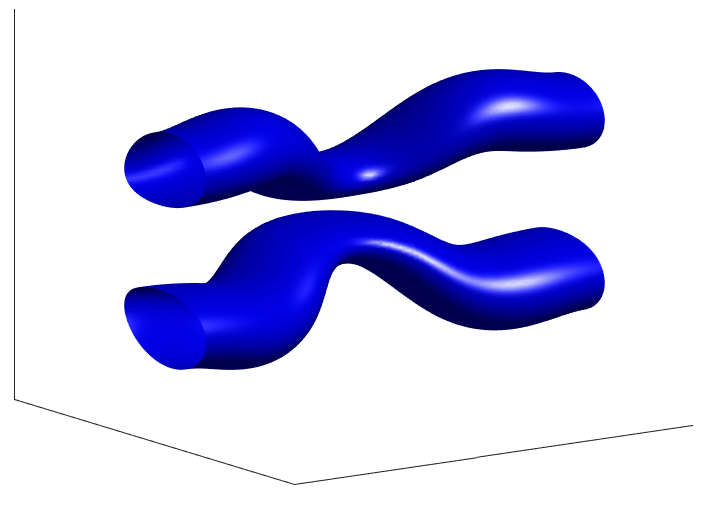}
\caption{$t=8$}
\end{subfigure}
\begin{subfigure}{0.24\linewidth}
\centering
\includegraphics[width = \linewidth]{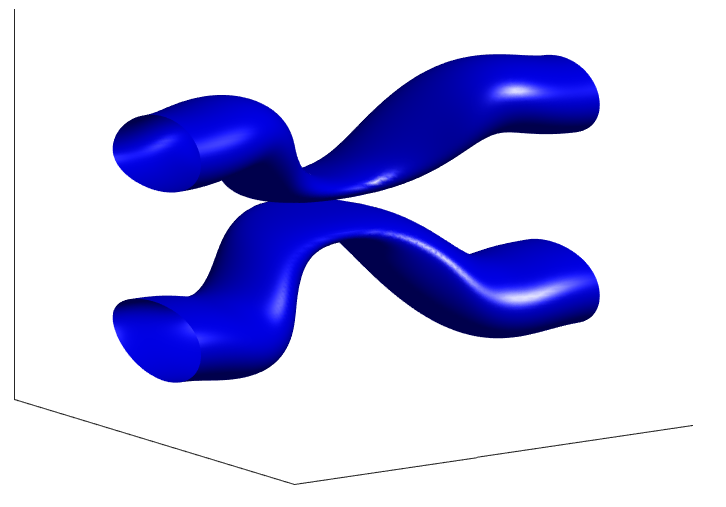}
\caption{$t=12$}
\end{subfigure}
\begin{subfigure}{0.24\linewidth}
\centering
\includegraphics[width = \linewidth]{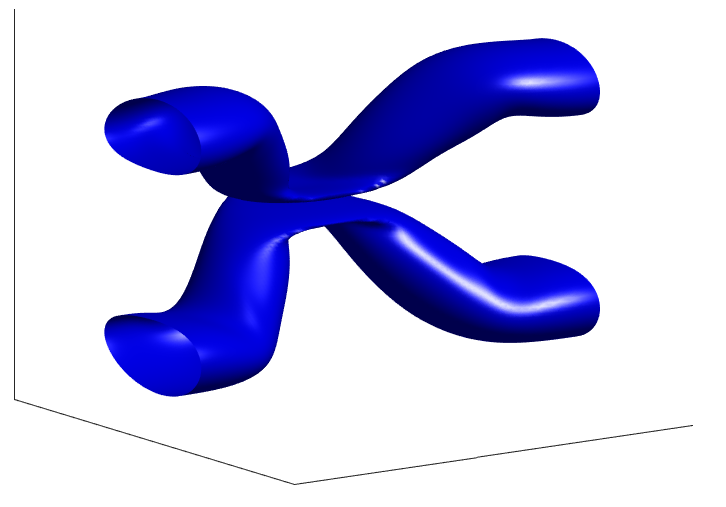}
\caption{$t=16$}
\end{subfigure}
\caption{Transport of the initial vortex core. We show level-set surfaces of $\phi$ at 0.4015, which is $60\%$ of the maximum value. Figures are generated using a $128^3$ grid.}
\label{fig:advKerrEvolve}
\end{figure}

Another feature of the CM method is that we have access to the solution operator $\vX_B$ of the advection under the velocity field $\vu$, this means that we can evolve passively advected quantities at no additional computational cost. This has several applications such as tracking passively transported fluid quantities or solute densities or visualization of the fluid flow. As example we solve the following scalar advection problem using the initial vorticity strength $| \vw_0 |$ as initial condition:
\begin{subequations} \label{eqs:vtxAdv}
\begin{gather}
(\partial_t  + \vu \cdot \grad) \phi = 0 \\
\phi(\vx, 0) = | \vw_0(\vx) | .
\end{gather}
\end{subequations}
From \ref{sec:CMM}, the solution to this advection equation is given by $\phi(\vx, t) = | \vw_0 \circ \vX_{[t, 0]} |$. This gives us the evolution of the initial vortex strength as a passively advected quantity. In figure \ref{fig:advKerrEvolve}, we show a level-set surface of $\phi$ at $60\%$ the maximum value; this allows us the track the motion of the initial vortex core transported under the fluid flow. We note that this does not correspond to the evolution of the actual vortex core as the vortex stretching can play an important role in moving the location of the vortex core.

\subsection{Perturbed Perpendicular Vortex Tubes} \label{sec:numTestCross}
\begin{figure}
\centering
\includegraphics[width = 0.45\linewidth]{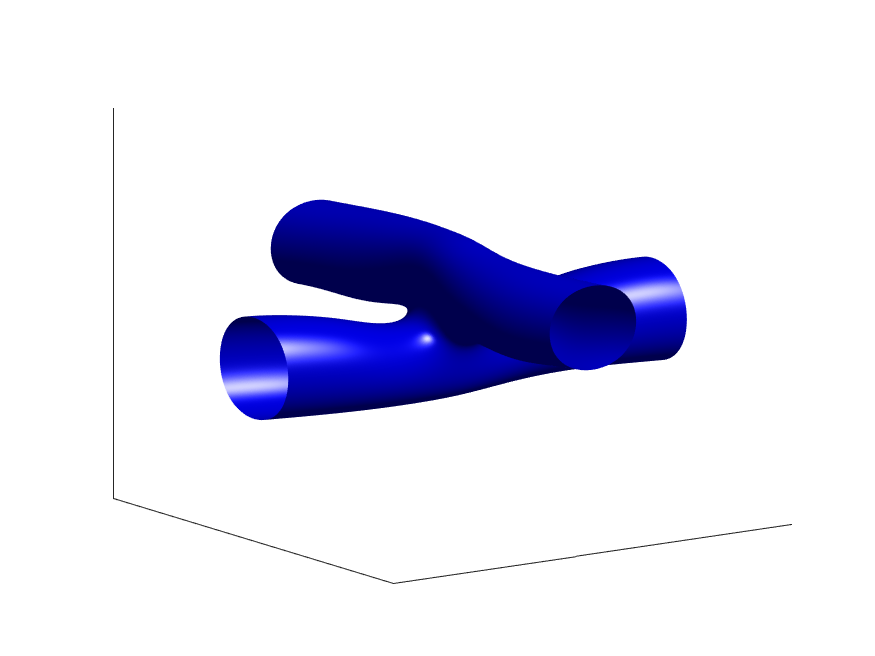}
\caption{Initial condition for the perpendicular vortex tubes test. The figure shows the level-set surface of $|\vw_0|$ at 0.5402, which is $60\%$ of the maximum value.}
\label{fig:crossInit}
\end{figure}

Another test we performed is the merging of two perturbed perpendicular vortex tubes inspired by the tests in \cite{pelz1997locally}. The initial condition is constructed in a similar fashion as in section \ref{sec:numTestKerr}. The unperturbed vortex tube is given by \eqref{eq:vortexTubeIC} with $R = 0.5, \, x_0 = 0, \, z_0 = -1$. We apply the perturbation transformations $\vT : [-2\pi, 2\pi]^3 \to [-2\pi, 2\pi]^3$ given by a sinusoidal shear deformation $\vT : x \mapsto x -0.5 \sin(0.5 y)$ and a reflection and translation map $\vR : (x,y,z) \mapsto (y, x, z+2)$. The two vortex tubes are then defined as
\begin{gather}
\vphi = (\vT^{-1})^* \vphi_+ + (\vR^{-1})^*(\vT^{-1})^* \vphi_+ ,
\end{gather}
which corresponds to a sinusoidal vortex tube in the $y$-direction through $(x,z) = (0, -1)$ combined with a reflected tube in the $x$-direction through $(y,z)  = (0, 1)$. The initial vorticity field is given as a scaled and filtered version of $\vphi$ where the filter is the same as the one used in section \ref{sec:numTestKerr},
\begin{gather}
\vw_0 = 24K*\vphi .
\end{gather}
Figure \ref{fig:crossInit} shows a level-set surface of this initial condition.

\begin{figure}
\centering
\begin{subfigure}{0.32\linewidth}
\centering
\includegraphics[width = \linewidth]{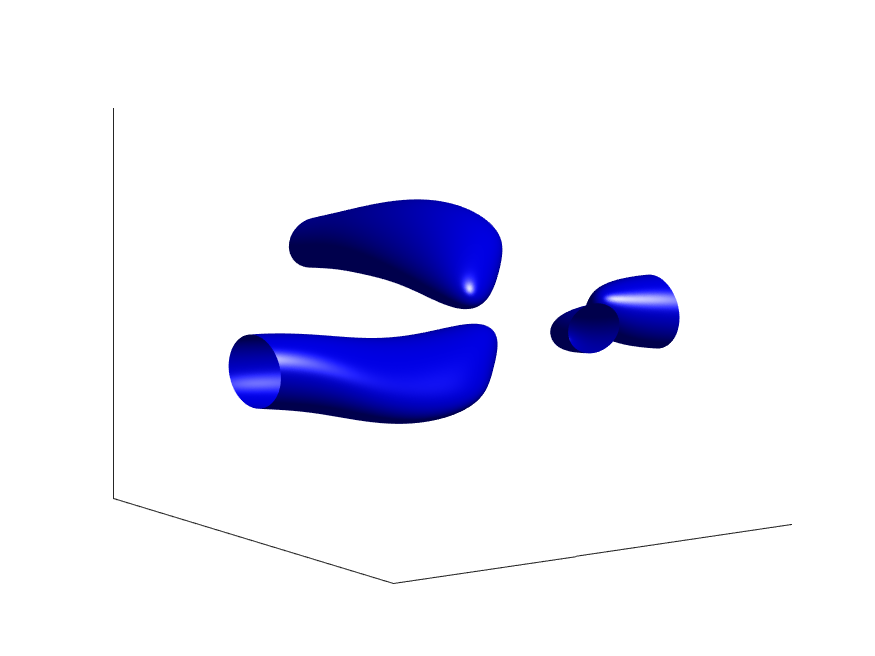}
\caption{$t=3$}
\end{subfigure}
\begin{subfigure}{0.32\linewidth}
\centering
\includegraphics[width = \linewidth]{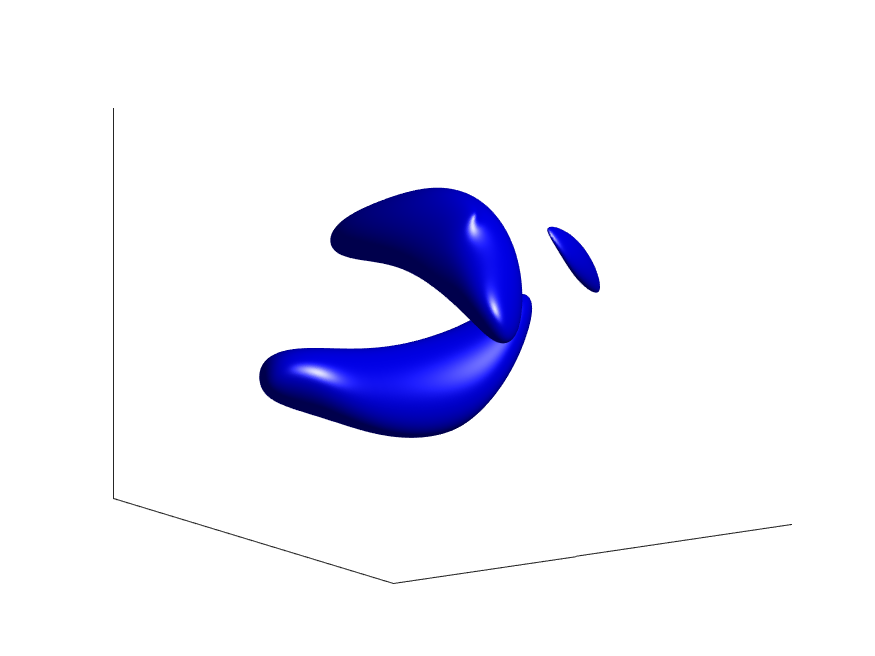}
\caption{$t=6$}
\end{subfigure}
\begin{subfigure}{0.32\linewidth}
\centering
\includegraphics[width = \linewidth]{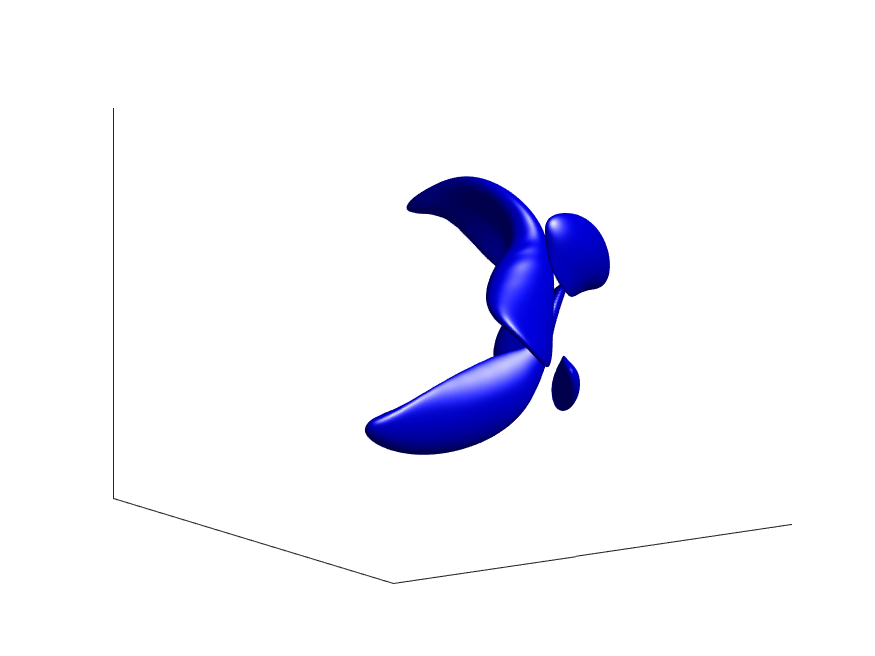}
\caption{$t=9$}
\end{subfigure}
\caption{Evolution of the vortex cores for the perpendicular vortex tubes test. We show the level-set surfaces of $|\vw|$ at 0.6817, 0.8192 and 1.2450 for times 3, 6 and 9 respectively which is $60\%$ of the maximum value. Figures are generated using a $324^3$ grid.}
\label{fig:vtxCrossEvolve}
\end{figure}

\begin{figure}
\centering
\begin{subfigure}{0.32\linewidth}
\centering
\includegraphics[width = \linewidth]{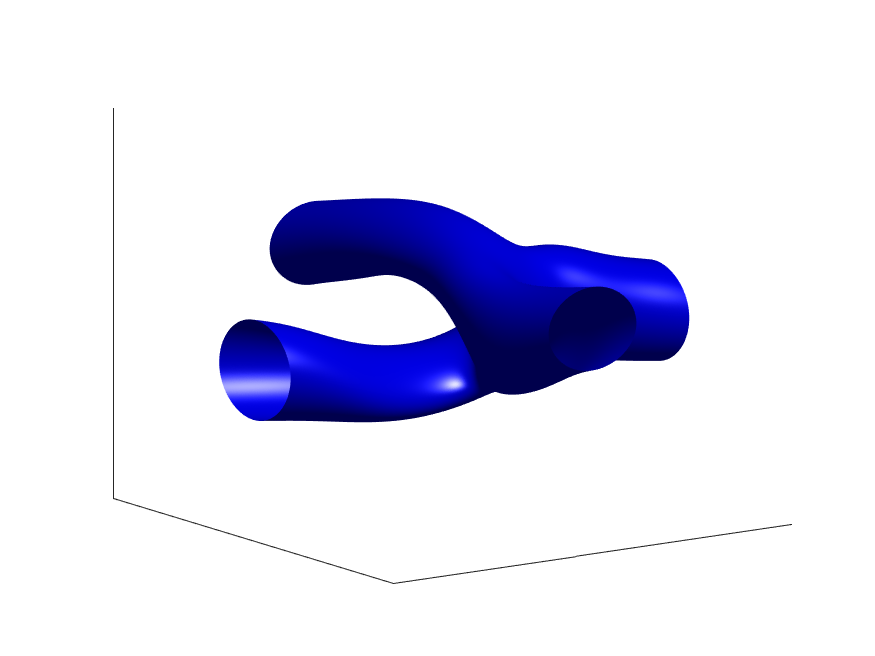}
\caption{$t=3$}
\end{subfigure}
\begin{subfigure}{0.32\linewidth}
\centering
\includegraphics[width = \linewidth]{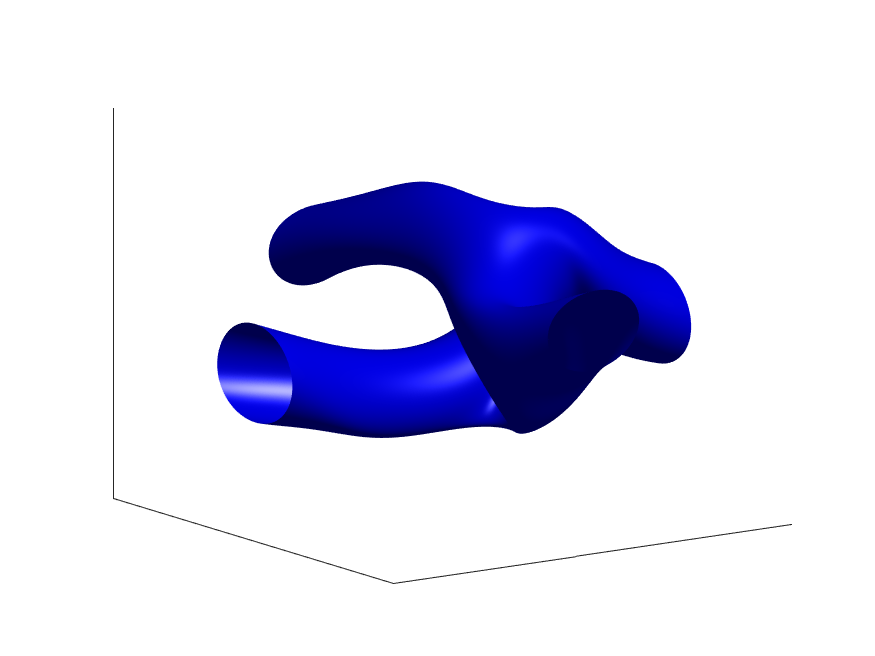}
\caption{$t=6$}
\end{subfigure}
\begin{subfigure}{0.32\linewidth}
\centering
\includegraphics[width = \linewidth]{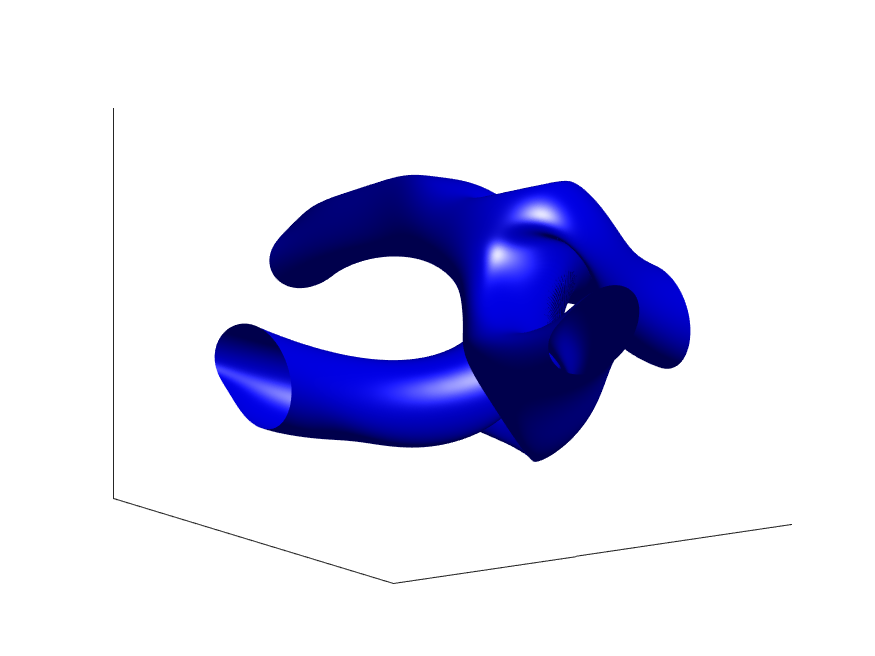}
\caption{$t=9$}
\end{subfigure}
\caption{Transport of the initial vortex core. We show level-set surfaces of $\phi$ at 0.5402, which is $60\%$ of the maximum value. Figures are generated using a $324^3$ grid.}
\label{fig:advCrossEvolve}
\end{figure}

These two vortex tubes will start rotating around each other, creating high vorticity gradients and significant small scale features. In the viscous case, the vortices are expected to collide at the intersection. Figure \ref{fig:vtxCrossEvolve} shows the evolution of the vortex tubes computed from a simulation using a $48^3$ grid for both $\gM$ and $\gV$, the Hermite cubic interpolation of the velocity field uses a $96^3$ grid to ensure sufficient smoothness. The time step $\incr{t}$ is fixed at $1/50$, Fourier truncation radius at $32$ and the Jacobian determinant error tolerance at $10^{-3}$. Due to the significant localized small scale features, the vorticity computation uses the adaptive sampling described in the appendix \ref{append:vtxSmp}. The mollifier $\mu_h$ is given by a tensor product of a 1D function $\cos^2 \left(\frac{\pi s}{2 h} \right)$ in each cell. The cell integral \eqref{eq:mollIntCell} is computed by numerical quadrature using equidistributed sample points in each cell. The number of sample points per cell is at minimum 2 points per dimension and for each cell, this number is allowed to increases adaptively depending on the range and total variation of $\vw$ in each cell. The total number of sample points is capped at $192^3$ at which point all cells have their sample number rescaled down proportionally. The results of this test until time $t=12$ are presented in table \ref{tab:crossSim}. Figure \ref{fig:vtxCrossEvolve} shows the evolution of the vortex cores until time $t=9$. To better visualize the flow, we also include in figure \ref{fig:advCrossEvolve} the scalar advection of the initial vorticity strength $\phi_0 = | \vw_0 |$ as given by equation \eqref{eqs:vtxAdv}.

\begin{table}[h]
\centering
\begin{tabular}{| r || c | c | c | c | c | r | r |}
\hline
& & & & & &  & \\[-1em]
$t$ & $\| \vw \|_{L^2}^2$ & $\| \vu \|_{L^2}^2 / \| \vu_0 \|_{L^2}^{2} -1 $ & $H /H_0 - 1$ & $ \| \vw \|_{L^{\infty}} $ & $ \| \vu \|_{L^{\infty}} $ & $n_{maps}$ & time (s) \\[0.2em]
\hline
\hline
% 0 & 125.7910 & 0.000e+00 & 0.000e+00 & 0.9004 & 0.9684 & 1 & 0 \\ 
% 1 & 126.9290 & -7.816e-04 & 2.296e-03 & 0.9183 & 0.9644 & 1 & 187 \\ 
% 2 & 130.3522 & -3.043e-03 & 8.889e-03 & 1.0317 & 0.9478 & 1 & 381 \\ 
% 3 & 136.0808 & -6.604e-03 & 1.901e-02 & 1.1361 & 0.9331 & 2 & 586 \\ 
% 4 & 144.1831 & -1.117e-02 & 3.160e-02 & 1.2187 & 0.9431 & 3 & 801 \\ 
% 5 & 154.8535 & -1.636e-02 & 4.557e-02 & 1.2885 & 0.9794 & 4 & 1027 \\ 
% 6 & 168.4914 & -2.178e-02 & 5.988e-02 & 1.3653 & 1.0208 & 8 & 1568 \\ 
% 7 & 185.7193 & -2.707e-02 & 7.354e-02 & 1.4806 & 1.0495 & 14 & 2386 \\ 
% 8 & 207.4032 & -3.201e-02 & 8.561e-02 & 1.6472 & 1.0636 & 22 & 3380 \\ 
% 9 & 234.9474 & -3.648e-02 & 9.518e-02 & 2.0750 & 1.0697 & 30 & 4691 \\ 
% 10 & 271.2813 & -4.048e-02 & 1.015e-01 & 3.3802 & 1.0602 & 40 & 6442 \\ 
% 11 & 322.7244 & -4.403e-02 & 1.038e-01 & 6.7721 & 1.0433 & 51 & 8489 \\ 
% 12 & 400.7262 & -4.733e-02 & 1.024e-01 & 14.1254 & 1.0929 & 64 & 10391 \\  
0 & 125.7910 & $0.000\times 10^{0}$ & $0.000\times 10^{0}$ & 0.9004 & 0.9684 & 1 & 0 \\ 
1 & 126.9290 & $-7.816\times 10^{-4}$ & $2.296\times 10^{-3}$ & 0.9183 & 0.9644 & 1 & 187 \\ 
2 & 130.3522 & $-3.043\times 10^{-3}$ & $8.889\times 10^{-3}$ & 1.0317 & 0.9478 & 1 & 381 \\ 
3 & 136.0808 & $-6.604\times 10^{-3}$ & $1.901\times 10^{-2}$ & 1.1361 & 0.9331 & 2 & 586 \\ 
4 & 144.1831 & $-1.117\times 10^{-2}$ & $3.160\times 10^{-2}$ & 1.2187 & 0.9431 & 3 & 801 \\ 
5 & 154.8535 & $-1.636\times 10^{-2}$ & $4.557\times 10^{-2}$ & 1.2885 & 0.9794 & 4 & 1027 \\ 
6 & 168.4914 & $-2.178\times 10^{-2}$ & $5.988\times 10^{-2}$ & 1.3653 & 1.0208 & 8 & 1568 \\ 
7 & 185.7193 & $-2.707\times 10^{-2}$ & $7.354\times 10^{-2}$ & 1.4806 & 1.0495 & 14 & 2386 \\ 
8 & 207.4032 & $-3.201\times 10^{-2}$ & $8.561\times 10^{-2}$ & 1.6472 & 1.0636 & 22 & 3380 \\ 
9 & 234.9474 & $-3.648\times 10^{-2}$ & $9.518\times 10^{-2}$ & 2.0750 & 1.0697 & 30 & 4691 \\ 
10 & 271.2813 & $-4.048\times 10^{-2}$ & $1.015\times 10^{-1}$ & 3.3802 & 1.0602 & 40 & 6442 \\ 
11 & 322.7244 & $-4.403\times 10^{-2}$ & $1.038\times 10^{-1}$ & 6.7721 & 1.0433 & 51 & 8489 \\ 
12 & 400.7262 & $-4.733\times 10^{-2}$ & $1.024\times 10^{-1}$ & 14.1254 & 1.0929 & 64 & 10391 \\ 
\hline
\end{tabular}
\caption{Evolution of total enstrophy, energy conservation relative error (divided by initial energy), helicity conservation relative error (divided by initial helicity), maximum vorticity and velocity, number of remaps and wallclock computation time of the perturbed perpendicular vortex tubes initial condition using CM method for 3D Euler. Grid resolutions: $48^3$ for $\gM$, $48^3$ for $\gV$, $\incr{t} = 1/50$, adaptive sampling with mollifier convolution is used for the vorticity, Fourier truncation at radius $32$, remapping Jacobian determinant tolerance at $10^{-3}$. All data in this table are evaluated using a grid of resolution $256^3$.}
\label{tab:crossSim}
\end{table}

\begin{figure}
\centering
\begin{subfigure}{0.8\linewidth}
\centering
\includegraphics[width = 0.32 \linewidth]{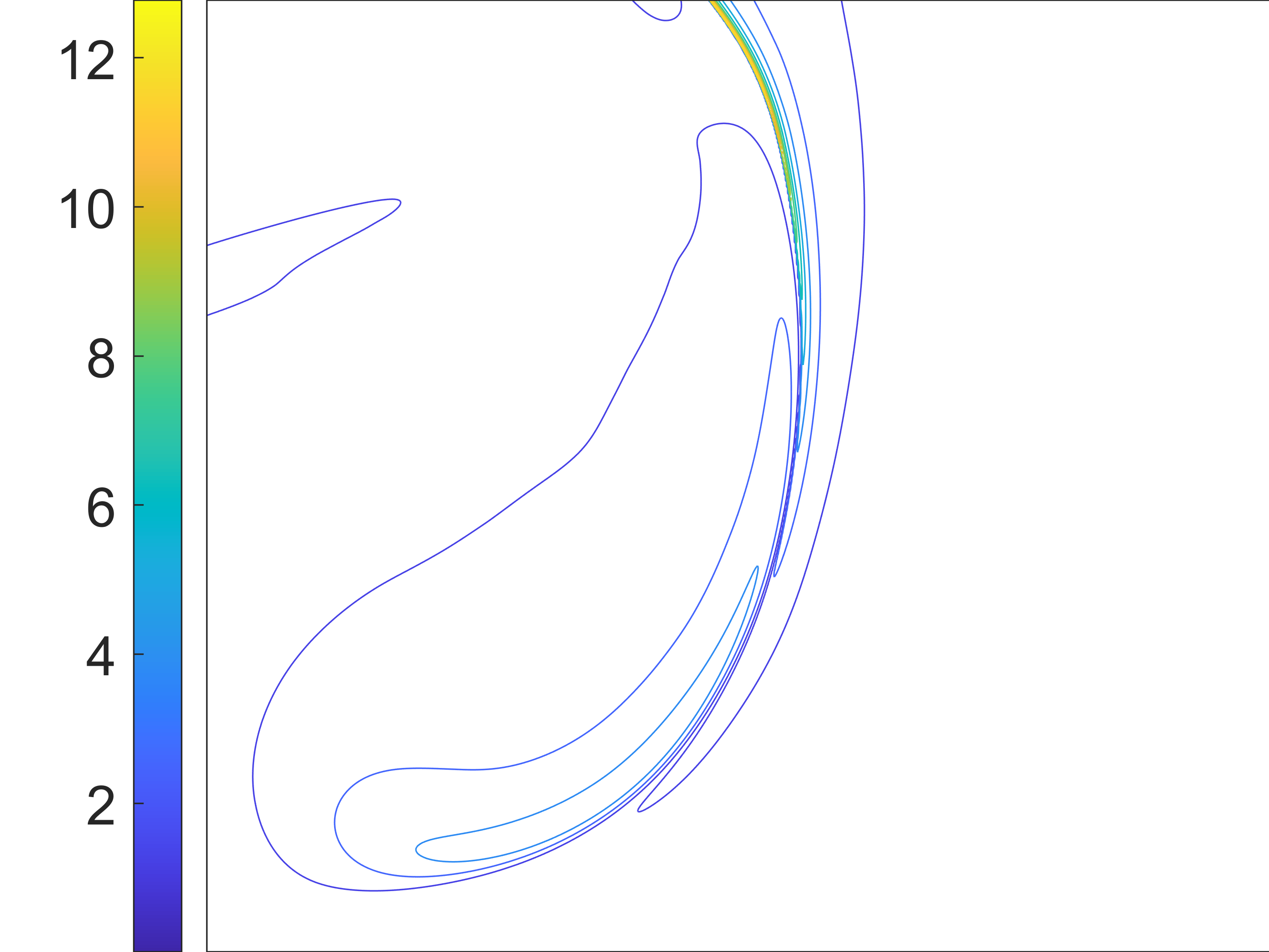}
\includegraphics[width = 0.32 \linewidth]{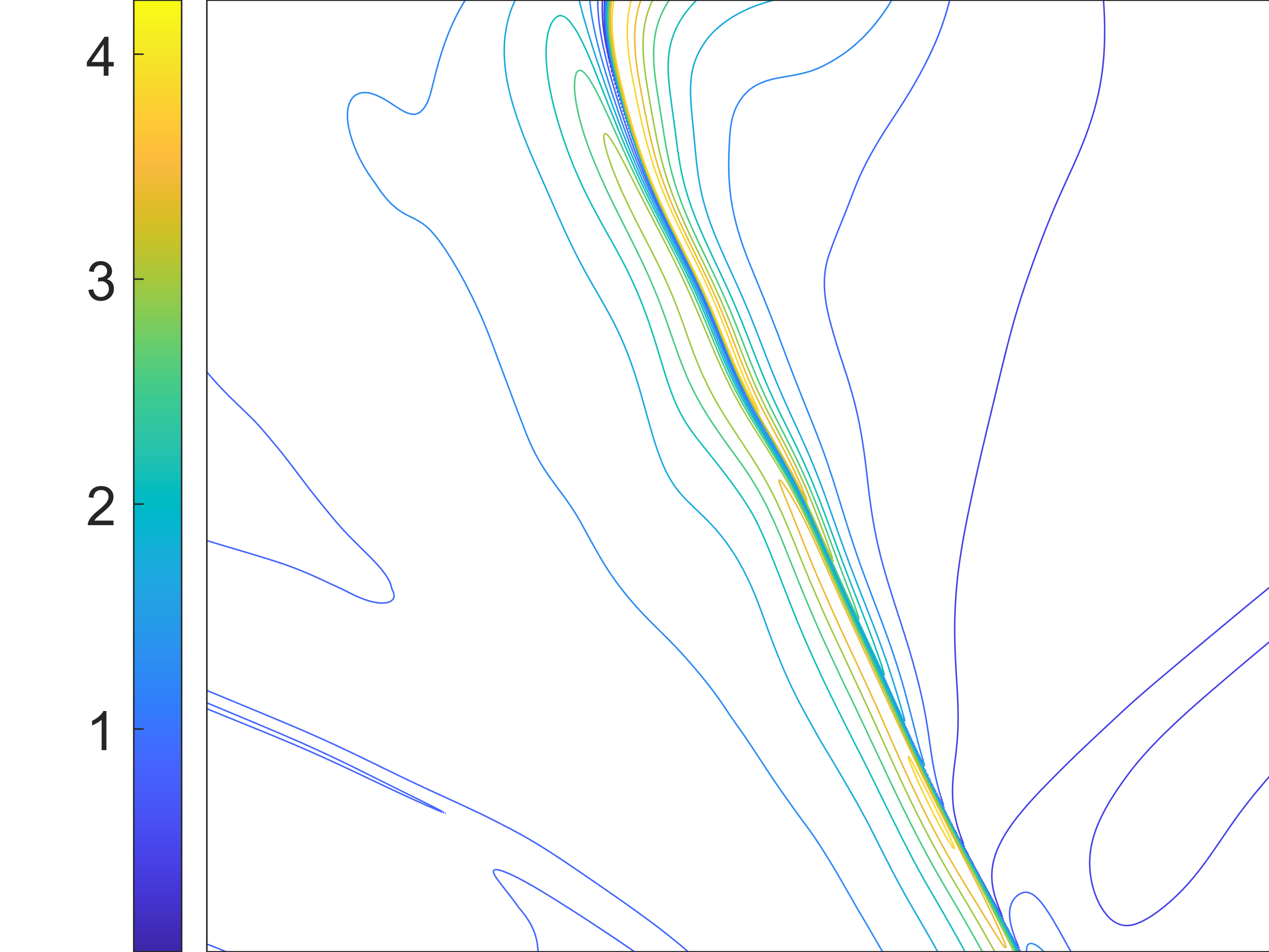}
\includegraphics[width = 0.32 \linewidth]{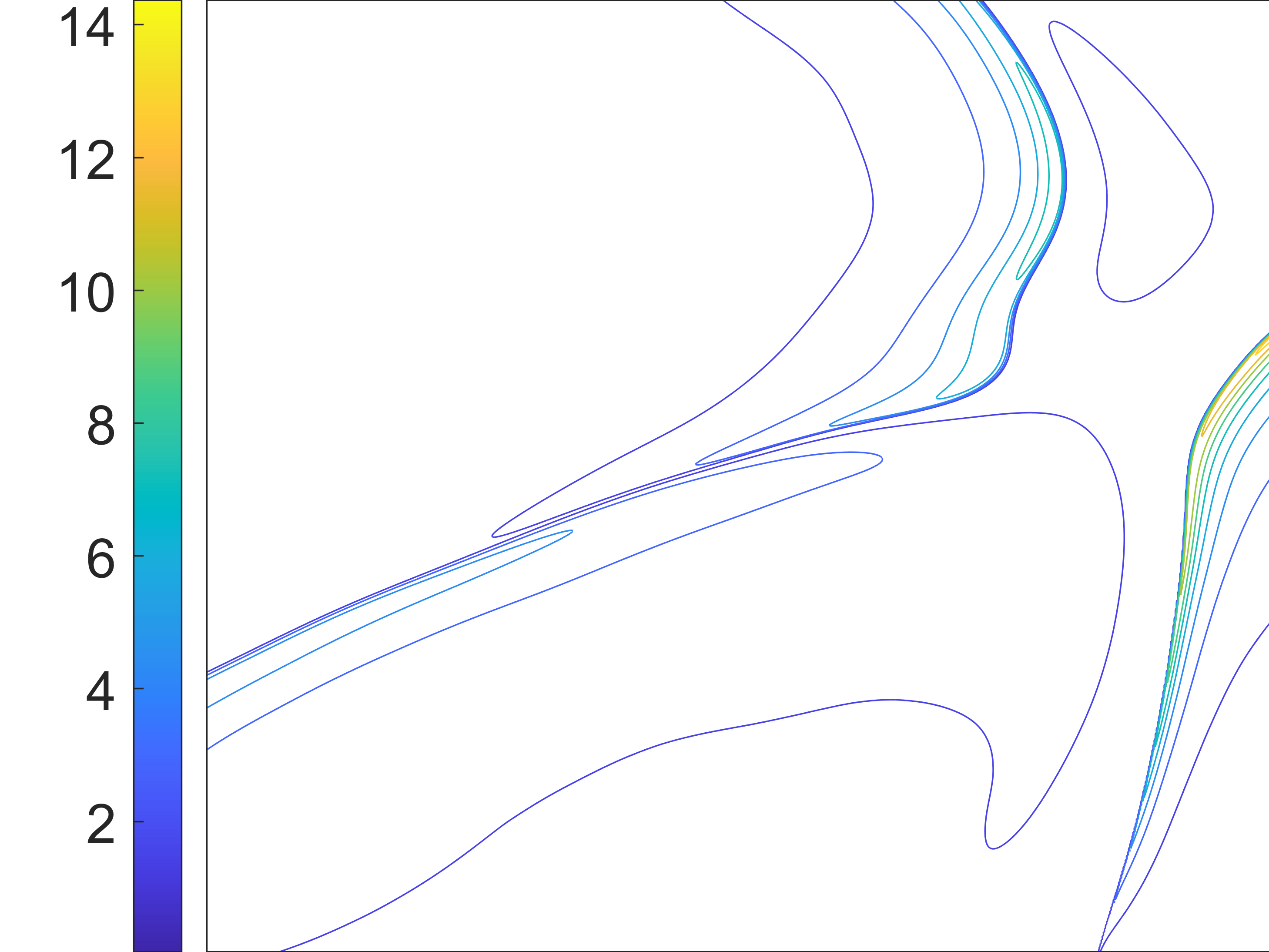}
\caption{Box width $ = \pi$.}
\end{subfigure}
\begin{subfigure}{0.8\linewidth}
\centering
\includegraphics[width = 0.32 \linewidth]{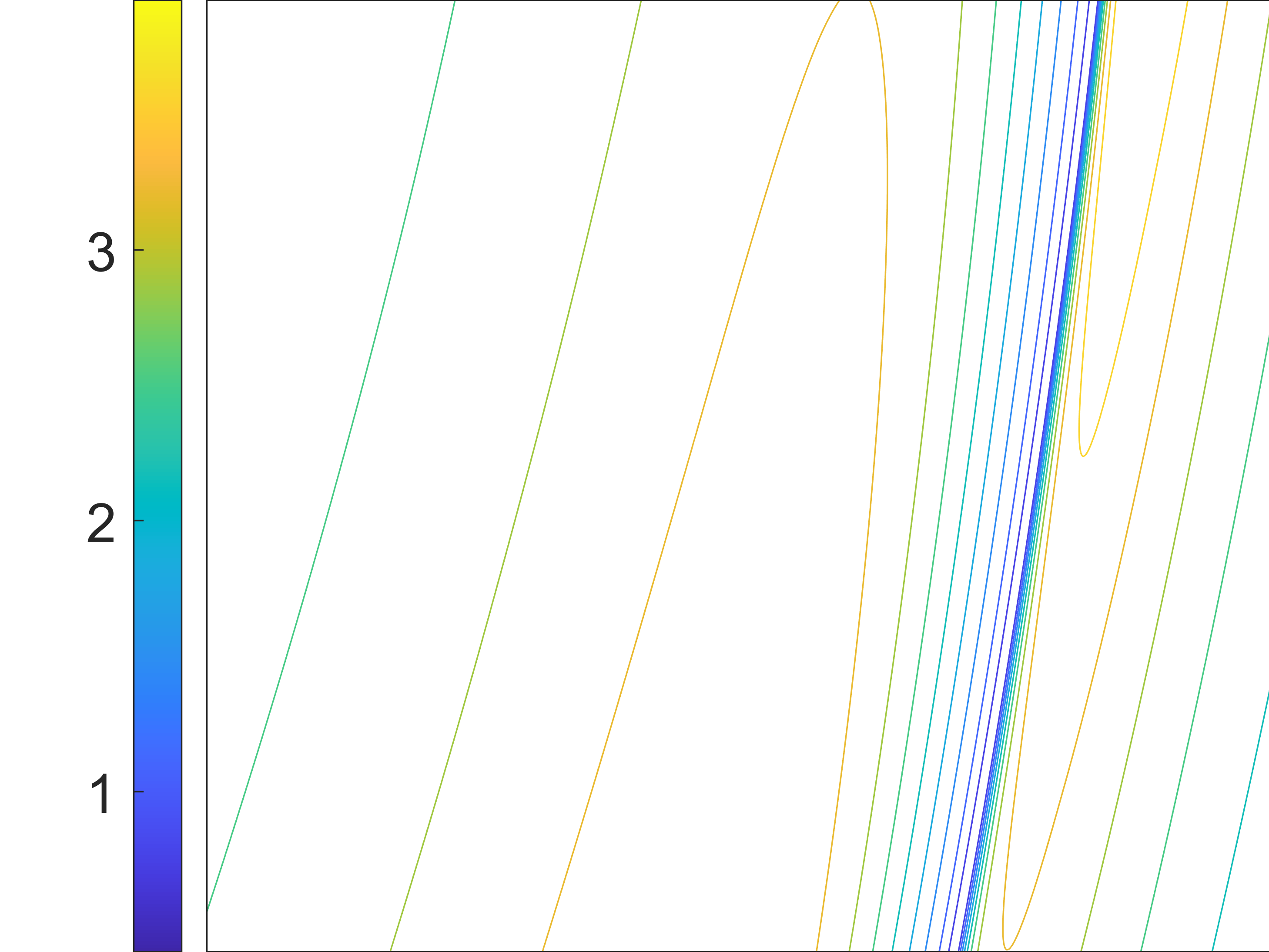}
\includegraphics[width = 0.32 \linewidth]{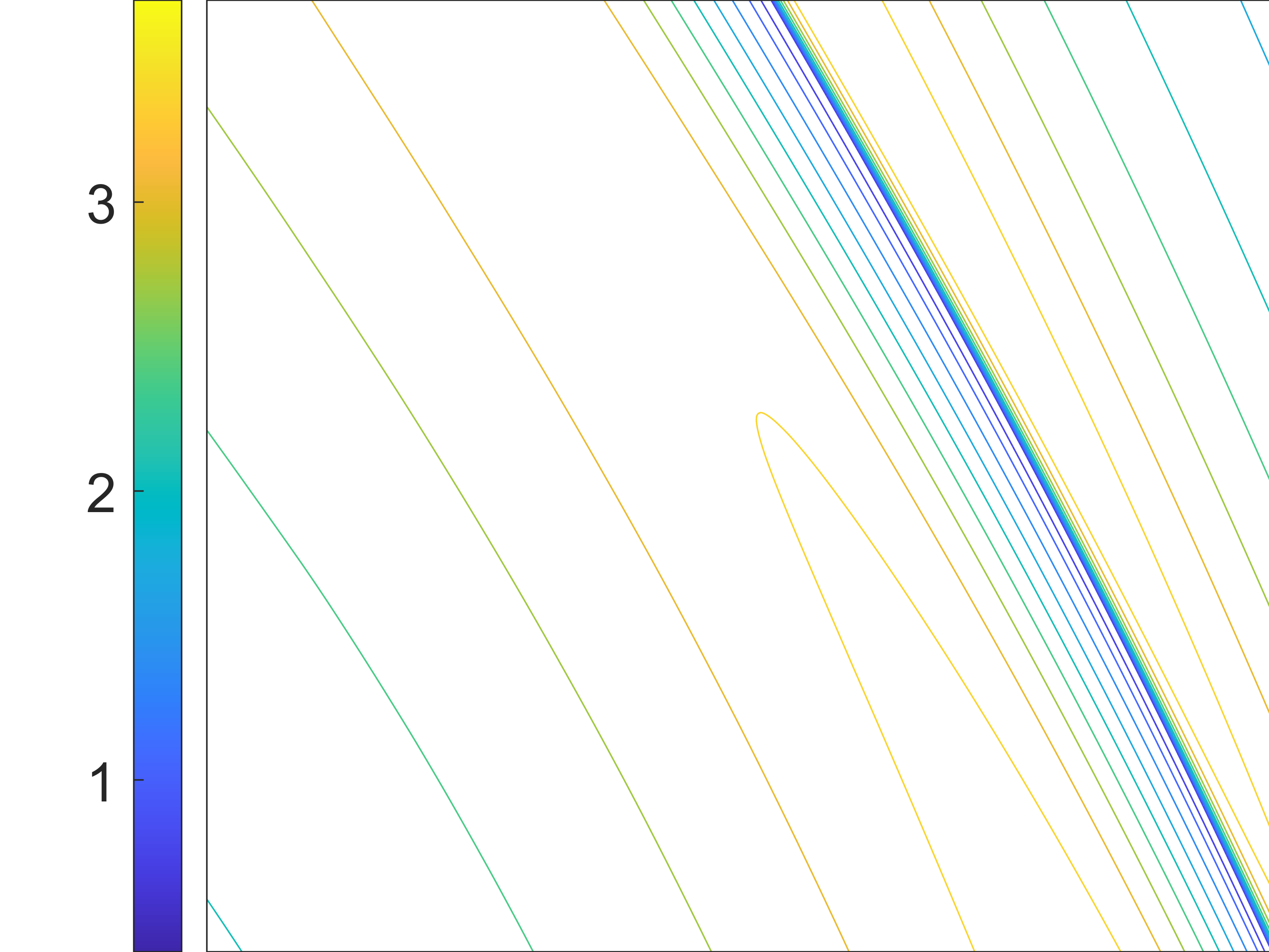}
\includegraphics[width = 0.32 \linewidth]{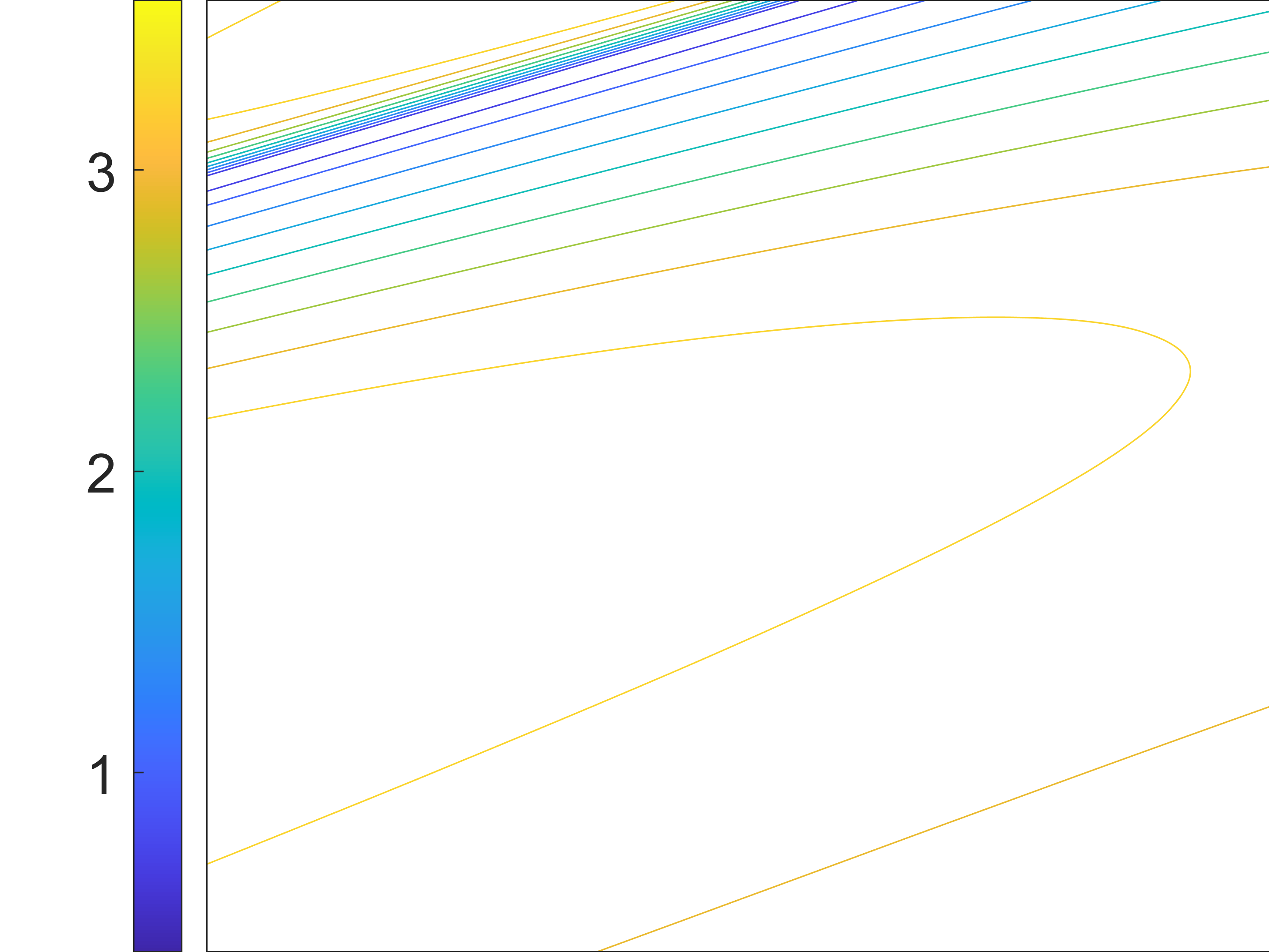}
\caption{Box width $ = \frac{\pi}{16}$.}
\end{subfigure}
\caption{Contour plot of vorticity strength $|\vw|$ at time $t=12$, centered at $(x,y,z) = (2.15, 1.5, 1.425)$; 10 isoline values evenly distributed in the range of each frame are shown. From left to right, contours on: $xy$- plane, $xz$-plane and $yz$-plane. The isoline ranges are $[0.0047, 12.7748]$, $[0.0075, 4.2412]$, $[0.0351, 14.3722]$, $[0.4072, 3.9227]$, $[0.4049, 3.6989]$ and $[0.4049, 3.5640]$. Figures are produced using a $1024^2$ 2D grid. We note that this is a much finer local sampling of the vorticity field compared to the data in table \ref{tab:crossSim} and thus able to resolve a higher vorticity maximum. }
\label{fig:vtxCrossZoom}
\end{figure}

\begin{figure}
\centering
\begin{subfigure}{0.8\linewidth}
\centering
\includegraphics[width = 0.32 \linewidth]{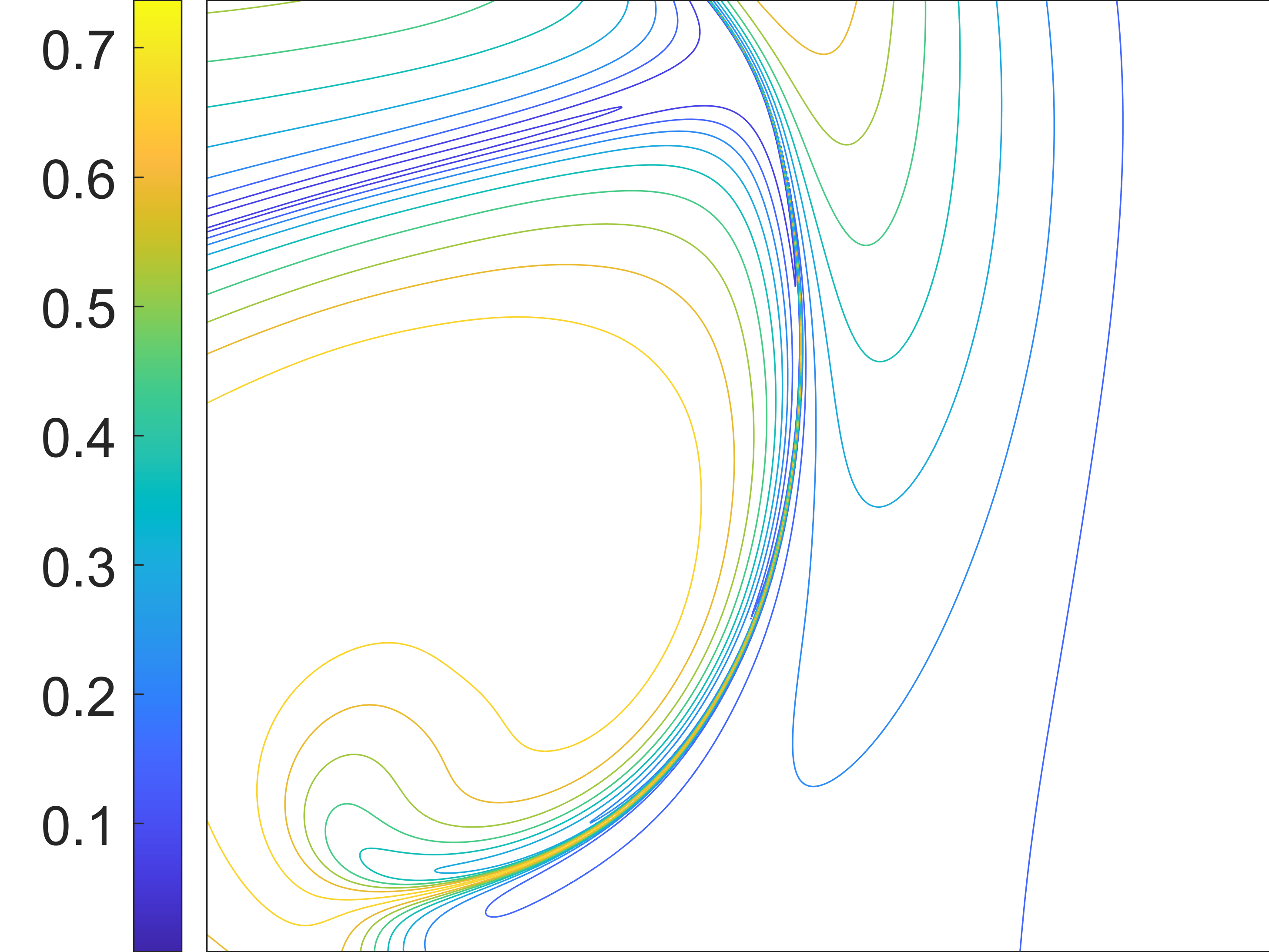}
\includegraphics[width = 0.32 \linewidth]{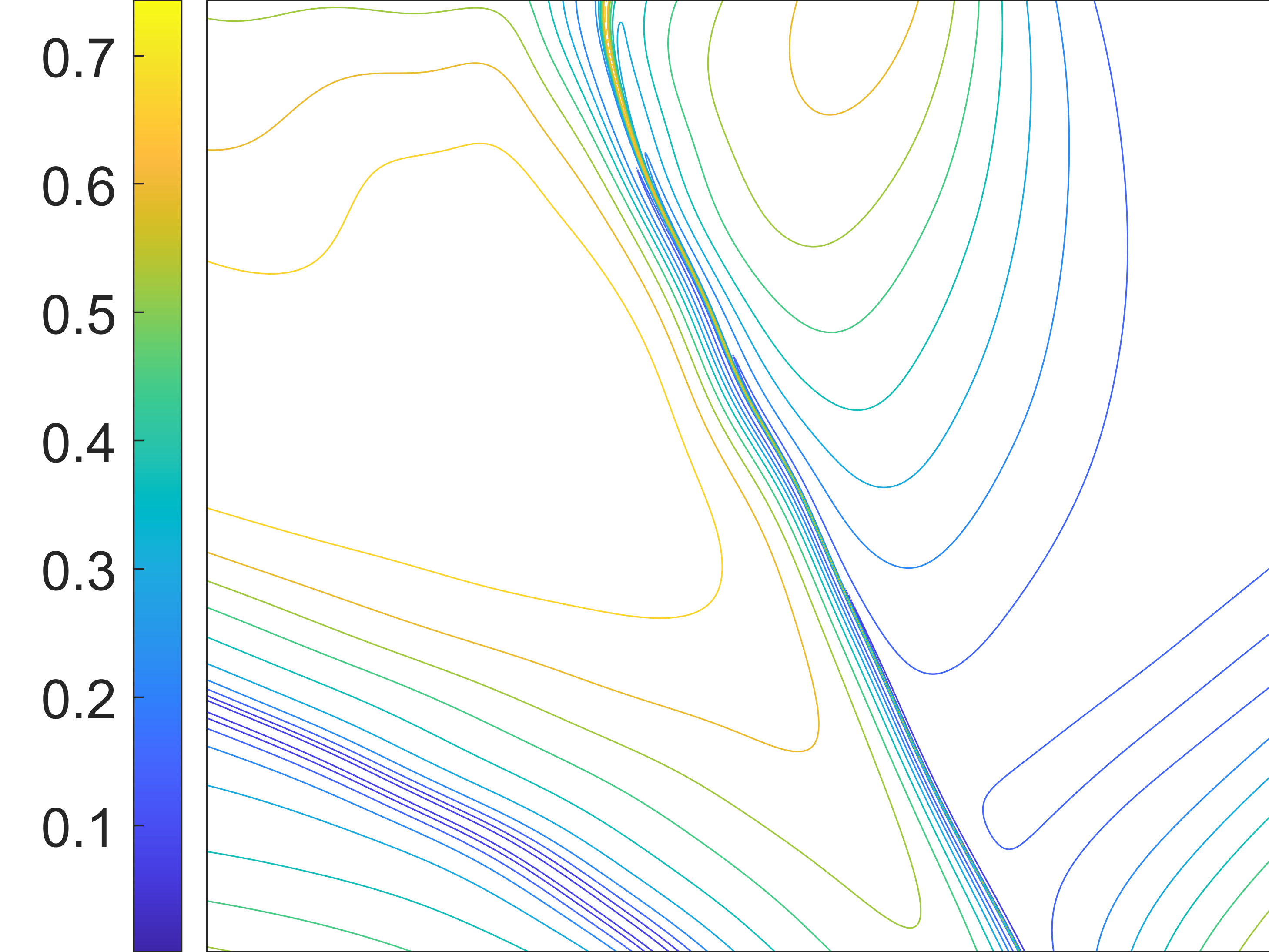}
\includegraphics[width = 0.32 \linewidth]{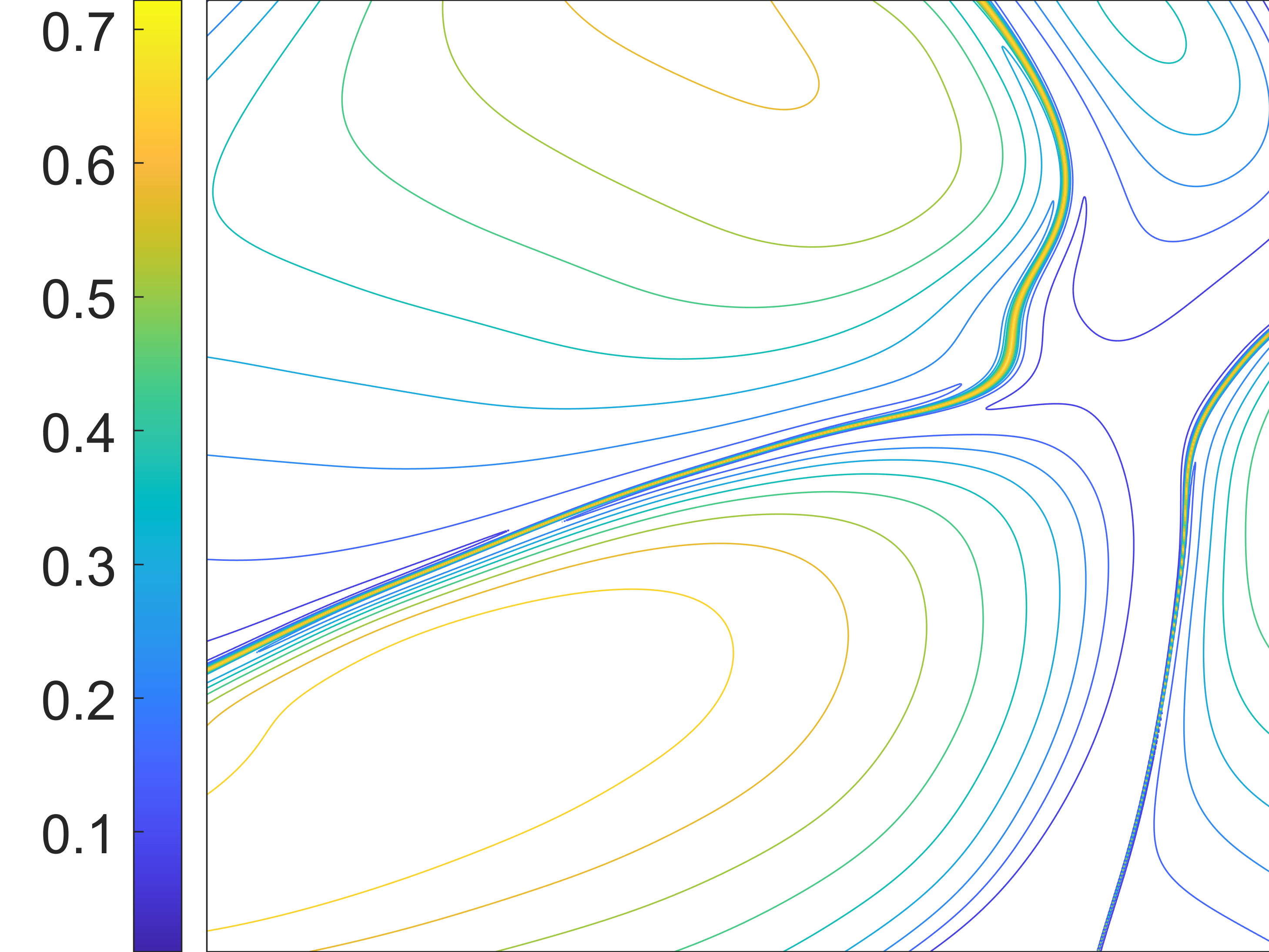}
\caption{Box width $ = \pi$.}
\end{subfigure}
\begin{subfigure}{0.8\linewidth}
\centering
\includegraphics[width = 0.32 \linewidth]{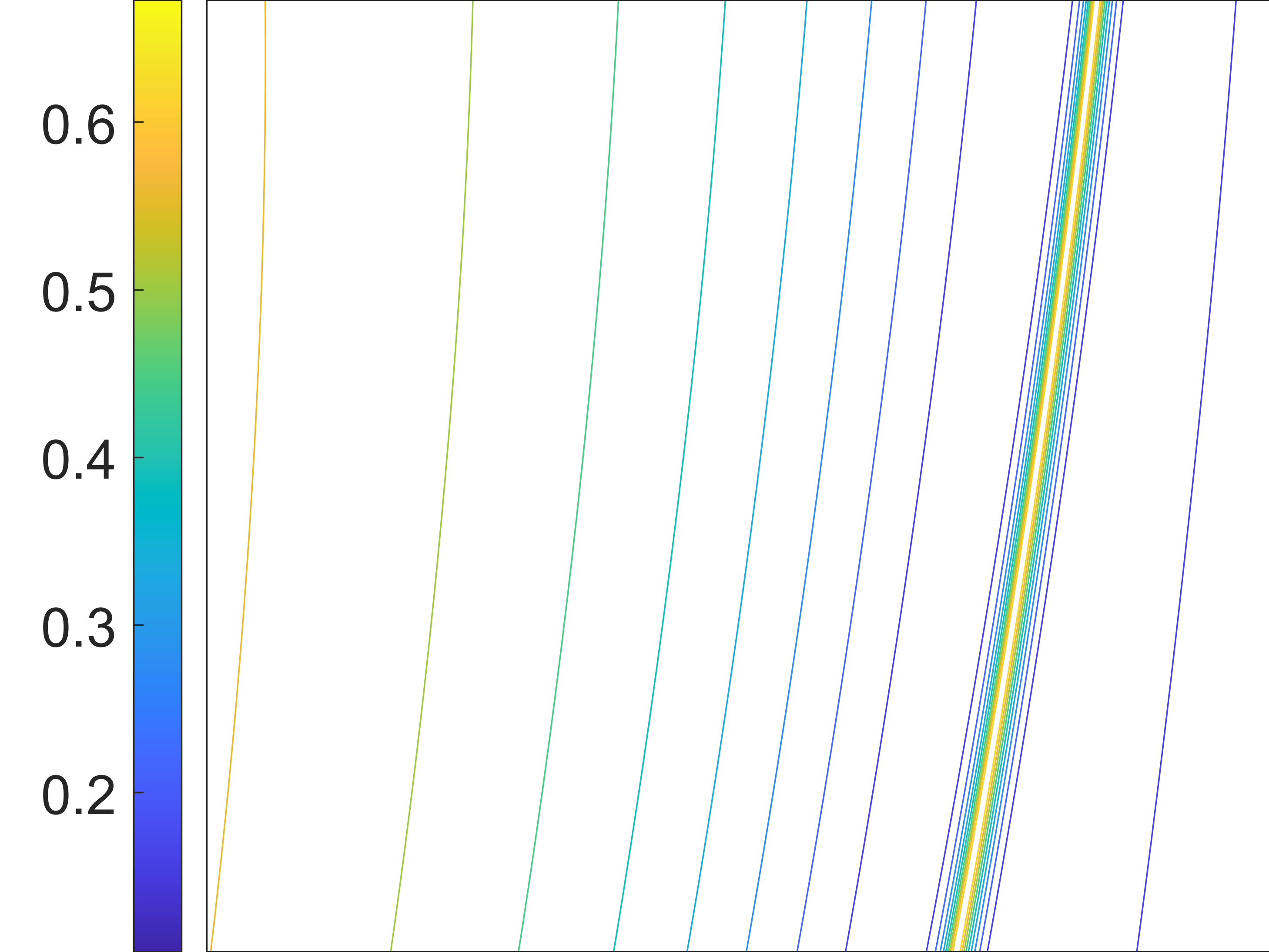}
\includegraphics[width = 0.32 \linewidth]{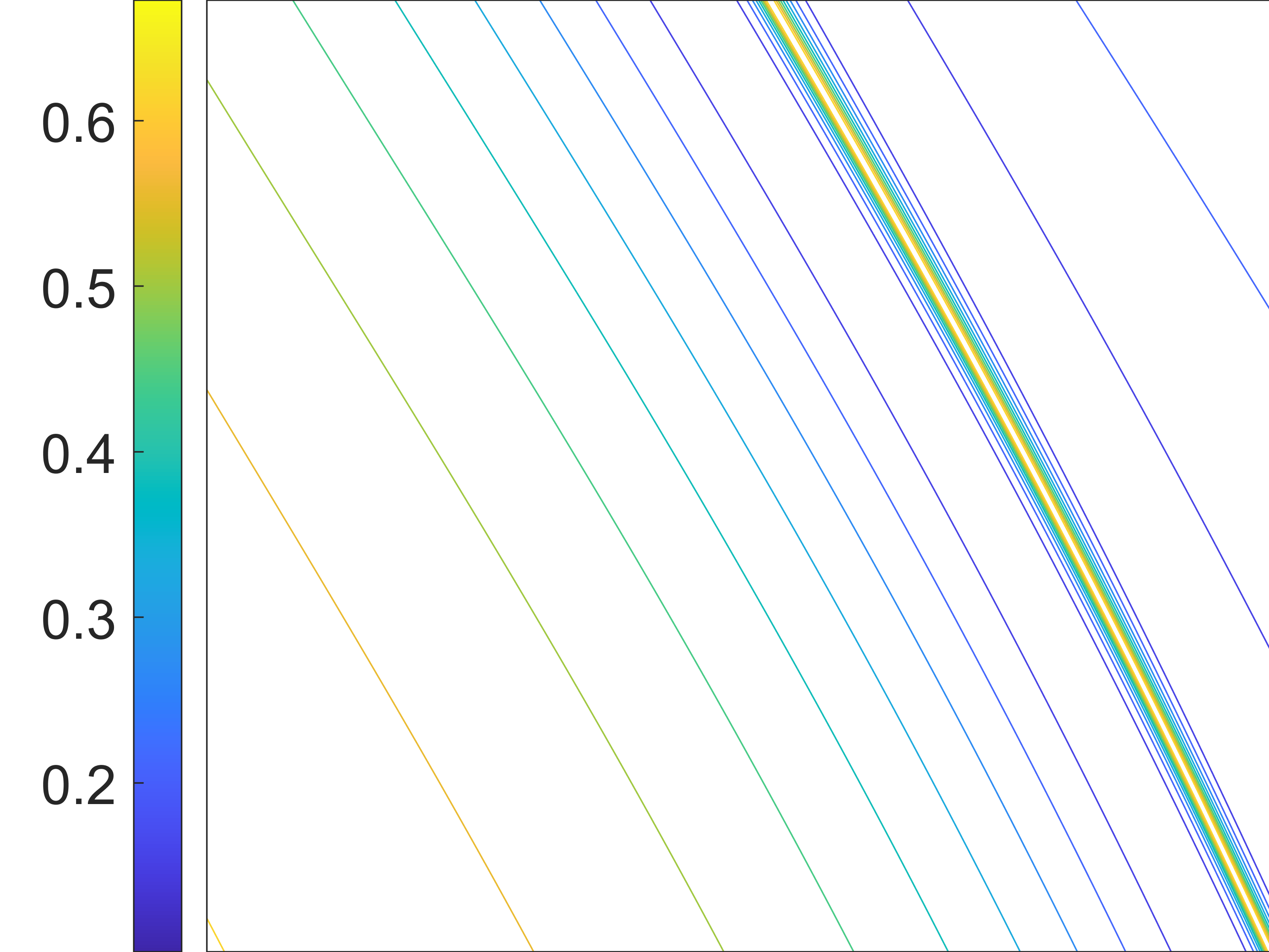}
\includegraphics[width = 0.32 \linewidth]{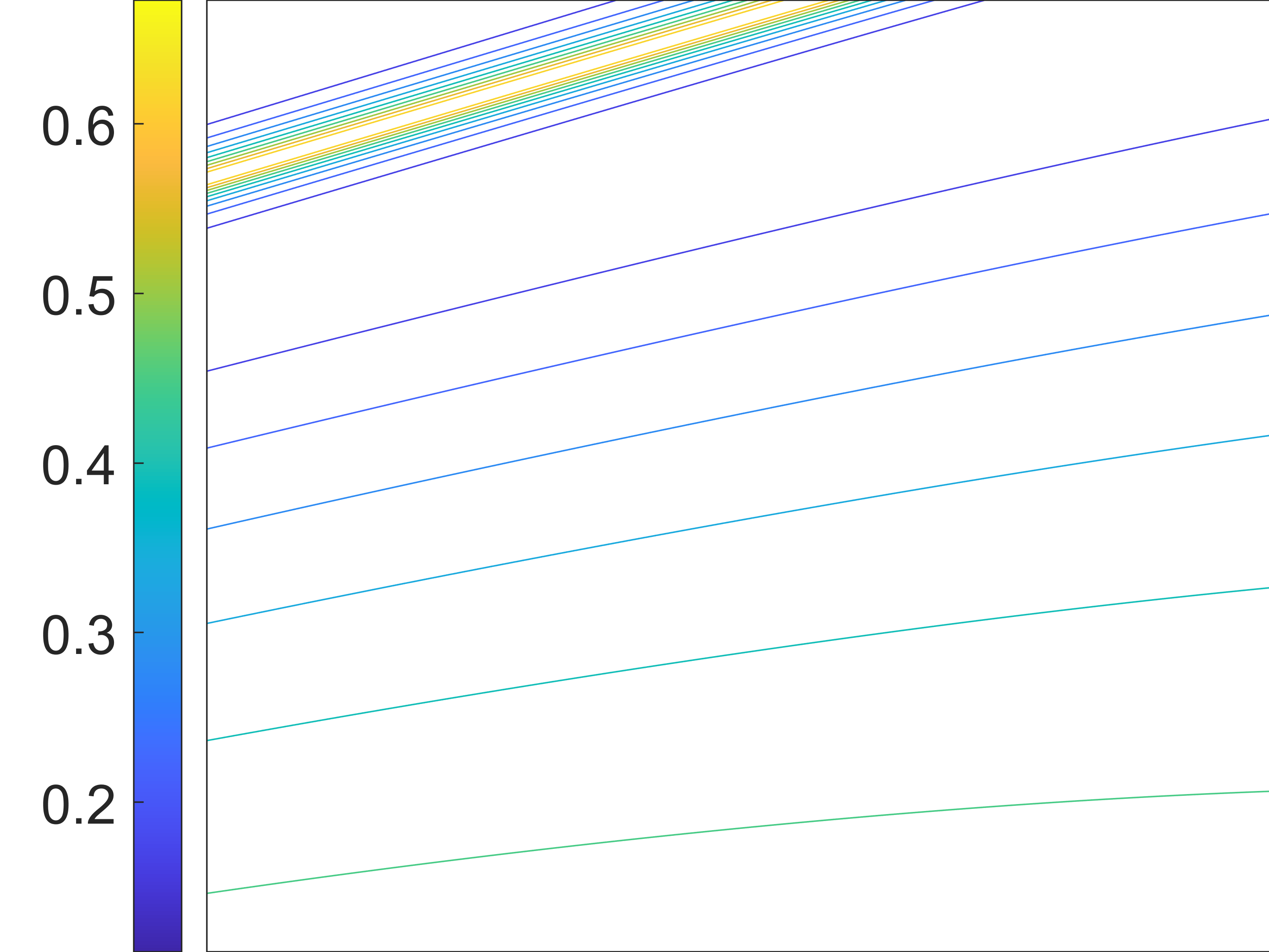}
\caption{Box width $ = \frac{\pi}{16}$.}\label{fig:advCrossZoom2}
\end{subfigure} 
\caption{Contour plot of the transported initial vorticity $\phi$ at time $t=12$, centered at $(x,y,z) = (2.15, 1.5, 1.42)$, 10 isoline values evenly distributed in the range of each frame are shown. From left to right, contours on: $xy$- plane, $xz$-plane and $yz$-plane. The isoline ranges are $[0.0007, 0.7435]$, $[0.0013, 0.7435]$, $[0.0102, 0.7220]$, $[0.1049, 0.6729]$, $[0.0977, 0.6729]$ and $[0.1117, 0.6729]$. Figures are produced using a $1024^2$ 2D grid.}
\label{fig:advCrossZoom}
\end{figure}

The small scale features that develop in the flow after $t=9$ become very fine and make 3D visualization difficult. We present instead three 2D contour plots on the $xy$-, $xz$- and $yz$- planes, respectively. Figures \ref{fig:vtxCrossZoom} and \ref{fig:advCrossZoom} show zoomed views on the lateral cuts of the vorticity magnitude and the advected initial vorticity strength. The highest zoom shows a domain of width $\pi/16$ corresponding to $1/64$ of the computational domain, i.e. smaller than a single cell of the advection grid. In figure \ref{fig:vtxCrossZoom}, we can see a presence of a high vorticity gradient at time $t=12$ and the formation of a vortex sheet. Figure \ref{fig:advCrossZoom} shows more clearly the material deformations which lead to this high gradient. Indeed, the advected initial vorticity shows that two separate level-sets of the vorticity were pushed close together by the flow. The yellow region squeezed between the blue curves in \ref{fig:advCrossZoom2} was formed from the flattening of the initial vortex tubes. However, the absence of a highly concentrated vorticity peak in \ref{fig:vtxCrossZoom} in contrast to \ref{fig:advCrossZoom} suggests that the vorticity direction was not fully aligned with the strain tensor eigenvector of the largest eigenvalue indicating that the vortex streching term attenuated the vorticity gradient in this region. This nicely illustrates that vortex streching can be quantified locally.

\section{Conclusion} \label{sec:conclu}

In this paper, we presented a novel numerical method for solving the 3D incompressible Euler equations. This method is based on the Gradient-Augmented Level-Set and Jet-Scheme frameworks and extends the previous work on the CM method for the 2D Euler equations, studied in \cite{CME}. Taking a more geometric approach, we proposed a reformulation of the CM framework which allows for the vortex stretching term in 3D to be seamlessly incorporated. The approach in this paper can be summarized as follows: we evolve numerically the backward characteristic map, i.e. the backward-in-time flow map generated by the fluid velocity. This map is in fact the transition map between the Eulerian and Lagrangian coordinate charts. From the Kelvin circulation theorem, the time $t$ vorticity field, expressed as a differential 2-form, can be computed from the pullback of the initial vorticity by the characteristic map, thereby allowing for vortex line deformation and stretching. The Biot-Savart law of the vorticity is then computed using Fourier spectral methods to provide the velocity field needed to further evolve the map.

This reformulation was motivated by the fact that direct integration of the vortex stretching term as a generic source term does not preserve the characteristic structure of the method and nullifies some of its numerical advantages. The geometric approach we proposed here allows us to retain all numerical properties of the CM method previously shown in \cite{CME}. Firstly, the functional representation of the vorticity field through pullback by the characteristic map preserves the fine scales generated by the inviscid flow. This is demonstrated by the tests in section \ref{sec:numTests}, notably, the energy and enstrophy spectra plots in section \ref{sec:numTestKerr} show that by oversampling the vorticity on a fine grid, high frequency features of the solution can be accurately reconstructed. Secondly, the rapid growth in the solution gradient can be efficiently resolved using the group structure of the characteristic map; the submap decomposition allows us to achieve multiplicative growth in the spatial resolution using only fast coarse grid map computations. As seen in section \ref{sec:numTestKerr}, the flow of the vortex reconnection test can be evolved up to time $t=17$ using 79 coarse grid maps of size $64\times 48\times 32 $, whereas purely Eulerian pseudospectral methods would require computations on grid sizes of order $1024^3$ to achieve comparable results. Lastly, by defining the vorticity 2-form as the pullback of the initial condition by the map, we preserve the non-dissipative property of the CM method, this was shown theoretically in the error estimates in section \ref{sec:ErrorEst} where we found that the error is in fact advective and that the numerical vorticity is related to the exact vorticity by pullback by an error map. Similar to the Lagrangian-Averaged Euler-$\alpha$ equations and the Kelvin-filtered turbulence models \cite{foias2001navier}, the solution is obtained from a nonlinear dispersive perturbation of the equation. This was also supported by the numerical tests which showed that subgrid structures are preserved and that the spectrum plots of the solutions do not exhibit the exponential decay associated to artificial viscous dissipation; as a matter of fact, higher frequency modes can be reconstructed by finer samplings of the solution.

This work constitutes a first investigation of the CM method for the 3D incompressible Euler equations, we demonstrated here the key properties of the method and studied its numerical accuracy. This opens numerous directions for future research. For instance, the inclusion of forcing terms which conform to the characteristic structure of the method could be of interest. Extensions to the framework to take into account different geometries and domain boundaries is also an important subject of further investigation. Furthermore, the numerical tests in section \ref{sec:numTestKerr} suggest that, with additional extensions in terms of spatial and temporal adaptivity or improvements in the discretization spaces and basis functions, more efficient higher resolution simulations using the CM method could produce new insights into the blow-up question for the Euler equations. We believe that the CM method provides a novel and unique approach for the simulation of inviscid flows, and offers a suitable framework for the numerical study of fluid dynamics.

\section*{Acknowledgements}
The authors acknowledge partial funding from the Agence Nationale de la Recherche (ANR), grant ANR-20-CE46-0010-01. Additionally, JCN acknowledges partial support from the NSERC Discovery Grant program.

% \newpage
\appendix
\section{Vorticity Sampling} \label{append:vtxSmp}
The map error studied in section \ref{sec:ErrorEst} can roughly be split into two main contributing parts: the map representation error arising from the interpolation of the flow generated by $\vtu$ and the velocity representation error $\vtu - \vu$, which correspond to the $\vphi$ and $\vpsi$ errors respectively.

The $\vpsi$ velocity representation error depend in part on the resolution of the $\gV$ grid used for the vorticity sampling. This grid needs to be fine enough to avoid Fourier aliasing from lack of resolution. One possible solution would be to use a dynamic sampling of the vorticity field such as through the use of oct-tree structured meshgrids. Here we will consider another adaptive sampling method: we directly evaluate a mollified version of the vorticity, that is defining $\vtw = \mu_h * \vw^n$ where $\mu_h$ is a mollifier supported in a neighborhood of size $h$, this approach was studied to some extent in \cite{CME}. We pick $h$ to be smaller than the cellwidth of $\gV$, the evaluation of the mollified vorticity at grid points $\vx_\vi$ can then be expressed as the sum of the convolution integrals in all 8 cells $C_{\vi + \vr}$ adjacent to $\vx_\vi$:
\begin{gather} \label{eq:mollIntCell}
(\mu_h * \vw^n)(\vx_\vi) = \sum_{\vr \in \{-1, 0\}^3} \int_{C_{\vi + \vr}} \mu_h(\vx_\vi-\vy) \vw^n(\vy) \d \vy .
\end{gather}
For instance, in the tests in section \ref{sec:numTestCross}, the mollifier $\mu_h$ is chosen to be a 3D tensor of $\cos^2(x/h)$ supported in the cells adjacent to $\vx_\vi$. The integral in each cell is computed using numerical quadrature and the number of quadrature points is then adapted to the local oscillations of $\vw^n$ to ensure the accuracy of the mollification. The resulting algorithm effectively computes a mollified vorticity field where subgrid oscillations are filtered out, by choosing an appropriate mollification scale, the $\gV$ grid can resolve the mollified vorticity without aliasing errors and the pointwise evaluation of $\mu_h * \vw^n$ on $\gV$ is accurate as long as the sampling in each cell is sufficiently dense. This approach is also related to the Kelvin-filtered Euler equations. From the Kelvin circulation theorem, total circulation along a closed curve evolving under the flow is conserved. Both in the Kelvin-filtered equations and the CM method, circulation is conserved for closed curves evolving under a modified flow; in the Kelvin-filtered case, the modified flow arises from a filtered velocity field and in the CM case, the numerical flow map is modified by the $\vphi$ and $\vpsi$ errors, of which $\vpsi$ contains the velocity filtering. In any case, the mollification of the transport velocity still deteriorates the accuracy of the method and, although the vorticity evolves by pulllback, dispersion is introduced in the form of a less energetic transport flow. The proper sampling of the vorticity is still subject of our current work, and future directions may include the use of adaptive meshes and wavelet transforms to reduce sampling errors.

\newpage
%\vspace*{1cm}
\bibliographystyle{siamplain}
\bibliography{CME3D}

\end{document}